\documentstyle[12pt]{amsart}

\begin{document}

\def\Resetstrings{
    \def\present{ }\let\bgroup={\let\egroup=}
    \def\Astr{}\def\astr{}\def\Atest{}\def\atest{}%
    \def\Bstr{}\def\bstr{}\def\Btest{}\def\btest{}%
    \def\Cstr{}\def\cstr{}\def\Ctest{}\def\ctest{}%
    \def\Dstr{}\def\dstr{}\def\Dtest{}\def\dtest{}%
    \def\Estr{}\def\estr{}\def\Etest{}\def\etest{}%
    \def\Fstr{}\def\fstr{}\def\Ftest{}\def\ftest{}%
    \def\Gstr{}\def\gstr{}\def\Gtest{}\def\gtest{}%
    \def\Hstr{}\def\hstr{}\def\Htest{}\def\htest{}%
    \def\Istr{}\def\istr{}\def\Itest{}\def\itest{}%
    \def\Jstr{}\def\jstr{}\def\Jtest{}\def\jtest{}%
    \def\Kstr{}\def\kstr{}\def\Ktest{}\def\ktest{}%
    \def\Lstr{}\def\lstr{}\def\Ltest{}\def\ltest{}%
    \def\Mstr{}\def\mstr{}\def\Mtest{}\def\mtest{}%
    \def\Nstr{}\def\nstr{}\def\Ntest{}\def\ntest{}%
    \def\Ostr{}\def\ostr{}\def\Otest{}\def\otest{}%
    \def\Pstr{}\def\pstr{}\def\Ptest{}\def\ptest{}%
    \def\Qstr{}\def\qstr{}\def\Qtest{}\def\qtest{}%
    \def\Rstr{}\def\rstr{}\def\Rtest{}\def\rtest{}%
    \def\Sstr{}\def\sstr{}\def\Stest{}\def\stest{}%
    \def\Tstr{}\def\tstr{}\def\Ttest{}\def\ttest{}%
    \def\Ustr{}\def\ustr{}\def\Utest{}\def\utest{}%
    \def\Vstr{}\def\vstr{}\def\Vtest{}\def\vtest{}%
    \def\Wstr{}\def\wstr{}\def\Wtest{}\def\wtest{}%
    \def\Xstr{}\def\xstr{}\def\Xtest{}\def\xtest{}%
    \def\Ystr{}\def\ystr{}\def\Ytest{}\def\ytest{}%
}
\Resetstrings\def\Ztest{}\def\ztest{}

\def\Refformat{
         \if\Jtest\present
             {\if\Vtest\present\journalarticleformat
                  \else\conferencereportformat\fi}
            \else\if\Btest\present\bookarticleformat
               \else\if\Rtest\present\technicalreportformat
                  \else\if\Itest\present\bookformat
                     \else\otherformat\fi\fi\fi\fi}

\def\Rpunct{
   \def\Lspace{ }%
   \def\Lperiod{ }
   \def\Lcomma{ }
   \def\Lquest{ }
   \def\Lcolon{ }
   \def\Lscolon{ }
   \def\Lbang{ }
   \def\Lquote{ }
   \def\Lqquote{ }
   \def\Lrquote{ }
   \def\Rspace{}%
   \def\Rperiod{.}
   \def\Rcomma{,}
   \def\Rquest{?}
   \def\Rcolon{:}
   \def\Rscolon{;}
   \def\Rbang{!}
   \def\Rquote{'}
   \def\Rqquote{"}
   \def\Rrquote{`}
   }

\def\Lpunct{
   \def\Lspace{}%
   \def\Lperiod{\unskip.}
   \def\Lcomma{\unskip,}
   \def\Lquest{\unskip?}
   \def\Lcolon{\unskip:}
   \def\Lscolon{\unskip;}
   \def\Lbang{\unskip!}
   \def\Lquote{\unskip'}
   \def\Lqquote{\unskip"}
   \def\Lrquote{\unskip`}
   \def\Rspace{\spacefactor=1000}%
   \def\Rperiod{\spacefactor=3000}
   \def\Rcomma{\spacefactor=1250}
   \def\Rquest{\spacefactor=3000}
   \def\Rcolon{\spacefactor=2000}
   \def\Rscolon{\spacefactor=1250}
   \def\Rbang{\spacefactor=3000}
   \def\Rquote{\spacefactor=1000}
   \def\Rqquote{\spacefactor=1000}
   \def\Rrquote{\spacefactor=1000}
   }

\def\Refstd{
     \def\Acomma{\unskip, }
     \def\Aand{\unskip\ and }
     \def\Aandd{\unskip\ and }
     \def\Ecomma{\unskip, }
     \def\Eand{\unskip\ and }
     \def\Eandd{\unskip\ and }
     \def\acomma{\unskip, }
     \def\aand{\unskip\ and }
     \def\aandd{\unskip\ and }
     \def\ecomma{\unskip, }
     \def\eand{\unskip\ and }
     \def\eandd{\unskip\ and }
     \def\Namecomma{\unskip, }
     \def\Nameand{\unskip\ and }
     \def\Nameandd{\unskip\ and }
     \def\Revcomma{\unskip, }
     \def\Initper{.\ }
     \def\Initgap{\dimen0=\spaceskip\divide\dimen0 by 2\hskip-\dimen0}%
   }

\def\Smallcapsaand{
     \def\Aand{\unskip\bgroup{\Smallcapsfont\ AND }\egroup}%
     \def\Aandd{\unskip\bgroup{\Smallcapsfont\ AND }\egroup}%
     \def\eand{\unskip\bgroup\Smallcapsfont\ AND \egroup}%
     \def\eandd{\unskip\bgroup\Smallcapsfont\ AND \egroup}%
   }

\def\Smallcapseand{
     \def\Eand{\unskip\bgroup\Smallcapsfont\ AND \egroup}%
     \def\Eandd{\unskip\bgroup\Smallcapsfont\ AND \egroup}%
     \def\aand{\unskip\bgroup\Smallcapsfont\ AND \egroup}%
     \def\aandd{\unskip\bgroup\Smallcapsfont\ AND \egroup}%
   }

\def\Refstda{
    \chardef\Ampersand=`\&
    \def\Acomma{\unskip, }
    \def\Aand{\unskip\ \Ampersand\ }
    \def\Aandd{\unskip\ \Ampersand\ }
    \def\Ecomma{\unskip, }
    \def\Eand{\unskip\ \Ampersand\ }
    \def\Eandd{\unskip\ \Ampersand\ }
    \def\acomma{\unskip, }
    \def\aand{\unskip\ \Ampersand\ }
    \def\aandd{\unskip\ \Ampersand\ }
    \def\ecomma{\unskip, }
    \def\eand{\unskip\ \Ampersand\ }
    \def\eandd{\unskip\ \Ampersand\ }
    \def\Namecomma{\unskip, }
    \def\Nameand{\unskip\ \Ampersand\ }
   \def\Nameandd{\unskip\ \Ampersand\ }
    \def\Revcomma{\unskip, }
    \def\Initper{.\ }
    \def\Initgap{\dimen0=\spaceskip\divide\dimen0 by 2\hskip-\dimen0}%
  }

   \def\Citefont{}
   \def\ACitefont{}
   \def\Authfont{}
   \def\Titlefont{}
   \def\Tomefont{\sl}
   \def\Volfont{}
   \def\Flagfont{}
   \def\Reffont{\rm}
   \def\Smallcapsfont{\sevenrm}
   \def\Flagstyle#1{\hangindent\parindent\indent\hbox to0pt
       {\hss[{\Flagfont#1}]\kern.5em}\ignorespaces}

\def\Underlinemark{\vrule height .7pt depth 0pt width 3pc}

\def\Citebrackets{\Rpunct
   \def\Lcitemark{\def\Cfont{\Citefont}[\bgroup\Cfont}
   \def\Rcitemark{\egroup]}
   \def\LAcitemark{\def\Cfont{\ACitefont}\bgroup\ACitefont}%
   \def\RAcitemark{\egroup}
   \def\LIcitemark{\egroup}
   \def\RIcitemark{\bgroup\Cfont}
   \def\Citehyphen{\egroup--\bgroup\Cfont}
   \def\Citecomma{\egroup,\hskip0pt\bgroup\Cfont}%
   \def\Citebreak{}
   }

\def\Citeparen{\Rpunct
   \def\Lcitemark{\def\Cfont{\Citefont}(\bgroup\Cfont}
   \def\Rcitemark{\egroup)}
   \def\LAcitemark{\def\Cfont{\ACitefont}\bgroup\ACitefont}%
   \def\RAcitemark{\egroup}
   \def\LIcitemark{\egroup}
   \def\RIcitemark{\bgroup\Cfont}
   \def\Citehyphen{\egroup--\bgroup\Cfont}
   \def\Citecomma{\egroup,\hskip0pt\bgroup\Cfont}%
   \def\Citebreak{}
   }

\def\Citesuper{\Lpunct
   \def\Lcitemark{\def\Cfont{\Citefont}\raise1ex\hbox\bgroup\bgroup\Cfont}%
   \def\Rcitemark{\egroup\egroup}
   \def\LAcitemark{\def\Cfont{\ACitefont}\bgroup\ACitefont}%
   \def\RAcitemark{\egroup}
   \def\LIcitemark{\egroup\egroup}
   \def\RIcitemark{\raise1ex\hbox\bgroup\bgroup\Cfont}%
   \def\Citehyphen{\egroup--\bgroup\Cfont}
   \def\Citecomma{\egroup,\hskip0pt\bgroup%
      \Cfont}
   \def\Citebreak{}
   } 

\def\Citenamedate{\Rpunct
   \def\Lcitemark{
      \def\Citebreak{\egroup\ [\bgroup\Citefont}
      \def\Citecomma{\egroup]; 
         \bgroup\let\uchyph=1\Citefont}(\bgroup\let\uchyph=1\Citefont}%
   \def\Rcitemark{\egroup])}
   \def\LAcitemark{
      \def\Citebreak{\egroup\ [\bgroup\Citefont}\def\Citecomma{\egroup], %
         \bgroup\ACitefont }\bgroup\let\uchyph=1\ACitefont}%
   \def\RAcitemark{\egroup]}
  \def\Citehyphen{\egroup--\bgroup\Citefont}
   \def\LIcitemark{\egroup}
   \def\RIcitemark{\bgroup\Citefont}
   }

\Refstd\Citebrackets
\def\Citefont{\bf}\def\Titlefont{\sl}\def\Volfont{\bf}\def\Tomefont{\Reffont}

\def\journalarticleformat{\Reffont\let\uchyph=1\parindent=1.25pc\def\Comma{}%
                \sfcode`\.=1000\sfcode`\?=1000\sfcode`\!=1000\sfcode`\:=1000\sfcode`\;=1000\sfcode`\,=1000
                \par\vfil\penalty-200\vfilneg
      \if\Ftest\present\Flagstyle\Fstr\fi%
       \if\Atest\present\bgroup\Authfont\Astr\egroup\def\Comma{\unskip, }\fi%
        \if\Ttest\present\Comma\bgroup\Titlefont\Tstr\egroup\def\Comma{, }\fi%
         \if\etest\present\hskip.2em(\bgroup\estr\egroup)\def\Comma{\unskip, }\fi%
          \if\Jtest\present\Comma\bgroup\Tomefont\Jstr\/\egroup\def\Comma{, }\fi%
           \if\Vtest\present\if\Jtest\present\hskip.2em\else\Comma\fi\bgroup\Volfont\Vstr\egroup\def\Comma{, }\fi%
            \if\Dtest\present\hskip.2em(\bgroup\Dstr\egroup)\def\Comma{, }\fi%
             \if\Ptest\present\bgroup, \Pstr\egroup\def\Comma{, }\fi%
              \if\ttest\present\Comma\bgroup\Titlefont\tstr\egroup\def\Comma{, }\fi%
               \if\jtest\present\Comma\bgroup\Tomefont\jstr\/\egroup\def\Comma{, }\fi%
                \if\vtest,\present\if\jtest\present\hskip.2em\else\Comma\fi\bgroup\Volfont\vstr\egroup\def\Comma{, }\fi%
                 \if\dtest\present\hskip.2em(\bgroup\dstr\egroup)\def\Comma{, }\fi%
                  \if\ptest\present\bgroup, \pstr\egroup\def\Comma{, }\fi%
                   \if\Gtest\present{\Comma Gov't ordering no. }\bgroup\Gstr\egroup\def\Comma{, }\fi%
                    \if\Mtest\present\Comma MR \#\bgroup\Mstr\egroup\def\Comma{, }\fi%
                     \if\Otest\present{\Comma\bgroup\Ostr\egroup.}\else{.}\fi%
                      \vskip3ptplus1ptminus1pt}

\def\conferencereportformat{\Reffont\let\uchyph=1\parindent=1.25pc\def\Comma{}%
                \sfcode`\.=1000\sfcode`\?=1000\sfcode`\!=1000\sfcode`\:=1000\sfcode`\;=1000\sfcode`\,=1000
                \par\vfil\penalty-200\vfilneg
      \if\Ftest\present\Flagstyle\Fstr\fi%
       \if\Atest\present\bgroup\Authfont\Astr\egroup\def\Comma{\unskip, }\fi%
        \if\Ttest\present\Comma\bgroup\Titlefont\Tstr\egroup\def\Comma{, }\fi%
         \if\Jtest\present\Comma\bgroup\Tomefont\Jstr\/\egroup\def\Comma{, }\fi%
          \if\Ctest\present\Comma\bgroup\Cstr\egroup\def\Comma{, }\fi%
           \if\Dtest\present\hskip.2em(\bgroup\Dstr\egroup)\def\Comma{, }\fi%
            \if\Mtest\present\Comma MR \#\bgroup\Mstr\egroup\def\Comma{, }\fi%
             \if\Otest\present{\Comma\bgroup\Ostr\egroup.}\else{.}\fi%
              \vskip3ptplus1ptminus1pt}

\def\bookarticleformat{\Reffont\let\uchyph=1\parindent=1.25pc\def\Comma{}%
                \sfcode`\.=1000\sfcode`\?=1000\sfcode`\!=1000\sfcode`\:=1000\sfcode`\;=1000\sfcode`\,=1000
                \par\vfil\penalty-200\vfilneg
      \if\Ftest\present\Flagstyle\Fstr\fi%
       \if\Atest\present\bgroup\Authfont\Astr\egroup\def\Comma{\unskip, }\fi%
        \if\Ttest\present\Comma\bgroup\Titlefont\Tstr\egroup\def\Comma{, }\fi%
         \if\etest\present\hskip.2em(\bgroup\estr\egroup)\def\Comma{\unskip, }\fi%
          \if\Btest\present\Comma in \bgroup\Tomefont\Bstr\/\egroup\def\Comma{\unskip, }\fi%
           \if\otest\present\ \bgroup\ostr\egroup\def\Comma{, }\fi%
            \if\Etest\present\Comma\bgroup\Estr\egroup\unskip, \ifnum\Ecnt>1eds.\else ed.\fi\def\Comma{, }\fi%
             \if\Stest\present\Comma\bgroup\Sstr\egroup\def\Comma{, }\fi%
              \if\Vtest\present\Comma vol. \bgroup\Vstr\egroup\def\Comma{, }\fi%
               \if\Ntest\present\Comma no. \bgroup\Nstr\egroup\def\Comma{, }\fi%
                \if\Itest\present\Comma\bgroup\Istr\egroup\def\Comma{, }\fi%
                 \if\Ctest\present\Comma\bgroup\Cstr\egroup\def\Comma{, }\fi%
                  \if\Dtest\present\Comma\bgroup\Dstr\egroup\def\Comma{, }\fi%
                   \if\Ptest\present\Comma\Pstr\def\Comma{, }\fi%
                    \if\ttest\present\Comma\bgroup\Titlefont\Tstr\egroup\def\Comma{, }\fi%
                     \if\btest\present\Comma in \bgroup\Tomefont\bstr\egroup\def\Comma{, }\fi%
                       \if\atest\present\Comma\bgroup\astr\egroup\unskip, \if\acnt\present eds.\else ed.\fi\def\Comma{, }\fi%
                        \if\stest\present\Comma\bgroup\sstr\egroup\def\Comma{, }\fi%
                         \if\vtest\present\Comma vol. \bgroup\vstr\egroup\def\Comma{, }\fi%
                          \if\ntest\present\Comma no. \bgroup\nstr\egroup\def\Comma{, }\fi%
                           \if\itest\present\Comma\bgroup\istr\egroup\def\Comma{, }\fi%
                            \if\ctest\present\Comma\bgroup\cstr\egroup\def\Comma{, }\fi%
                             \if\dtest\present\Comma\bgroup\dstr\egroup\def\Comma{, }\fi%
                              \if\ptest\present\Comma\pstr\def\Comma{, }\fi%
                               \if\Gtest\present{\Comma Gov't ordering no. }\bgroup\Gstr\egroup\def\Comma{, }\fi%
                                \if\Mtest\present\Comma MR \#\bgroup\Mstr\egroup\def\Comma{, }\fi%
                                 \if\Otest\present{\Comma\bgroup\Ostr\egroup.}\else{.}\fi%
                                  \vskip3ptplus1ptminus1pt}

\def\bookformat{\Reffont\let\uchyph=1\parindent=1.25pc\def\Comma{}%
                \sfcode`\.=1000\sfcode`\?=1000\sfcode`\!=1000\sfcode`\:=1000\sfcode`\;=1000\sfcode`\,=1000
                \par\vfil\penalty-200\vfilneg
      \if\Ftest\present\Flagstyle\Fstr\fi%
       \if\Atest\present\bgroup\Authfont\Astr\egroup\def\Comma{\unskip, }%
            \else\if\Etest\present\bgroup\def\Eand{\Aand}\def\Eandd{\Aandd}\Authfont\Estr\egroup\unskip, \ifnum\Ecnt>1eds.\else ed.\fi\def\Comma{, }%
                  \else\if\Itest\present\bgroup\Authfont\Istr\egroup\def\Comma{, }\fi\fi\fi%
          \if\Ttest\present\Comma\bgroup\Titlefont\Tstr\/\egroup\def\Comma{\unskip, }%
                \else\if\Btest\present\Comma\bgroup\Titlefont\Bstr\/\egroup\def\Comma{\unskip, }\fi\fi%
            \if\otest\present\ \bgroup\ostr\egroup\def\Comma{, }\fi%
             \if\etest\present\hskip.2em(\bgroup\estr\egroup)\def\Comma{\unskip, }\fi%
              \if\Stest\present\Comma\bgroup\Sstr\egroup\def\Comma{, }\fi%
               \if\Vtest\present\Comma vol. \bgroup\Vstr\egroup\def\Comma{, }\fi%
                \if\Ntest\present\Comma no. \bgroup\Nstr\egroup\def\Comma{, }\fi%
                 \if\Atest\present\if\Itest\present
                         \Comma\bgroup\Istr\egroup\def\Comma{\unskip, }\fi%
                      \else\if\Etest\present\if\Itest\present
                              \Comma\bgroup\Istr\egroup\def\Comma{\unskip, }\fi\fi\fi%
                     \if\Ctest\present\Comma\bgroup\Cstr\egroup\def\Comma{, }\fi%
                      \if\Dtest\present\Comma\bgroup\Dstr\egroup\def\Comma{, }\fi%
                       \if\ttest\present\Comma\bgroup\Titlefont\tstr\egroup\def\Comma{, }%
                             \else\if\btest\present\Comma\bgroup\Titlefont\bstr\egroup\def\Comma{, }\fi\fi%
                          \if\stest\present\Comma\bgroup\sstr\egroup\def\Comma{, }\fi%
                           \if\vtest\present\Comma vol. \bgroup\vstr\egroup\def\Comma{, }\fi%
                            \if\ntest\present\Comma no. \bgroup\nstr\egroup\def\Comma{, }\fi%
                             \if\itest\present\Comma\bgroup\istr\egroup\def\Comma{, }\fi%
                              \if\ctest\present\Comma\bgroup\cstr\egroup\def\Comma{, }\fi%
                               \if\dtest\present\Comma\bgroup\dstr\egroup\def\Comma{, }\fi%
                                \if\Gtest\present{\Comma Gov't ordering no. }\bgroup\Gstr\egroup\def\Comma{, }\fi%
                                 \if\Mtest\present\Comma MR \#\bgroup\Mstr\egroup\def\Comma{, }\fi%
                                  \if\Otest\present{\Comma\bgroup\Ostr\egroup.}\else{.}\fi%
                                   \vskip3ptplus1ptminus1pt}

\def\technicalreportformat{\Reffont\let\uchyph=1\parindent=1.25pc\def\Comma{}%
                \sfcode`\.=1000\sfcode`\?=1000\sfcode`\!=1000\sfcode`\:=1000\sfcode`\;=1000\sfcode`\,=1000
                \par\vfil\penalty-200\vfilneg
      \if\Ftest\present\Flagstyle\Fstr\fi%
       \if\Atest\present\bgroup\Authfont\Astr\egroup\def\Comma{\unskip, }%
            \else\if\Etest\present\bgroup\def\Eand{\Aand}\def\Eandd{\Aandd}\Authfont\Estr\egroup\unskip, \ifnum\Ecnt>1eds.\else ed.\fi\def\Comma{, }%
                  \else\if\Itest\present\bgroup\Authfont\Istr\egroup\def\Comma{, }\fi\fi\fi%
          \if\Ttest\present\Comma\bgroup\Titlefont\Tstr\egroup\def\Comma{, }\fi%
           \if\Atest\present\if\Itest\present
                   \Comma\bgroup\Istr\egroup\def\Comma{, }\fi%
                \else\if\Etest\present\if\Itest\present
                        \Comma\bgroup\Istr\egroup\def\Comma{, }\fi\fi\fi%
            \if\Rtest\present\Comma\bgroup\Rstr\egroup\def\Comma{, }\fi%
             \if\Ctest\present\Comma\bgroup\Cstr\egroup\def\Comma{, }\fi%
              \if\Dtest\present\Comma\bgroup\Dstr\egroup\def\Comma{, }\fi%
               \if\ttest\present\Comma\bgroup\Titlefont\tstr\egroup\def\Comma{, }\fi%
                \if\itest\present\Comma\bgroup\istr\egroup\def\Comma{, }\fi%
                 \if\rtest\present\Comma\bgroup\rstr\egroup\def\Comma{, }\fi%
                  \if\ctest\present\Comma\bgroup\cstr\egroup\def\Comma{, }\fi%
                   \if\dtest\present\Comma\bgroup\dstr\egroup\def\Comma{, }\fi%
                    \if\Gtest\present{\Comma Gov't ordering no. }\bgroup\Gstr\egroup\def\Comma{, }\fi%
                     \if\Mtest\present\Comma MR \#\bgroup\Mstr\egroup\def\Comma{, }\fi%
                      \if\Otest\present{\Comma\bgroup\Ostr\egroup.}\else{.}\fi%
                       \vskip3ptplus1ptminus1pt}

\def\otherformat{\Reffont\let\uchyph=1\parindent=1.25pc\def\Comma{}%
                \sfcode`\.=1000\sfcode`\?=1000\sfcode`\!=1000\sfcode`\:=1000\sfcode`\;=1000\sfcode`\,=1000
                \par\vfil\penalty-200\vfilneg
      \if\Ftest\present\Flagstyle\Fstr\fi%
       \if\Atest\present\bgroup\Authfont\Astr\egroup\def\Comma{\unskip, }%
            \else\if\Etest\present\bgroup\def\Eand{\Aand}\def\Eandd{\Aandd}\Authfont\Estr\egroup\unskip, \ifnum\Ecnt>1eds.\else ed.\fi\def\Comma{, }%
                  \else\if\Itest\present\bgroup\Authfont\Istr\egroup\def\Comma{, }\fi\fi\fi%
          \if\Ttest\present\Comma\bgroup\Titlefont\Tstr\egroup\def\Comma{, }\fi%
            \if\Atest\present\if\Itest\present
                    \Comma\bgroup\Istr\egroup\def\Comma{, }\fi%
                 \else\if\Etest\present\if\Itest\present
                         \Comma\bgroup\Istr\egroup\def\Comma{, }\fi\fi\fi%
                 \if\Ctest\present\Comma\bgroup\Cstr\egroup\def\Comma{, }\fi%
                  \if\Dtest\present\Comma\bgroup\Dstr\egroup\def\Comma{, }\fi%
                   \if\Gtest\present{\Comma Gov't ordering no. }\bgroup\Gstr\egroup\def\Comma{, }\fi%
                    \if\Mtest\present\Comma MR \#\bgroup\Mstr\egroup\def\Comma{, }\fi%
                     \if\Otest\present{\Comma\bgroup\Ostr\egroup.}\else{.}\fi%
                      \vskip3ptplus1ptminus1pt}

\def\Flagstyle#1{\Flagfont#1. }

\def\frak{\underline}

\def\cases#1{\left\{\,\vcenter{\normalbaselines
    \ialign{$##\hfil$&\quad##\hfil\crcr#1\crcr}}\right.}
\def\pmatrix#1{\left(\matrix{#1}\right)}
\def\matrix#1{\null\,\vcenter{\normalbaselines
    \ialign{\hfil$##$\hfil&&\quad\hfil$##$\hfil\crcr
      \mathstrut\crcr\noalign{\kern-\baselineskip}
      #1\crcr\mathstrut\crcr\noalign{\kern-\baselineskip}}}\,}

\def\Fr{\hbox{Fr}}
\def\Cl{\hbox{Cl}}
\def\Sing{\hbox{Sing}}
\def\Nonsing{\hbox{Nonsing}}
\def\set\{#1|#2\}{\{\,#1\mid #2\,\}}
\def\cite#1{{\bf [#1]}}

\def\proclaim #1. #2\endproclaim{
	\def\procpar{\par}
	\medbreak
	\noindent{\bf#1.\enspace}{\sl#2}\par\medbreak}
\def\proof{\medskip\bf\noindent Proof: \rm}
\def\qed{{\hskip 0pt plus 1filll}
 \vbox{\hrule height 4pt width 4pt}\hfil \par \medbreak}
\def\nusa{{ \nu _{\alpha }}}
\def\iv*{{i{\frak v}^*}}
 \def\it*{{i{\frak t}^*}}    
\def\a*{{{\frak a }^*}}
\def\calD{{{\cal D}}}
\def\calH{{{\cal H}}}
\def\calC{{{\cal C}}}

\def\dbyd{{\partial / \partial }}
\def\fraka{{{\frak a}}}

\def\suba{{_{\alpha }}}
\def\uofg{{ {\Cal U}({\frak g}_{\bf C})}}
\def\bsc{\hbox{\rm :\kern-2pt
\vbox{\hrule height .25pt width .2pt depth .85pt}
\kern-3.4pt
 \vbox{\hrule height -.86pt width .2pt
    depth 1.3pt} 
\kern-3.4pt
\vbox{\hrule height -1.31pt width .2pt
		depth 1.7pt}
\kern-3.4pt
\vbox{\hrule height -1.71pt width .2pt
		depth 2pt}
\kern-3.4pt
\vbox{\hrule height -2pt width .2pt
		depth 2.15pt}
 }}

\def\calA{{ {\cal A}}}
\def\calH {{ {\cal H} }}
\def\calE {{ {\cal E} }}
\def\calO {{ {\cal O} }}
\def\calT {{ {\cal T} }}
\def\chk{{ ^{\vee } }}
\def\barN{{ {\overline N} }}
\def\barn{{ {\overline n} }}
\def\barP{{ {\overline P} }}
\def\fraka{{ {\frak a} }}
\def\bfC{{ {\bf C} }}
\def\bfR{{ {\bf R} }}
\def\calR{{ {\cal R}}}
\def\barw{{ {\overline w} }}
\def\barx{{ {\overline x} }}
\def\barchi{{ {\overline {\chi }} }}
\def\bfZ{{ {\bf Z} }}
\def\s{\sigma}
\def\Ind{\operatorname{Ind}}

\title[$\!$Admissible  representations of non-connected reductive $p$-adic groups$\!$]
{Some results on the admissible  representations of
non-connected reductive $p$-adic groups}

\author{David Goldberg}
\address{\hskip-\parindent David Goldberg,
Department of Mathematics, Purdue University, West Lafayette, IN 47907}
\email{goldberg@@math.purdue.edu}

\author{Rebecca A. Herb}
\address{\hskip-\parindent Rebecca A. Herb,
Department of Mathematics, University of Maryland, College Park, MD 20742}
\email{rah@@math.umd.edu}

\thanks{Goldberg was partially supported by National Science
Foundation Fellowship DMS-9206246 and National Science Foundation
Career Grant DMS-9501868.  Research at MSRI supported in part by
National Science Foundation Grant DMS-9022140.
Herb was partially supported by National Science Foundation
Grant DMS-9400797.}

\begin{abstract}  
 We examine the theory of induced representations
for non-connected reductive $p$-adic groups for which $G/G^0$
is abelian.  We first examine the structure of those representations
of the form $\Ind_{P^0}^G(\sigma),$ where $P^0$ is a 
parabolic subgroup of $G^0$ and $\s$ is a discrete series representation
of the Levi component of $P^0.$  Here we develop a theory of $R$--groups,
extending the theory in the connected case.
We then prove some general results
in the theory of representations of non-connected $p$-adic groups
whose component group is abelian.
We  define
the notion of cuspidal parabolic for $G$ in order to give
a context for this discussion.  
Intertwining operators for the non-connected case are examined
and the notions of supercuspidal and discrete series are defined.
Finally, we examine
parabolic induction from a cuspidal parabolic subgroup of $G.$  Here we
also develop a theory of $R$--groups, and show that these groups parameterize
the induced representations in a manner that is consistent with
the connected case and with the first set of results as well.
\end{abstract}

\maketitle

\section {Introduction}

The theory of induced representations plays a fundamental role within representation theory
in general.  Within the theory of admissible representations of connected reductive
algebraic groups over local fields, parabolic induction is used to complete classification
theories, once certain families of representations are understood 
\Lcitemark 3\Citecomma
8\Citecomma
9\Citecomma
10\Rcitemark \Rspace{}.  The theory of
admissible representations on non-connected reductive groups over nonarchimedean local
fields has been addressed in part in\Lspace \Lcitemark 1\Citecomma
4\Citecomma
6\Citecomma
11\Rcitemark \Rspace{}, among other places.  We will study certain aspects
of parabolic induction for disconnected groups whose component group
is abelian.

Let $F$ be a locally compact, non-discrete, nonarchimedean field of
characteristic zero.  Let $G$ be a (not necessarily connected)
reductive $F$-group.  Thus $G$ is the set of $F$-rational points of a
reductive algebraic group defined over
$F$. Let
$G^0$ be the connected component of the identity in $G$.  We assume that $G/G^0$ is finite
and abelian.

Our goal is to address three major points.  The first is an extension of 
the results of\Lspace \Lcitemark 6\Rcitemark \Rspace{} to the case at hand.  
This entails a study
of induction from a parabolic subgroup of $G^\circ$
to $G.$  In particular, suppose that $P^\circ=M^\circ N$ 
is a parabolic subgroup of $G^\circ,$ and let $\s_0$
be an irreducible discrete series representation of $M^\circ.$
We are interested in the structure of $\pi_0=\Ind_{P^\circ}^G(\s_0).$
In\Lspace \Lcitemark 1\Rcitemark \Rspace{} Arthur suggests a
construction,
in terms of the
conjectural local Langlands parameterization, of a finite group whose
representation theory should describe the structure of $\pi_0,$
when $G/G^\circ$ is cyclic.  In\Lspace \Lcitemark 6\Rcitemark \Rspace{} the
case where $G/G^\circ$ is of prime order is studied, and
there is a construction, on the group side,
of a finite group $R_G(\s_0)$
which (along with an appropriate $2$--cocycle)
parameterizes the components of $\pi_0.$  It is also shown there
that $R_G(\s_0)$ must be isomorphic to Arthur's group $R_{\psi,\s_0},$
if the latter exists.  One cannot confirm the existence
of $R_{\psi,\s_0}$ without  proofs of  both the local Langlands conjecture
and Shelstad's conjecture\Lspace \Lcitemark 12\Rcitemark \Rspace{} that
$R_{\psi_0,\s_0}$  is isomorphic to $R_{G^\circ}(\s_0).$
(See\Lspace \Lcitemark 1\Rcitemark \Rspace{} for the precise definitions
of $R_{\psi,\s_0}$ and $R_{\psi_0,\s_0}.$)
Here, by extending the definition of the standard intertwining operators
(cf Section 4) we show we can construct a group $R_G(\s_0)$
in a manner analogous to\Lspace \Lcitemark 6\Rcitemark \Rspace{}, 
and show that it has the correct parameterization properties (cf Theorems 4.16 and 4.17).
An argument, similar to the one given in\Lspace \Lcitemark 6\Rcitemark \Rspace{}
shows that if $G/G^\circ$ is cyclic, then
$R_G(\s_0)$ must be isomorphic to $R_{\psi,\s_0},$ assuming the latter exists
(cf Remark 4.18).

The second collection of results is an extension of some standard results in admissible
representation theory to the disconnected group $G$.  In order to develop a theory
consistent with the theory for connected groups, one needs to determine an appropriate
definition of parabolic subgroup.  There are several definitions in the literature already,
yet they do not always agree.  We use a definition of parabolic subgroup which is well
suited to our purposes.  Among the parabolic subgroups of $G$ we single out a collection
of parabolic subgroups which we call cuspidal.  
They have the property that they support discrete series
and supercuspidal representations and can be described as follows. Let $P^0$ be a parabolic subgroup of
$G^0$ with Levi decomposition
$M^0N$ and let $A$ be the split component of
$M^0$.  Let $M=C_G(A)$.  Then $P=MN$ is a cuspidal parabolic subgroup of $G$ lying over
$P^0$. We also say in this case that $M$ is a cuspidal Levi subgroup of $G$.  Thus cuspidal
parabolic subgroups of $G$ are in one to one correspondence with parabolic subgroups of
$G^0$.  

Using our definitions we can prove the following.  Let $M$ be a Levi subgroup of $G$ and let
$M^0 = M
\cap G^0$.  

\proclaim Lemma 1.1.  (i)  If $M$ is not cuspidal, then $M$ has no
supercuspidal representations, i.e., admissible 
representations with matrix coefficients which are
compactly supported modulo the center of $M$ and have zero constant term along the nil
radical of any proper parabolic subgroup of $M$. 
\hfil \break 
(ii)  If $M$ is cuspidal and $\pi $ is an irreducible admissible representation of $M$,
then $\pi$ is supercuspidal if and only if the restriction of $\pi $ to $M^0$ is
supercuspidal.
\hfil \break
(iii)  If $M$ is not cuspidal, then $M$ has no discrete series representations, ie.
unitary representations with matrix coefficients which are square-integrable modulo the
center of $M$.  \hfil \break 
(iv)  If $M$ is cuspidal and $\pi $ is an irreducible unitary representation of $M$, then
$\pi $ is discrete series if and only if the restriction of $\pi $ to $M^0$ is discrete
series.
\endproclaim 

Using Lemma 1.1 it is easy to extend the following theorem from the connected case to our
class of disconnected groups.

\proclaim Theorem 1.2. Let $\pi $ be an irreducible admissible (respectively
tempered) representation of
$G$.  Then there are a cuspidal parabolic subgroup $P=MN$ of $G$ and an irreducible
supercuspidal (respectively discrete series) representation $\sigma $ of $M$ such that $\pi $
is a subrepresentation of
$\Ind_P^G(\sigma )$. 
\endproclaim

Let $P_1=M_1N_1$ and $P_2=M_2N_2$ be cuspidal parabolic subgroups and let $\sigma _i$ be
irreducible representations of $M_i, i = 1,2$, which are either both supercuspidal or both
discrete series. By studying the orbits for the action of $P_1\times P_2$
on
$G$, we are able to extend the proof for the connected case to our situation and obtain the
following theorem.

\proclaim Theorem 1.3.  Let $P_1= M_1N_1,P_2=M_2N_2,\sigma _1, \sigma _2$ be as above.  Then
if $\pi _1 = \Ind _{P_1}^G(\sigma _1)$ and $\pi _2 = \Ind _{P_2}^G(\sigma _2)$ have a
nontrivial intertwining, then there is $y \in G$ so that $$M_2 = yM_1y^{-1} {\rm \ and \ }
\sigma _2
\simeq y\sigma _2y^{-1}.$$\endproclaim
 
The third question of study is the structure of $\pi = \Ind_P^G(\sigma )$ when $P=MN$ is a
cuspidal parabolic subgroup of $G$ and $\sigma $ is a discrete series representation of
$M$.  We show that, as in the connected case, the components of $\pi $ are naturally
parameterized using a finite group $R$.  As in the connected
case we first describe a collection of standard intertwining operators $R(w,\sigma )$ which
are naturally indexed by $w \in W_G(\sigma ) = N_G(\sigma )/M$, where $$N_G(\sigma ) = \{x
\in N_G(M):
\sigma ^x \simeq \sigma \}.$$   
 We prove that there is a cocycle $\eta $ so that
$$R(w_1w_2,\sigma ) = \eta (w_1,w_2) R(w_1,\sigma )R(w_2,\sigma ), w_1,w_2 \in W_G(\sigma
).$$

 Let $\sigma _0$ be an irreducible component of the
restriction of
$\sigma $ to $M^0$, and $P^0=M^0N = P \cap G^0$.  Then $\sigma \subset \Ind_{M^0}^M(\sigma
_0)$ so that 
$$\pi = \Ind _P^G(\sigma ) \subset \Ind _P^G(\Ind _{M^0}^M(\sigma _0)) \simeq \Ind
_{P^0}^G(\sigma _0).$$ 
Using the intertwining operators and $R$-group theory developed
earlier for the representation
$\Ind_{P^0}^G(\sigma _0)$, we can prove the following results.
First, the collection
$$\{ R(w,\sigma ), w \in W_G(\sigma )\}$$ spans the commuting algebra of $\pi $.   Second, let
$\Phi _1^+$ be the set of positive restricted roots for which the
rank one Plancherel measures of $\sigma _0$ are zero and
let $W(\Phi
_1)$ be the group generated by the reflections corresponding to the roots in $\Phi _1$.  
Then 
$W(\Phi _1)$ is naturally embedded as a normal subgroup of $W_G(\sigma )$ and 
$W_G(\sigma )$ is the semidirect product of $W(\Phi _1)$ and the group
$$R_{\sigma } = \{ w \in W_G(\sigma ): w\alpha >0 {\rm \ for \ all \ } \alpha \in \Phi
_1^+\}.$$ Finally, $R(w,\sigma )$ is scalar if $w \in W(\Phi _1)$.   This proves that
the operators $R(w,\sigma ), w \in R_{\sigma }$, span the intertwining algebra.  But we can
compute the dimension of the space of intertwining operators for $\Ind_P^G(\sigma )$,
again by comparison with  that of $\Ind_{P^0}^G(\sigma _0)$, and we find that it is equal to
$[R_{\sigma }]$.  Thus we have the following theorem.

\proclaim Theorem 1.4.  The $R(w,\sigma ), w \in R_{\sigma }$, form a basis for the algebra
of intertwining operators of $\Ind_P^G(\sigma )$.  \endproclaim

Just as in the connected case, we show that there are a finite central extension 
$$1 \rightarrow Z_{\sigma } \rightarrow \tilde R _{\sigma } \rightarrow R_{\sigma }
\rightarrow 1$$ over which $\eta $ splits and a character $\chi _{\sigma }$ of $Z_{\sigma
}$ so that the irreducible constituents of $\Ind_P^G(\sigma )$ are naturally parameterized
by the irreducible representations of $\tilde R _{\sigma }$ with $Z_{\sigma }$-central
character $\chi _{\sigma }$.  

Finally, we give a few examples which point out some of the subtleties
involved in working with disconnected groups.  For instance, we show that
if we do not restrict ourselves to cuspidal parabolic subgroups,
then the standard disjointness theorem for induced 
representations fails.
Examples such as these show why one must restrict to
induction from cuspidal parabolic subgroups in order
to develop a theory which is  consistent with that
for connected groups.

Many interesting problems involving disconnected groups remain.  
For example, the question of a Langlands classification is still unresolved,
and some of the results here on intertwining operators
may help in this direction.  One also hopes to remove
the condition that $G/G^\circ$ is abelian, and
extend all the results herein to that case.  Problems such as these
we leave to further consideration.

The organization of the paper is as follows.  In \S 2 we give the definition of parabolic
subgroup and prove Lemma 1.1 and Theorem 1.2.  The proof of Theorem 1.3 is in \S 3. 
The results on induction from a parabolic subgroup of $G^0$ to $G$ are in \S 4, and the
results on induction from a parabolic subgroup of $G$ to $G$, including Theorem 1.4, are in
\S 5.  Finally,
\S 6 contains examples that show what can go wrong when we induce from parabolic subgroups
of $G$ which are not cuspidal.

The first named author would like to thank the Mathematical
Sciences Research Institute in Berkeley, California,
for the pleasant and rich atmosphere in which some of the work
herein was completed.

\section {Basic Definitions}

Let $F$ be a locally compact, non-discrete, nonarchimedean field of characteristic zero. 
Let $G$ be a (not necessarily connected)
reductive $F$-group.  Thus $G$ is the set of $F$-rational points of a
reductive algebraic group over
$F$. Let
$G^0$ be the connected component of the identity in $G$.  We assume that $G/G^0$ is finite
and abelian.

The split component of $G$ is defined to be the maximal $F$-split torus lying in the
center of $G$.  Let $A$ be any $F$-split torus in $G$ and let $M = C_G(A)$. 
Then $M$ is a reductive $F$-group.  Now $A$ is called a
{\bf special torus} of $G$ if $A$ is the split component of $M$.  (Of course
$A$ is an $F$-split torus lying in the center of $M$.  The only question is
whether or not $A$ is maximal with respect to this property.)

\proclaim Lemma 2.1.  Let $A$ be a special torus of $G^0$.  Then $A$ is a
special torus of $G$.  \endproclaim

\proof   Let $M = C_G(A)$
and $M^0 = C_{G^0}(A) = M \cap G^0$.  Write $Z(M)$ and $Z(M^0)$ for the centers
of $M$ and $M^0$ respectively.   Now $A$ is the maximal $F$-split torus lying
in $Z(M^0)$ and $A \subset Z(M)$.  Suppose $A'$ is the maximal $F$-split torus
lying in $Z(M)$.  Then $A \subset A'$.  But $A'$ is a torus, so it is
connected.  Hence $A' \subset Z(M) \cap M^0 \subset Z(M^0)$.  Thus $A' \subset
A$ and so $A' = A$ is the split component of $M$.  \qed

\medbreak
\noindent {\bf Remark 2.2.}  The converse of Lemma 2.1 is not true.  For example, let $G =
O(2) = SO(2) \cup wSO(2)$ where $SO(2) \simeq F^{\times }$ is the group of $2\times 2$
matrices 
$$d(a) = \pmatrix{ a & 0\cr 0 & a^{-1} \cr}, a \in F^{\times}, 
{ \rm \ \ and \ \  }
 w =  \pmatrix{ 0 & 1\cr 1 & 0 \cr} $$ satisfies $wd(a)w^{-1} = d(a^{-1}), a \in
F^{\times }$.  Let $A = \{ d(1) \}$.  Then $M = C_G(A) = G$ and $Z(M) = \{ \pm d(1) \}
$. Thus $A$ is a special torus of $G$.  However
$M^0 = C_{G^0}(A) = G^0$ and $Z(M^0) = G^0$ is an $F$-split torus.  Hence $A$ is not the
maximal $F$-split torus in $Z(M^0)$ and so is not special in $G^0$. 

If $G$ is connected, then $A$ is a special torus of $G$ by the above
definition if and only if $A$ is the split component of a Levi component $M$ of
a parabolic subgroup of $G$.  We will define parabolic subgroups in the
non-connected case so that we have this property in the non-connected case also.

Let $A$ be a special torus of $G$ and let $M = C_G(A)$.  Then $M$ is called a Levi subgroup
of $G$.  The Lie algebra $L(G)$ can be decomposed into root spaces with respect to the roots
$\Phi $ of $L(A)$: $$ L(G) = L(G)_0 \oplus \sum _{\alpha \in \Phi } L(G)_{\alpha }$$
where $L(G)_0$ is the Lie algebra of $M$.  Let $\Phi ^+$ be a choice of positive
roots, and let $N$ be the connected subgroup of $G$ corresponding to $\sum
_{\alpha \in \Phi ^+} L(G)_{\alpha }.$   Since elements of $M$ centralize $A$ and $L(A)$,
they preserve the root spaces with respect to $L(A)$.  Thus $M$ normalizes $N$.  Now $P=MN$
is called a parabolic subgroup of $G$ and $(P,A)$ is called a $p$-pair of $G$.  The
following lemma is an immediate consequence of this definition and Lemma 2.1.  

\proclaim Lemma 2.3.  Let
$P^0 = M^0N$ be a parabolic subgroup of $G^0$ and let $A$ be the split component of $M^0$.
Let $M = C_G(A)$.  Then $P=MN$ is a parabolic subgroup of $G$ and $P \cap G^0 = P^0$. 
\endproclaim

\proclaim Lemma 2.4.  Let $P$ be a parabolic subgroup of $G$.  Then $P^0 = P \cap G^0$ is a
parabolic subgroup of $G^0$.  
\endproclaim

\proof  Let $P=MN$ be a parabolic subgroup of $G$ and let $A$ be the split component of $M$.  Let $M^0
= C_{G^0}(A) = M \cap G^0$.  Let $A_1$ be the split component of $M^0$.  Then $A \subset
A_1$ so that $C_{G^0}(A_1) \subset C_{G^0}(A) = M^0$.  But $A_1$ is in the center of
$M^0$, so that $M^0 \subset C_{G^0}(A_1)$.  Thus $C_{G^0}(A_1) = M^0$ so that $A_1$ is a
special torus in $G^0$ and $M^0$ is a Levi subgroup of $G^0$.  Let $\Phi $ and $\Phi _1$
denote the sets of roots of $L(A)$ and $L(A_1)$ respectively.  For each $\alpha _1\in \Phi
_1$, the restriction $r\alpha _1$ of $\alpha _1$ to $L(A)$ is non-zero since
$C_{G^0}(A_1) = C_{G^0}(A) = M^0$.  Let $\Phi ^+$ be
the set of positive roots used to define $N$.  Then $\Phi _1^+ = \{ \alpha _1
\in \Phi _1 : r\alpha _1 \in \Phi ^+ \}$ is a set of positive roots for $\Phi _1$ and 
$$\sum
_{\alpha \in \Phi ^+} L(G)_{\alpha } =\sum
_{\alpha _1 \in \Phi _1 ^+} L(G)_{\alpha _1 }.$$  Thus $P^0 = M^0N$ is a parabolic subgroup
of $G^0$.  \qed

We say the parabolic subgroup $P$ of $G$ lies over the parabolic subgroup $P^0$ of $G^0$
if $P^0 = P \cap G^0$.   We will also say the Levi subgroup $M$ of $G$ lies over the Levi
subgroup $M^0$ of $G^0$ if $M^0 = M \cap G^0$.   Lemma 2.4 and its proof show that every
parabolic (resp. Levi) subgroup of $G$ lies over a parabolic (resp. Levi) subgroup of
$G^0$.  

\medbreak
\noindent {\bf Remark 2.5.}  There can be more than one parabolic subgroup $P$ of $G$
lying over a parabolic subgroup $P^0$ of $G^0$.  For example, define $G=O(2)$ and
$G^0=SO(2)$ as in Remark 2.2.  Then $A = \{ d(1)\}$ and $A_0 = SO(2)$ are special vector
subgroups of $G$ corresponding to parabolic subgroups $P = O(2)$ and $P_0 = SO(2)$
respectively.  Both lie over the unique parabolic subgroup $SO(2)$ of $G^0$.  

\proclaim Lemma 2.6.  Let $P^0 = M^0N$ be a parabolic subgroup of $G^0$ and let $A$ be the
split component of $M^0$.  Let $M = C_G(A)$ and let $P = MN$.  Then if $P_1$ is any
parabolic subgroup of $G$ lying over $P^0$ we have $P \subset P_1$.  Further, $M$ is the
unique Levi subgroup lying over $M^0$ such that the split component of $M$ is equal to $A$.
\endproclaim

\proof  Write $P_1 = M_1N$ where $M_1$ lies over $M^0$.  Let $A_1$ be the split component
of $M_1$.  Then $A_1 \subset A $ so that $M = C_G(A) \subset C_G(A_1) = M_1$.   Clearly
$M_1 = M$ if and only if $A_1 = A$.  \qed

\medbreak
\noindent {\bf Remark 2.7.}  Lemma 2.6 shows that there is a unique smallest parabolic
subgroup $P$ of $G$ lying over $P^0$.  Although it is defined using a Levi decomposition
$P^0 = M^0N$ of $P^0$, it is independent of the Levi decomposition.  Recall that if
$M_1^0$ and $M_2^0$ are two Levi components of $P^0$ with split components $A_1$ and $A_2$
respectively, then there is $n \in N$ such that $A_2 = nA_1n^{-1}$ and $M^0_2 =
nM^0_1n^{-1}$.  Now if $M_i = C_G(A_i), i = 1,2$, we have
$M_2 = C_G(A_2) = nC_G(A_1)n^{-1} = nM_1n^{-1}$, and $M_2N =  nM_1n^{-1}N = M_1N$ since
$M_1$ normalizes $N$.  

Let $Z$ be the split component of $G$.  We let $ C_c^{\infty }(G,Z)$ denote the space of
all smooth complex-valued functions on $G$ which are compactly supported modulo $Z$.   We
say $f \in C_c^{\infty }(G,Z)$ is a cusp form if for every proper parabolic subgroup $P=MN$
of $G$,  $$\int _N f(xn)dn = 0 \ \ \ \ \forall x \in G.$$
Let $\ ^0 \calA (G)$ denote the set of cusp forms on $G$.  We say $G$ is cuspidal if 
$\ ^0 \calA (G) \not = \{ 0 \} $.  We know that every connected $G$ is cuspidal.

\proclaim Lemma 2.8.   $G$ is cuspidal if and only if the split component of $G$ is equal to
the split component of $G^0$.  Moreover, if $G$ is cuspidal, then a subgroup $N$ of $G$ is
the nilradical of a proper parabolic subgroup of $G$ if and only if $N$ is the nilradical
of a proper parabolic subgroup of $G^0$.   If $G$ is not cuspidal, then $G$ has a proper
parabolic subgroup $G_1$ with nilradical $N_1 = \{ 1 \}$.   \endproclaim

\proof  First suppose that $G$ and $G^0$ have the same split component $Z$.   Let $f \not =
0 \in \ ^0 \calA (G^0)$.  Define $F:G \rightarrow \bfC$ by $F(x) = f(x), x \in G^0$, and
$F(x) = 0, x \not \in G^0$.  Then $F \in C_c^{\infty }(G,Z)$ and is non-zero.  Let $P=MN$ be
any proper parabolic subgroup of $G$.  Then $N \subset G^0$, so for all $n \in N, x \in G$,
$xn \in G^0 $ if and only if  $x  \in G^0$.  Thus for $x \not \in G^0$,
 $$\int _N F(xn)dn = 0 $$
while for $x \in G^0$,
 $$\int _N F(xn)dn = \int _N f(xn)dn .$$

Now $P^0 = P \cap G^0 = M^0N$ is a parabolic subgroup of
$G^0$.  Suppose that $P^0 = G^0$.  Then $P$ lies over $G^0$ so that by Lemma 2.6, $G \subset
P$.  This contradicts the fact that $P$ is a proper parabolic subgroup of $G$.  Thus $P^0 =
M^0N $ is a proper parabolic subgroup of $G^0$. Since $f$ is a cusp form for $G^0$ we
have
 $$\int _N f(xn)dn = 0 \ \ \ \ \forall x \in G^0.$$ 
Thus $F $ is a non-zero cusp form for $G$.  

The above argument also showed that if $P=MN$ is a proper parabolic subgroup of $G$, then
$P^0 = M^0N$ is a proper parabolic subgroup of $G^0$.  Conversely, if $P^0 = M^0N$ is a
proper parabolic subgroup of $G^0$ and $P=MN$ is any parabolic subgroup of $G$ lying over
$P^0$, then clearly $P \not = G$.

Conversely, suppose that $G$ and $G^0$ do not have the same split component.  Let $Z$ be
the split component of $G^0$ and define $G_1 = C_G(Z)$.  By Lemma 2.6, $G_1$ is a proper
parabolic subgroup of $G$.  Further since $G_1$ lies over $G^0$ its nilradical is $N_1 =
\{1 \}$.  Now if $F$ is any cusp form on $G$ and $x \in G$, we  have
$$F(x) = \int _{N_1} F(xn)dn = 0.$$  Thus $G$ has no non-zero cusp forms and so is not
cuspidal.  \qed
 
\medbreak
\noindent {\bf Example 2.9.}  Let $G = O(2)$ as in Remarks 2.2 and 2.5.  Then $SO(2)$ is a
cuspidal parabolic subgroup of $G$ and $O(2)$ is not cuspidal.  

We can sum up the proceeding lemmas in the following proposition.

\proclaim Proposition 2.10.  Let $P^0 = M^0N$ be a parabolic subgroup of $G^0$.  Then there
is a unique cuspidal parabolic subgroup $P=MN$ of $G$ lying over $P^0$.  It is contained in
every parabolic subgroup of $G$ lying over $P^0$, and is defined by $M =
C_G(A)$ where $A$ is the split component of $M^0$.  \endproclaim 

Now that we have parabolic subgroups of $G$, we want to study parabolic induction of
representations.  Many of the most basic notions of representation theory are defined in
\Lcitemark 13\LIcitemark{}, chapter 1\RIcitemark \Rcitemark \Rspace{} for any totally disconnected group.  In particular, admissible
representations of $G$ are defined and the following is an easy consequence of the
definition.

\proclaim Lemma 2.11.  Let $\Pi $ be a representation of $G$.  Then $\Pi $ is admissible if
and only if $\Pi |_{G^0}$, the restriction  of $\Pi $ to $G^0$, is admissible.  \endproclaim

Further, the results of Gelbart and Knapp regarding induction and restriction between a
totally disconnected group $G$ and an open normal subgroup $H$ with $G/H$ finite abelian
can be applied to $G$ and $G^0$.  
If $\pi $ is any admissible representation of $G^0$ on $V$, we will let
$\Ind^G_{G^0}(\pi )$ denote the representation of $G$ by left translations on
${\cal H} = \{ f:G \rightarrow V: f(gg_0) = \pi (g_0)^{-1} f(g),\forall  g \in G,
g_0 \in G^0 \}$.  

\proclaim Lemma 2.12. 
{\rm{\bf (Gelbart-Knapp\Lspace \Lcitemark 5\Rcitemark \Rspace{})}}
Let $\Pi $ be an irreducible admissible
representation of $G$.  Then $\Pi |_{G^0}$ is a finite direct sum of irreducible
admissible representations of $G^0$.  Let $\pi $ be an irreducible constituent of $\Pi
|_{G^0}$ which occurs with multiplicity $r$.  Then
$$\Pi |_{G^0} \simeq r \sum _{g \in G/G_{\pi } } \pi ^g$$ where
$G_{\pi } = \{ g \in G: \pi ^g \simeq \pi \}$.  \endproclaim

\proclaim Lemma 2.13. {\rm{\bf (Gelbart Knapp\Lspace \Lcitemark 5\Rcitemark \Rspace{})}}   
Let $\pi $ be an irreducible admissible
representation of $G^0$.  Then there is an irreducible admissible representation $\Pi $ of
$G$ such that
$\pi $ occurs in the restriction of $\Pi $ to $G^0$ with multiplicity $r > 0$.  Let $X $
denote the group of unitary characters of $G/G^0$ and let
$X(\Pi ) = \{ \chi \in X: \Pi \otimes \chi \simeq \Pi \}$.  Then
$$\Ind _{G^0}^G(\pi ) \simeq r \sum _{\chi \in X/X(\Pi )} \Pi \otimes \chi $$ is the
decomposition of
$\Ind _{G^0}^G(\pi )$ into irreducibles and $r^2 [X/X(\Pi )] = [G_{\pi }/G^0]$.  \endproclaim

The following result was proved by Gelbart and Knapp in the case
where the restriction is multiplicity one\Lspace \Lcitemark 5\Rcitemark \Rspace{}.  
Tadic\Lspace \Lcitemark 14\Rcitemark \Rspace{} refined their result in the connected case.  We now prove
the more general result.

\proclaim Lemma 2.14. 
Suppose that $G$ is a totally disconnected group,
and $H$ is a closed normal subgroup, with $G/H$ a finite abelian
group.  If $\Pi_1$ and $\Pi_2$ are irreducible admissible reprsentations of
$G,$ which have a common constituent upon restriction to $H,$
then $\Pi_2\simeq\Pi_1\otimes\chi,$ for some character $\chi$
with $\chi|_H\equiv 1.$
\endproclaim

\proof  If the multiplicity of the restrictions is one,
then this result holds by Gelbart-Knapp\Lspace \Lcitemark 5\Rcitemark \Rspace{}.  
In particular, if
$|G/H|$ is prime, the statement is true.  
We proceed by induction.  We know the Lemma holds
when $|G/H|=2.$
Suppose
the statement is true whenever $|G_1/H|<n.$  
Suppose $|G/H|=n.$
We may assume $n$
is composite, so write $n=km,$ with $1<k<n.$
Let $H\subset G_1\subset G$ with $|G/G_1|=k.$  If
$\Pi_1|_{G_1}$ and $\Pi_2|_{G_2}$ have a common constituent,
then, by our inductieve hypothesis, there is a $\chi$
with $\chi|_{G_1}\equiv 1$
with $\Pi_2\simeq\Pi_1\otimes\chi.$
Since $(G/G_1\hat)\subset(G/H\hat),$ we are done, in this case.

Now suppose that $\tau$
is an irreducible subrepresentarion of both $\Pi_1|_{H}$ and
$\Pi_2|_{H}.$
Then,
there are constituents $\Omega_i\subset\Pi_i|_{G_1}$
so that $\tau\subset\Omega_i|_{H}.$ By the inductive hypothesis
$\Omega_2=\Omega_1\otimes\chi,$ for some $\chi$
of $G_1$ whose restriction to $H$ is trivial.
But, since $G_1/H\subset G/H$ is abelian we can extend
$\chi$ to a character $\eta$ of $G/H.$  Note that $(\Pi_1\otimes\eta)|{G_1}$
has $\Omega_1\otimes\chi\simeq\Omega_2$ as a constituent,
so, as we have seen above, $\Pi_1\otimes\eta\otimes\omega\simeq\Pi_2,$
for some character  $\omega$ of $G$ whose
restriction to $G_1$ is trivial.  Thus, the statement holds by induction.
\qed

Let $(\pi ,V)$ be an admissible representation of $G$ and let $\calA (\pi )$ denote its
space of matrix coefficients.  We say $\pi $ is supercuspidal if $\calA (\pi ) \subset
\ ^0 \calA (G)$.  Of course if $G$ is not cuspidal, then $\ ^0 \calA (G) = \{ 0 \}$ so
that $G$ has no supercuspidal representations.  If $P=MN$ is any parabolic subgroup of
$G$, define $V(P) = V(N)$ to be the subspace of $V$ spanned by vectors of the form
$\pi (n) v -v, v \in V, n \in N$.  Then we say $\pi $ is J-supercuspidal if $V(P) =V$ for
every proper parabolic subgroup $P$ of $G$.  If $G$ is not cuspidal, then by Lemma 2.8
there is a proper parabolic subgroup $G_1$ of $G$ with nilradical $N_1 = \{ 1 \}$.  For
any admissible representation $(\pi ,V)$ of $G$, $V(G_1) = V(N_1) = \{ 0 \} \not = V$ so
that $\pi $ is not J-supercuspidal.  Thus $G$ has no J-supercuspidal representation.  

Suppose now that $G$ is cuspidal and let $(\pi ,V)$ be an irreducible admissible
representation of $G$.  Let $(\pi _0, V)$ denote the restriction of $\pi $ to $G^0$.

\proclaim Lemma 2.15.  Assume that $G$ is cuspidal.  Then $\pi $ is J-supercuspidal if and
only if $\pi _0$ is J-supercuspidal if and only if any irreducible constituent of $\pi
_0$ is J-supercuspidal.  \endproclaim

\proof  Since $G$ is cuspidal, by Lemma 2.8 the set of nilradicals of proper parabolic
subgroups is the same for both $G$ and $G^0$.  Thus $\pi $ is  J-supercuspidal if and only
if $\pi _0$ is J-supercuspidal.  Moreover, since by Lemma 2.12 the irreducible constituents
of
$\pi _0$ are all conjugate via elements of $G$, it is clear that $\pi _0$ is J-supercuspidal
if and only if every irreducible constituent of $\pi _0$ is J-supercuspidal if and only if
any  irreducible constituent of $\pi
_0$ is J-supercuspidal.  \qed

\proclaim Lemma 2.16.  Assume that $G$ is cuspidal.  Then $\pi $ is supercuspidal if and
only if $\pi _0$ is supercuspidal if and only if any irreducible constituent of $\pi
_0$ is supercuspidal.  \endproclaim

\proof  Let $Z$ denote the split component of $G$.  By Lemma 2.8 it is also the split
component of $G^0$.

 Assume that $\pi $ is supercuspidal.  Thus $\calA (\pi ) \subset \ ^0
\calA (G)$. Let $f_0$ be a matrix coefficient of $\pi _0$.  Then there is a matrix
coefficient $f $ of $\pi $ so that $f_0$ is the restriction of $f$ to $G^0$.  Now $f \in \
^0 \calA (G)$. Since $f $ smooth and compactly supported modulo $Z$, so is $f_0$. 
Further, by Lemma 2.8 the nilradicals of proper parabolic subgroups are the same for both
$G$ and $G^0$.  Thus $f_0$ will satisfy the integral condition necessary to be a cusp form
on $G^0$.  Hence $\calA (\pi _0 ) \subset \ ^0 \calA (G^0)$ so that $\pi _0$ is
supercuspidal.

Conversely, suppose that $\pi _0$ is supercuspidal.  Let $\pi _1$ be an irreducible
constituent of $\pi _0$.  Then $\pi \subset \Ind_{G^0}^G(\pi _1)$ so that every matrix
coefficient of $\pi $ is a matrix coefficient of the induced representation.  But since
$G^0$ is a normal subgroup of finite index in $G$, the restriction of $\Ind_{G^0}^G(\pi _1)$ to $G^0$ is
equivalent to $\sum _{x \in G/G^0} \pi _1^x$.  Thus  matrix coefficients of the induced
representation can be described as follows.   
Let $f$ be a matrix coefficient of 
$\Ind_{G^0}^G(\pi _1)$ and fix $g \in G$.  Then there are matrix coefficients $f_x$ of $\pi
_1^x, x \in G/G^0$, (depending on both $f$ and $g$) so that for all $g_0 \in G^0$,
$$f(gg_0) = \sum _{x \in G/G^0} f_x(g_0).$$  Since $\pi _1$ is supercuspidal, so is $\pi
_1^x$ for any $x \in G/G^0$, and so each $f_x \in \ ^0 \calA (G^0)$.  Thus the restriction
of $f$ to each connected component of $G$ is smooth and compactly supported modulo $Z$. 
Also, if $N$ is the nilradical of any proper parabolic subgroup of $G$, then
$$\int _N f(gn)dn = \sum _{x \in G/G^0} \int _N f_x(n)dn = 0$$ since by Lemma 2.8, $N$ is
also the nilradical of a proper parabolic subgroup of $G^0$.  \qed

\proclaim Proposition 2.17.   Assume that $G$ is cuspidal and let $\pi $ be an irreducible
admissible representation of $G$.  Then $\pi $ is supercuspidal if and
only if $\pi $ is J-supercuspidal.   \endproclaim

\proof   This is an immediate consequence of Lemmas 2.15 and 2.16 and the corresponding
result in the connected case.  \qed

We now drop the assumption that $G$ is cuspidal.
Let $P=MN$ be a parabolic subgroup of $G$ and let $\sigma $ be an admissible
representation of $M$.  Then we let $\Ind_P^G(\sigma )$ denote the
representation of $G$ by left translations on 
$${\cal H} = \{ f \in C^{\infty }(G,V): f(gmn) = \delta _P^{-{1\over 2}}(m)
\sigma (m)^{-1}f(g), \forall g \in G, m \in M, n \in N\}.$$
Here $\delta _P$ denotes the modular function of $P$.  

\proclaim Theorem 2.18.  Let $\pi $ be an irreducible
admissible representation of
$G$.  Then there are a cuspidal parabolic subgroup $P= MN$ of $G$ and an irreducible 
supercuspidal representation $\sigma $ of $M$ such that $\pi $ is a subrepresentation of
$\Ind _P^G(\sigma )$.  \endproclaim  

\medbreak
\noindent {\bf Remark 2.19.}  We will see in Corollary 3.2 that the group $M$ and
supercuspidal representation $\sigma $ in Theorem 2.18 are unique up to conjugacy.

\proof  Let $\rho $ be an irreducible constituent of the restriction of $\pi $ to $G^0$. 
Then $\pi \subset \Ind_{G^0}^G(\rho )$. 
Since $\rho $ is admissible, there are a parabolic subgroup $P^0= M^0N$ of $G^0$ and
an irreducible  supercuspidal representation $\tau $ of $M^0$ such that 
$\rho \subset \Ind _{P^0}^{G^0}(\tau )$.  Thus
$$\pi \subset \Ind_{G^0}^G (\rho ) \subset \Ind _{G^0}^G (\Ind _{P^0}^{G^0}(\tau )) \simeq
\Ind _{P^0}^G (\tau ).$$   

 Let $P=MN$ be the unique cuspidal parabolic subgroup of $G$
lying over $P^0$.  Let $\sigma $ be an irreducible admissible representation of $M$ such that
$\tau $ is contained in the restriction of $\sigma $ to $M^0$.  By Lemma 2.16, $\sigma $ is
supercuspidal.
By Lemma 2.13 applied to $M$ and $M^0$ we have
$$\Ind_{M^0}^M(\tau ) \simeq s \sum _{\eta \in Y/Y(\sigma )} \sigma \otimes \eta $$ where $Y$
is the group of unitary characters of $M/M^0$.  Since $\tau $ is contained in the
restriction of $\sigma \otimes \eta $ to $M^0$ for any $\eta $, all the representations
$\sigma \otimes \eta $ are supercuspidal.  Now
$$\pi \subset \Ind _{P^0}^G (\tau )  \simeq
 s \sum _{\eta \in Y/Y(\sigma )}\Ind _P^G( \sigma \otimes \eta ) .$$  Thus
$\pi \subset \Ind _P^G( \sigma \otimes \eta ) $ for some $\eta $. \qed

Let $\calA (G) = \cup _{\pi }\calA (\pi )$ where $\pi $ runs over the set of all admissible
representations of $G$.  Similarly we have $\calA (G^0)$ and because of Lemma 2.11 it is
clear that if $f \in \calA (G)$, then $f|_{G^0} \in \calA (G^0)$.  Define the subspace 
$\calA_T(G^0)\subset \calA (G^0)$ as in\Lspace \Lcitemark 13\LIcitemark{},\S 4.5\RIcitemark \Rcitemark \Rspace{}.  It is the set of functions in $\calA
(G^0)$ which satisfy the weak inequality.  Define
$$\calA _T(G) = \{ f \in \calA (G): l(x)f|_{G^0} \in \calA _T(G^0) {\rm \ for \ all \ } x
\in G \}$$ where $l(x)f$ denotes the left translate of $f$ by $x$.  In other words, $\calA
_T(G)$ is the set of functions in
$\calA (G)$ which satisfy the weak inequality on every connected component of $G$. 
If $\pi $ is an admissible representation of $G$, we say $\pi $ is tempered if $\calA (\pi
) \subset \calA _T(G)$.  The following lemma is easy to prove using the properties of
matrix coefficients of $\pi $ and $\pi _0$ from the proof of Lemma 2.16.

\proclaim Lemma 2.20.  Let $\pi $ be an irreducible admissible representation of $G$.  Then
$\pi $ is tempered if and only if $\pi _0$ is tempered if and only if any irreducible
constituent of $\pi _0$ is tempered.  \endproclaim

Let $\pi $ be an irreducible unitary representation of $G$ and let $Z$ be the split
component of $G$.  We say that $\pi $ is discrete series if $\calA (\pi ) \subset
L^2(G/Z)$.  Every unitary supercuspidal representation is discrete series since its
matrix coefficients are compactly supported modulo $Z$.  

\proclaim Lemma 2.21.  If $G$ is not cuspidal, then $G$ has no discrete series
representations.  If $G$ is cuspidal, then $\pi $ is discrete series if and only if $\pi
_0$ is discrete series if and only if any irreducible constituent of $\pi _0$ is
discrete series.  \endproclaim

\proof  Suppose that $G$ is not cuspidal.  Then the split component $Z$ of $G$ is a proper
subgroup of the split component $Z_0$ of $G^0$.  Let $\pi $ be any irreducible unitary
representation of $G$.  
Then there is an irreducible unitary representation $\pi _1$ of
$G^0$ so that $\pi $ is contained in $\Ind _{G^0}^G(\pi _1)$.  Thus as in the proof of Lemma
2.16, for any matrix coefficient $f$ of $\pi $ and any $g \in G$, we have matrix coefficients
$f_x$ of $\pi
_1^x, x \in G/G^0$, so that for all $g_0 \in G^0$,
$$f(gg_0) = \sum _{x \in G/G^0} f_x(g_0).$$     
Let $\omega $ be the $Z_0$-character of $\pi
_1$.  Then for any $z \in Z_0, g_0 \in G^0$, we have
$$f(gg_0z) = \sum _{x \in G/G^0} f_x(g_0z) = \sum _{x \in G/G^0}\omega ^x(z) f_x(g_0).$$ 
Thus $z \mapsto f(gg_0z), z \in Z_0,$ is a finite linear combination of unitary characters
of $Z_0$, and cannot be square-integrable on $Z_0/Z $ unless it is zero.  Now if $f$ is
square-integrable on $G/Z$,  then
$g_0 \mapsto f(gg_0)$ is square-integrable on
$G_0/Z$ for all coset representatives $g \in G/G^0$,
and $z \mapsto f(gg_0z)$ must be square-integrable on $Z_0/Z$ for
almost all $g_0$, so that $f(gg_0z)$ must be zero for
almost all $g_0,z$, and $f = 0$.  

Suppose that $G$ is cuspidal.  Let $\pi $ be a discrete series representation of $G$. 
Let $f_0$ be a matrix coefficient of $\pi _0$.  Then there is a matrix
coefficient $f $ of $\pi $ so that $f_0$ is the restriction of $f$ to $G^0$.  Since $f$ is
square-integrable on $G/Z$, certainly $f_0$ is square-integrable on $G^0/Z$. 

 Conversely,
suppose that $\pi _0$ is discrete series.  Let $\pi _1$ be an irreducible constituent of
$\pi _0$.  For any matrix coefficient $f$ of $\pi $ we have
$$\int _{G/Z} |f(g)|^2 d(gZ) = \sum _{g \in G/G^0} \int _{G^0/Z} |f(gg_0)|^2 d(g_0Z).$$
As above, for fixed $g \in G$, we have
matrix coefficients $f_x$ of $\pi
_1^x, x \in G/G^0$, so that for all $g_0 \in G^0$,
$$f(gg_0) = \sum _{x \in G/G^0} f_x(g_0).$$   Thus 
$$ (\int _{G^0/Z} |f(gg_0)|^2 d(g_0Z))^{1\over 2} \leq  
\sum _{x \in G/G^0} (\int _{G^0/Z} |f_x(g_0)|^2 d(g_0Z))^{1\over 2} < \infty $$
since every $f_x$ is square-integrable on $G^0/Z$.  Thus $f$ is square-integrable on
$G/Z$.  \qed

The following theorem can be proven in the same way as Theorem 2.18 using Lemmas 2.20 and
2.21.

\proclaim Theorem 2.22.  Let $\pi $ be an irreducible
tempered representation of
$G$.  Then there are a cuspidal parabolic subgroup $P= MN$ of $G$ and an irreducible 
discrete series representation $\sigma $ of $M$ such that $\pi $ is a subrepresentation of
$\Ind _P^G(\sigma )$.  \endproclaim  

\medbreak
\noindent {\bf Remark 2.23.}  We will see in Corollary 3.2 that the group $M$ and
discrete series representation $\sigma $ in Theorem 2.22 are unique up to conjugacy.

\section {Intertwining Operators}

Let $G$ be a reductive $F$-group with $G/G^0$ finite and abelian as in \S 2.
If $\pi _1$ and $\pi _2$ are admissible representations of $G$, we let $J(\pi _1,\pi _2)$
denote the dimension of the space of all intertwining operators from $\pi _2$ to $\pi _1$.
We want to prove the following theorem.

Let $(P_i,A_i), i = 1,2$, be cuspidal parabolic pairs of $G$ with $P_i = M_iN_i$ the
corresponding Levi decompositions.  Let $\sigma _i$ be an irreducible admissible
representation of $M_i$ on a vector space $V_i$, and let $\pi _i = \Ind_{P_i}^G(\sigma _i)$,
where we use normalized induction as in \S 2.
Let $W = W(A_1|A_2)$ denote the set of mappings $s:A_2 \rightarrow
A_1$ such that there exists $y_s \in G$ such that $s(a_2) = y_sa_2y_s^{-1}$ for all $a_2
\in A_2$.

\proclaim Theorem 3.1.  Assume that $\sigma _1$ and $\sigma _2$ are either both
supercuspidal representations, or both discrete series representations.  If $A_1$ and
$A_2$ are not conjugate, then $J(\pi _1,\pi _2 ) = 0$.  Assume $A_1$ and $A_2$ are
conjugate.  Then
$$J(\pi _1, \pi _2) \leq \sum _{s \in W} J(\sigma _1, \sigma _2^{y_s}).$$
\endproclaim

\proclaim Corollary 3.2.  Let $\sigma _1$ and $\sigma _2$ be as above.  Then
$J(\pi _1, \pi _2) = 0$ unless there is $y \in G$ such that
$M_1 = M_2^y, \sigma _1 \simeq \sigma _2^y$.  
\endproclaim

\proclaim Corollary 3.3.  Let $\sigma _1$ and $\sigma _2$ be irreducible discrete series
representations, and suppose that $\pi _1$ and $\pi _2$ have an irreducible
constituent in common.   Then $\pi _1 \simeq \pi _2$.  \endproclaim

\proof  In this case $J(\pi _1, \pi _2) \not = 0$ so by Corollary 3.2 there is
$y \in G$ such that
$M_1 = M_2^y, \sigma _1 \simeq \sigma _2^y$.  Thus 
$$\pi _1 = \Ind _{M_1N_1}^G(\sigma _1) = \Ind _{M_2^yN_1}^G(\sigma _2^y) \simeq
\Ind _{M_2N_2'}^G (\sigma _2)$$ where $N_2' = N_1^{y^{-1}}$.  Now $P_2 = M_2N_2$ and
$P_2'=M_2N_2'$ are two cuspidal parabolic subgroups of $G$ with the same Levi component
$M_2$.  We will see in Corollary 5.9 that there is an equivalence
$R(P_2':P_2:\sigma _2)$ between $\pi _2 = \Ind _{P_2}^G(\sigma _2)$ and $\pi _2' =
\Ind _{P_2'}^G(\sigma _2)$.  Thus $$\pi _1 \simeq \pi _2' \simeq \pi _2.$$
\qed

Let $W(\sigma _1) = \{ s \in W(A_1,A_1): \sigma _1 ^s \simeq \sigma _1 \}$.  We say
$\sigma _1$ is unramified if $W(\sigma _1) = \{1 \}$.  

\proclaim Corollary 3.4.  Assume that $\sigma _1$ is discrete series.  Then $J(\pi _1, \pi
_1) \leq [W(\sigma _1)]$.  In particular, if $\sigma _1$ is unramified, then $\pi _1$ is
irreducible.  \endproclaim 

In order to prove Theorem 3.1 for discrete series representations we will need the following
results about dual exponents.  Let $(\pi ,V)$ be an irreducible representation of $G$ and let $V'$ denote
the algebraic dual of $V$.  For $x \in G$, let $\pi (x)^t:V' \rightarrow V'$ denote the transpose of
$\pi (x)$.
Let $(P,A)$ be a cuspidal parabolic pair in $G$.  Then a quasi-character $\chi $ of $A$ is
called a dual exponent of $\pi $ with respect to $(P,A)$ if there is a nonzero $\phi \in V'$
such that for all $\barn \in \barN, a \in A$,
 $$ \pi(\barn )^t \phi = \phi \ \ \ {\rm and } \ \ \ \ \pi (a)^t \phi = \delta _{\barP }
^{1\over 2}(a) \chi (a) \phi .\eqno(*)$$
We will write $Y_{\pi }(P,A)$ for the set of all dual exponents of $\pi $ with respect to
$(P,A)$.  

Let $V = \sum _{i=1}^k V_i$ be the decomposition of $V$ into $G^0$ invariant subspaces and
let $\pi _i$ be the irreducible representation of $G^0$ on $V_i, 1 \leq i \leq k$.

\proclaim Lemma 3.5. Let $(P,A)$ be a cuspidal parabolic pair
in $G$. Then $$Y_{\pi }(P,A) = \cup _{i=1}^k Y_{\pi
_i}(P^0,A).$$ \endproclaim

\proof   Let $\chi \in Y_{\pi }(P,A)$ and let $\phi \in V'$ be a nonzero functional
satisfying (*).  Then there is $1 \leq i \leq k$ such that $\phi _i$, the restriction of
$\phi $ to $V_i$, is nonzero.  Now for any $v_i \in V_i$ and $g_0
\in G^0$ we have
$$<\pi _i (g_0)^t \phi _i , v_i> = <\pi  (g_0)^t \phi , v_i>.$$
Now since $\barN$ and $A$ are both contained in $G^0$ it is easy to check that 
$$ <\pi _i (\barn )^t \phi _i, v_i>  = <\phi _i, v_i> \ \ \ {\rm and } \ \ \ <\pi _i (a)^t
\phi _i,v_i> = \delta _{\barP ^0} ^{1\over 2}(a) \chi (a) <\phi _i, v_i> \eqno(**)$$
for all $v_i \in V_i, \barn \in \barN , a \in A$.  Thus $\chi \in  Y_{\pi
_i}(P^0,A)$.

Now assume that $\chi \in  Y_{\pi _i }(P^0,A)$ for some $1 \leq i \leq k$.
and let $\phi _i \in V_i'$ be a nonzero functional
satisfying (**).   Now $\pi $ is contained in the induced representation $\Ind_{G^0}^G(\pi
_i)$ so we can realize $\pi $ on a subspace $V$ of
$$\calH = \{ f:G \rightarrow V_i: f(gg_0) = \pi _i(g_0)^{-1}f(g), g \in G, g_0 \in G^0
\}$$ with the action of $\pi $ given by left translation on the functions.
Now define $\phi \in V'$ by 
$$<\phi ,f> = <\phi _i, f(1)>, f \in V.$$  Then $\phi \not = 0$ and
it is easy to check that for all $g_0 \in G^0, f \in V$, we have
$$<\pi (g_0)^t\phi , f> = <\pi _i (g_0)^t \phi _i, f(1)>.$$  From this it easily follows
that since $\phi _i$ satisfies (**), $\phi $ satisfies (*).   Thus 
 $\chi \in Y_{\pi }(P,A)$.  \qed

\proclaim Lemma 3.6.  Assume that $G$ is cuspidal and let $\pi $ be a discrete series
representation of
$G$.  Then
$Y_{\pi }(P,A) \cap \hat A = \emptyset $ 
for every cuspidal parabolic pair $(P,A) \not = (G,Z)$.
  \endproclaim

\proof  This follows easily from Lemma 2.21, Lemma 3.5, and the corresponding result for
connected groups.  \qed 

Let $(P_0^0,A_0)$ be a minimal p-pair in $G^0$ and let $P_0$ be the cuspidal parabolic
subgroup lying over $P^0_0$.  We will call $(P_0,A_0)$ a minimal parabolic pair in $G$.

\proclaim Lemma 3.7.  Let $(P,A)$ be any parabolic pair in $G$.  Then there is $x\in G^0$
such that $P_0 \subset xPx^{-1}, xAx^{-1} \subset A_0$.  \endproclaim

\proof   Let $M=C_G(A)$ and let
$P^0 = P \cap G^0, M^0 = M \cap G^0$.  Let $A'$ be the split component of $M^0$.  Thus $A
\subset A'$.   Now $(P^0, A')$ is a p-pair in $G^0$ so there is $x \in G^0$ such that
$P_0^0 \subset xP^0x^{-1}$ and $xAx^{-1} \subset xA'x^{-1} \subset A_0$.   Now $M_0 =
C_G(A_0)
\subset C_G(xAx^{-1}) = xMx^{-1}$ so that $P_0 = M_0P_0^0 \subset (xMx^{-1})(x P^0x^{-1})
= xPx^{-1}$.  \qed

We will say a parabolic pair $(P,A)$ is standard with respect to the minimal
parabolic pair $(P_0,A_0)$ if $P_0 \subset P$ and $A \subset A_0$.  We say $(P,A)$ is
semi-standard with respect to $(P_0,A_0)$ if $A \subset A_0$.

   Fix a minimal parabolic pair
$(P_0,A_0)$ of
$G$.    
Let $N_G(P_0,A_0)$ denote the set of all elements of $G$ which normalize both $P_0$
and $A_0$.  Write $W_G(P_0,A_0) = 
N_G(P_0,A_0)/M^0_0$.  If $(P,A)$ is any parabolic pair of $G$ which is
standard with respect to $(P_0,A_0)$ and $M = C_G(A)$, write $N_M(P_0,A_0) = M \cap
N_G(P_0,A_0)$ and  $W_M(P_0,A_0) = 
N_M(P_0,A_0)/M^0_0$.

\proclaim Lemma 3.8.  We can write $P$ as a disjoint union 
$$P = \cup _{w \in W_M(P_0,A_0)} wP^0.$$  \endproclaim

\proof  We first prove the result when $P = G$.  Let $y \in G$.  Then $(yP_0^0y^{-1},
yA_0y^{-1})$ is a minimal parabolic pair in $G^0$ and hence is conjugate to $(P^0_0,A_0)$ via
an element of $G^0$.  Thus there is $g \in G^0$ such that $gy$ normalizes both $A_0$ and
$P^0_0$.  But then $gy$ also normalizes $M_0= C_G(A_0)$ and $P_0 = M_0P_0^0$.  Thus $gy \in
N_G(P_0, A_0) $. 

Thus the coset $yG_0$ has a representative $n \in N_G(P_0, A_0) $ which depends only on
the coset $w$ of $n$ in $W_G(P_0,A_0)$.   Since $N_G(P_0, A_0) \cap G^0 = M_0^0$,
$n_1,n_2 \in  N_G(P_0, A_0)$ determine the same coset of $G^0$ in $G$ just in case they
are in the same  coset in $W_G(P_0,A_0)$.  Thus we have the disjoint union
$$G = \cup _{w \in W_G(P_0,A_0)} wG^0.$$ 

Now let $(P,A)$ be arbitrary.  There is a subset $\theta $ 
of the simple roots $\Delta $ of $(P_0^0,A_0)$ corresponding to $P^0$ so that $P^0 =
P^0_{\theta }$.  Note that $W_G(P_0,A_0)$ acts on $\Delta $ and that  $w^{-1}P_{\theta
}^0w = P^0_{w \theta }$. 
 Let $w \in
W_G(P_0,A_0)$ and suppose that $P \cap wG^0 \not = \emptyset $.  Let $wg_0 \in P \cap
wG^0, g_0 \in G^0$.  Then $wg_0$ normalizes $P^0$ so that $w^{-1}P^0w = g_0P^0g_0^{-1}$ is
conjugate to
$P^0$ in $G^0$.  But $w^{-1}P^0w = P^0_{w \theta }$ is a standard parabolic subgroup 
of $G^0$ and so  $ P^0_{w \theta } = g_0P^0g_0^{-1} = P^0$.  Hence $g_0 \in N_{G^0}(P^0)
= P^0, w\theta = \theta $, and $w$ has a representative in $P$.  Hence $w$ normalizes
$A = \{ a \in A_0: \alpha (a) = 1$ for all $\alpha \in \theta \}$.  But $N_G(A) \cap P =
M$.  Thus $w \in W_M(P_0,A_0)$.  \qed

In order to prove Theorem 3.1 we extend the proof of Silberger in\Lspace \Lcitemark 13\LIcitemark{}, \S 2.5\RIcitemark \Rcitemark \Rspace{} to the
disconnected case.  Fortunately many of the technical results on intertwining forms needed
are proven in\Lspace \Lcitemark 13\LIcitemark{},  \S 1\RIcitemark \Rcitemark \Rspace{} for any totally disconnected group and so can be
directly used in our case.  We follow Silberger's notation.   Let $(P_1,A_1)$ and
$(P_2,A_2)$ be cuspidal parabolic pairs in $G$.  We may as well assume that they are
standard with respect to a fixed minimal parabolic pair $(P_0,A_0)$.

We need to study the orbits for the action of $P_1 \times P_2$ on $G$ given by
$y \cdot (p_1,p_2) = p_1^{-1}yp_2$.  Recall that $G^0 = \cup _v P_0^0vP^0_0$ where $v$
runs over $W(G^0,A_0) = N_{G^0}(A_0)/M^0_0$.  Thus,  $G = \cup _w wG^0 =
\cup _{w,v} wP^0_0vP_0^0 = \cup _{w,v}
P^0_0wvP_0^0 $ where $w \in W_G(P_0,A_0), v \in W(G^0,A_0)$.  Now since $P_0^0 \subset
P_i, i = 1,2$, each double coset $P_1yP_2$ can be represented by $y=wv,
w \in W_G(P_0,A_0), v \in W(G^0,A_0)$.  Write $W_i = W_{M_i}(P_0,A_0), i = 1,2$.

\proclaim Lemma 3.9.  Let $\calO = P_1w_0vP_2$ be an orbit of $P_1 \times P_2 $ in $G$
where $w_0 \in W_G(P_0,A_0)$ and $v \in W(G^0,A_0)$.  Then for $w \in W_G(P_0,A_0)$, 
$\calO \cap wG^0 $ is empty unless $w \in W_1w_0W_2$.   If
$w \in W_1w_0W_2$, then $\calO \cap wG^0  = w \calO _w$ where
$\calO _w$ is a finite union of orbits of $(w^{-1}P_1^0w) \times P_2^0$ in $G^0$, all of
which have fixed dimension $d_{\calO }$ equal to the dimension of the orbit 
$\calO _0 = (w_0^{-1}P_1^0w_0)vP_2^0$.  \endproclaim

\proof  Using Lemma 3.8 we can write 
$$P_1w_0vP_2 = \cup _{w_1,w_2} w_1P_1^0w_0vw_2P_2^0$$ where $w_i \in W_i =
W_{M_i}(P_0,A_0), i = 1,2$. 
But $$ w_1P_1^0w_0vw_2P_2^0 = w (w^{-1}P_1^0w)w_2^{-1}vw_2P_2^0$$ where  $w = w_1w_0w_2
\in W_G(P_0,A_0)$.  Thus $ w_1P_1^0w_0vw_2P_2^0 \subset w_1w_0w_2G^0$ and
for fixed $w \in W_1w_0W_2$, 
$$P_1w_0vP_2 \cap wG^0 = w \cup _{w_2} (w^{-1}P_1^0w)w_2^{-1}vw_2P_2^0$$ where $w_2$ runs
over elements of $W_2$ such that $ww_2^{-1}w_0^{-1} \in W_1$.  Finally, 
$$(w^{-1}P_1^0w)w_2^{-1}vw_2P_2^0 =  w_2^{-1}(w_0^{-1}P_1^0w_0)v P_2^0 w_2$$ so that
$$\dim (w^{-1}P_1^0w)w_2^{-1}vw_2P_2^0 = \dim (w_0^{-1}P_1^0w_0)v P_2^0 = d_{\calO}.$$
\qed

Because of Lemma 3.9 each orbit $\calO$ has a well-defined dimension $d_{\calO}$.  For
each integer $\nu \geq 0$, let $\calO (\nu )$ denote the union of all orbits of dimension
less than or equal to $\nu $.  Set $\calO (-1) = \emptyset $.  

\proclaim Lemma 3.10.  For each $\nu \geq 0$, $\calO (\nu )$ is a closed set in the
$p$-adic topology.  Further, if $\calO $ is an orbit of dimension $d$, then $\calO \cup
\calO (d-1)$ is also closed in the $p$-adic topology.  \endproclaim

\proof  Using Lemma 3.9 we see that $\calO (\nu ) \cap wG^0 = w\calO _w (\nu )$ where
$\calO _w (\nu )$ is the union of all orbits of $(w^{-1}P_1^0w) \times P_2^0$ in $G^0$
having dimension less than or equal to $\nu $.  This set is closed in the $p$-adic
topology by\Lspace \Lcitemark 13\LIcitemark{}, pg. 93\RIcitemark \Rcitemark \Rspace{}.  Thus $\calO (\nu )$ is a finite union of closed sets, hence is closed.
Similarly, if $\calO $ has dimension $d$, then $[\calO \cup \calO (d-1) ]\cap wG^0 =
w (\calO _w \cup \calO _w (d-1))$ where $\calO _w \cup \calO _w (d-1)$ is a finite union
of sets of the form $\calO ' \cup \calO _w(d-1)$ where $\calO '$ is an orbit of
$(w^{-1}P_1^0w) \times P_2^0$ in $G^0$ of dimension $d$.  These are also closed by\Lspace \Lcitemark 13\LIcitemark{}, page 94\RIcitemark \Rcitemark \Rspace{}. 
\qed

Let $E = V_1 \otimes V_2$ where $V_i$ is the space on which $\sigma _i$ acts, $i = 1,2$.
Let $\calT  $ be the space of all $E$-distributions $T^0$ on $G$ such that
$$T^0 (\lambda (p_1)\rho (p_2)\alpha ) = T^0(\delta _1(p_1)^{1\over 2}
\delta _2(p_2)^{1\over 2}\sigma _1 (p_1)^{-1} \otimes \sigma _2 (p_2)^{-1}\alpha )$$
for all $(p_1, p_2) \in P_1 \times P_2 , \alpha \in C_c^{\infty } (G:E)$.  
Here $\lambda $ and $\rho $ denote left and right translations respectively and
$\delta _i$ is the modular function of $P_i, i = 1,2$
The first step in the proof of Theorem 3.1 is the inequality\Lspace \Lcitemark 13\LIcitemark{},  1.9.4\RIcitemark \Rcitemark \Rspace{}
$$I(\pi _1, \pi _2) \leq \dim \calT .\eqno(3.1)$$
Here $I(\pi _1, \pi _2)$ is the dimension of the space of
``intertwining forms'' defined in\Lspace \Lcitemark 13\LIcitemark{},  \S 1.6\RIcitemark \Rcitemark \Rspace{}.  It is related to the dimension of the space
of intertwining operators by
 $I(\pi _1, \pi _2) = J(\tilde {\pi } _1, \pi _2)$ where $\tilde {\pi } _1$ is the
contragredient of $\pi _1$.\Lspace \Lcitemark 13\LIcitemark{},  1.6.2\RIcitemark \Rcitemark \Rspace{} 

If $\calO $ is
an orbit of dimension $d$, write $\calT (\calO )$ for the vector space of $T^0 \in \calT $
with support in $\calO \cup \calO (d-1)$ and $\calT _{\nu }$ for the space of those with
support in $\calO (\nu )$.  We have
$$\calT = \sum _{\calO } \calT (\calO ) $$ and
$$\dim \calT \leq
\sum _{\calO } \dim( \calT (\calO )/ \calT _{d(\calO ) -1}).\eqno(3.2)$$

\proclaim Lemma 3.11.  Suppose that $(P_1,A_1)$ and $(P_2,A_2)$ are semi-standard
cuspidal parabolic pairs in $G$.  Then $M_1 \cap P_2 = (M_1 \cap M_2)(M_1 \cap N_2) =
\ ^*P_1$ is a cuspidal parabolic subgroup of $M_1$ with split component $\ ^*A_1 = A_1A_2$.
\endproclaim

\proof  We know from the connected case\Lspace \Lcitemark 13\LIcitemark{}, p 94\RIcitemark \Rcitemark \Rspace{} that $\ ^*P_1^0 = M_1^0 \cap P_2^0$ is a
parabolic subgroup of $M_1^0$ with split component $\ ^*A_1 = A_1A_2$ and Levi
decomposition  $\ ^*P_1^0 = (M_1^0 \cap M_2^0)(M_1^0 \cap N_2)$.  Thus there is a cuspidal
parabolic subgroup $\ ^*P_1$ of $M_1$ with split component $\ ^*A_1$ and Levi decomposition
$\ ^*P_1= \ ^*M_1\  ^*N_1$.  Here $\ ^*M_1 = C_{M_1}(\ ^*A_1) =  C_{M_1}(A_1A_2) = M_1 \cap
C_G(A_2) = M_1 \cap M_2$ and $\ ^*N_1 = M_1^0 \cap N_2 =  M_1 \cap N_2$.  Clearly $\ ^*P_1
=  (M_1 \cap M_2)(M_1 \cap N_2) \subset M_1 \cap P_2$.  
Thus we need only show
that $ M_1 \cap P_2 \subset \ ^*P_1$. 

 Let $x \in M_1 \cap P_2$.  Using the Levi
decomposition of $P_2$ we can write $x = m_2n_2$ where $m_2 \in M_2, n_2 \in N_2$.   Since
$m_2 n_2 \in M_1= C_G(A_1)$ we have $m_2n_2a_1n_2^{-1}m_2^{-1} = a_1$ for any $a_1 \in
A_1$.  This implies that $a_1^{-1}n_2a_1n_2^{-1} = a_1^{-1}m_2^{-1}a_1m_2$.  But since 
$(P_1,A_1)$ and $(P_2,A_2)$ are semi-standard, we have $A_1 \subset A_0 \subset C_G(A_2) =
M_2$.  Thus $a_1^{-1}n_2a_1n_2^{-1} \in N_2$ and $a_1^{-1}m_2^{-1}a_1m_2 \in M_2$.  Hence
$a_1^{-1}n_2a_1n_2^{-1} = a_1^{-1}m_2^{-1}a_1m_2 \in N_2 \cap M_2 = \{ 1 \}$.  Thus $n_2$
and $m_2$ both commute with $a_1$ so that $n_2 \in M_1 \cap N_2$ and $m_2 \in M_1 \cap
M_2$. \qed

Return to the assumption that $(P_1,A_1)$ and $(P_2,A_2)$
are standard with respect to $(P_0,A_0)$ and that $y \in W(A_0)$ so that $(P_2^y,A_2^y)$
is semi-standard.  Using Lemma 3.11 we know that $\ ^*P_1 = M_1 \cap P_2^y$ and $\ ^*P_2 =
M_2
\cap P_1^{y^{-1}}$ are cuspidal parabolic subgroups of $M_1$ and $M_2$ respectively.  

 Let $E'(y)$ denote the space of linear functionals $\phi $ on $E = V_1\otimes
V_2$ such that $$\delta _1(p)^{1\over 2} \delta _2(p^{y^{-1}})^{1\over 2}<\phi , \sigma
_1(p) v_1 \otimes \sigma _2 (p^{y^{-1}})v_2> = \delta _{P_1 \cap P_2^y} (p) <\phi , v_1
\otimes v_2>$$ for all $p \in P_1 \cap P_2^y, v_i \in V_i, i = 1,2$.  

\proclaim Lemma 3.12.  Let $m \in \ ^*M_1 = M_1 \cap M_2^y$.  Then
$$\delta _{P_1 \cap P_2^y} (m) = \delta _{*1}(m)^{1\over 2}
\delta _{*2}(m^{y^{-1}})^{1\over 2} \delta _1(m)^{1\over 2}
\delta _2(m^{y^{-1}})^{1\over 2}$$
where $\delta _{*i}$ denotes the modular function for $\ ^*P_i$ and $\delta _i$ the
modular function for $P_i, i = 1,2$.  \endproclaim 

\proof  Define the homomorphism $\delta : \ ^*M_1 \rightarrow \bfR _+^*$ by
$$\delta (m) =  \delta _{*1}(m)^{1\over 2}
\delta _{*2}(m^{y^{-1}})^{1\over 2} \delta _1(m)^{1\over 2}
\delta _2(m^{y^{-1}})^{1\over 2}\delta _{P_1 \cap P_2^y} (m)^{-1}, m \in \ ^*M_1.$$
By\Lspace \Lcitemark 13\LIcitemark{}, \S 2.5.2\RIcitemark \Rcitemark \Rspace{}, the restriction of $\delta $ to $\ ^*M_1  \cap G^0$ is trivial.  Now since 
the quotient of $\ ^*M_1$ by  $\ ^*M_1  \cap G^0$ is a finite group, $\delta = 1$ on all of 
$\ ^*M_1 $.  \qed

The proof of the following lemma is the same as that of Lemma 2.5.1 in\Lspace \Lcitemark 13\Rcitemark \Rspace{}.

\proclaim Lemma 3.13.  If $\phi \in E'(y)$, then $\phi $ vanishes on $V_1(\ ^*P_1)\otimes
V_2 + V_1 \otimes V_2(\ ^*P_2)$.  Let $m \in \ ^*M_1, v_1 \otimes v_2 \in E$, and $\phi \in
E'(y)$.  Then
$$<\phi , \sigma _1 (m) v_1 \otimes \sigma _2 (m^{y^{-1}})v_2> =
\delta _{*1}(m)^{1\over 2} \delta _{*2}(m^{y^{-1}})^{1\over 2} <\phi , v_1 \otimes v_2>.$$
\endproclaim

\proclaim Corollary 3.14.  If $A_1 = A_2^y$, then
$$<\phi , \sigma _1 (m) v_1 \otimes \sigma _2 (m^{y^{-1}})v_2> =
<\phi , v_1 \otimes v_2>$$  for all $m \in \ ^*M_1, v_1 \otimes v_2 \in E$, and $\phi \in
E'(y)$. \endproclaim

\proof  As in\Lspace \Lcitemark 13\LIcitemark{}, \S 2.5.4\RIcitemark \Rcitemark \Rspace{}, $A_1 = A_2^y$ implies that $\ ^*N_1 = \ ^*N_2 = \{ 1 \}$ so that
$\ ^*P_1$ and $\ ^*P_2$ are reductive and $\delta _{*1} = \delta _{*2} = 1$.  \qed

\proclaim Corollary 3.15.  Assume that $\sigma _1$ and $\sigma _2$ are either both
supercuspidal or both discrete series.  Then $E'(y) \not = \{ 0 \}$ implies that $A_1 =
A_2^y$. \endproclaim

\proof   This follows as in\Lspace \Lcitemark 13\LIcitemark{}, \S\S 2.5.3, 2.5.5\RIcitemark \Rcitemark \Rspace{} using Proposition 2.17 and Lemma 3.6.  \qed

The following Lemma is now proven as in\Lspace \Lcitemark 13\LIcitemark{}, \S 2.5.7\RIcitemark \Rcitemark \Rspace{}.

\proclaim Lemma 3.16.  Assume that $\sigma _1$ and $\sigma _2$ are either both
supercuspidal or both discrete series.  Let $\calO $ be an orbit in $G$ of dimension $d$. 
Then $\dim \calT (\calO ) / \calT _{d-1} = 0$ unless there exists $y \in \calO $ such
that $A_1 = A_2^y$.  If $A_1 = A_2^y$ for some $y \in \calO $, then
 $$\dim (\calT (\calO ) / \calT _{d-1} ) \leq I(\sigma _1, \sigma _2^y).$$  \endproclaim

\medbreak
\noindent {\bf Proof of Theorem 3.1.}  Combine equations (3.1) and (3.2) to
obtain
$$I(\pi _1, \pi _2) \leq \dim \calT \leq 
\sum _{\calO } \dim( \calT (\calO )/ \calT _{d(\calO ) -1}).$$  Now using Lemma 3.16, if
$A_1$ and $A_2$ are not conjugate, $ \dim( \calT (\calO )/ \calT _{d(\calO ) -1}) = 0$ for
every orbit $\calO$ so that $I(\pi _1, \pi _2 ) = 0$.  If $A_1 $ and $A_2$ are conjugate we
have, again using Lemma 3.16,
$$I(\pi _1, \pi _2) \leq \dim \calT \leq \sum _{s \in W} 
 \dim( \calT (\calO _s)/ \calT _{d(\calO _s) -1}) \leq \sum _{s \in W} I(\sigma _1, \sigma
_2 ^{y_s}).$$
To complete the proof note that $J(\pi _1, \pi _2) = I(\tilde {\pi } _1, \pi
_2), J(\sigma _1, \sigma
_2 ^{y_s}) = I(\tilde {\sigma} _1, \sigma
_2 ^{y_s})$ and $\tilde {\pi } _1 = \Ind_{P_1}^G(\tilde {\sigma} _1)$.  Thus the statement
with dimensions of spaces of intertwining operators rather than intertwining forms follows
by substituting $\tilde {\pi _1}$ for $\pi _1$.  
 \qed

\section{$R$-groups for $\protect\Ind_{P^0}^G(\sigma _0)$}

In this section we will study representations of $G$ which are induced from a parabolic
subgroup $P^0 = M^0N$ of $G^0$.  Because in this section we will only work with parabolic subgroups
of
$G^0$, we will simplify notation by dropping the superscripts on $P^0$ and $M^0$.  

Let $P= MN$ be a parabolic subgroup of $G^0$ and fix an irreducible discrete series
representation
$(\sigma , V)$ of $M$.
Define 
$$\thickmuskip=.5\thickmuskip\calH _P(\sigma ) = \{ f\in C^{\infty }(G,V): f(xmn) = \delta _P^{-{1\over 2}}(m)\sigma
(m)^{-1} f(x) {\rm \ for \ all \ } x \in G, m \in M, n \in N \}$$ where $\delta _P$ is the
modular function on $P$. Then $G$ acts by left translation on $\calH _P(\sigma )$ and we call this
induced representation
$I_P(\sigma )$. We will also need to consider 
$\calH _P^0 (\sigma ) = $ $$ \{ \phi \in C^{\infty }(G^0,V): \phi (xmn) = \delta _P^{-{1\over
2}}(m)\sigma (m)^{-1} \phi (x) {\rm \ for \ all \ } x \in G^0, m \in M, n \in N \}.$$
$G^0$ acts by left
translation on $\calH ^0_P(\sigma )$ and we call this induced representation $I^0_P(\sigma )$. 

It is well known that the equivalence class of $I^0_P(\sigma ) = \Ind _P^{G^0}(\sigma )$ is
independent of
$P$.  But
$$I_P(\sigma ) = \Ind_P^G(\sigma ) \simeq \Ind _{G^0}^G \Ind _P^{G^0}(\sigma ),$$  so that the
equivalence class of $I_P(\sigma )$ is also independent of $P$.
We denote the equivalence classes of $I_P(\sigma )$ and $I_P^0(\sigma )$ by
$i_{G,M}(\sigma )$ and $i_{G^0,M}(\sigma )$ respectively.  If $\pi $ is a representation of $G^0$ we
will also write $i_{G,G^0}(\pi )$ for the equivalence class of the induced representation
$\Ind_{G^0}^G(\pi )$.

We first want to compare the dimensions of the intertwining algebras of $i_{G,M}(\sigma )$  and
$i_{G^0,M}(\sigma )$ .  For this we need the results of Gelbart and Knapp summarized in
\S 2 and the following
facts.

\proclaim Lemma 4.1.  Suppose $\pi _1$ and $\pi _2$ are irreducible representations of $G^0$. 
Then $i_{G,G^0}(\pi _1 )$  and $i_{G,G^0}(\pi _2)$ have an irreducible constituent in common if and
only if $\pi _2 \simeq \pi _1^g$ for some $g \in G$.  In this case they are equivalent. 
\endproclaim

\proof  For $i = 1,2$, write $$i_{G,G^0}(\pi _i ) = r_i \sum _{\chi \in X/X(\Pi _i)} \Pi _i\otimes
\chi $$ as in Lemma 2.13.  Suppose that $i_{G,G^0}(\pi _1 )$  and $i_{G,G^0}(\pi _2)$ have an
irreducible  constituent in common.  Then $\Pi _1 \otimes \chi _1 \simeq \Pi _2 \otimes \chi _2$ for
some $\chi _1, \chi _2 \in X$.  Now $$\pi _2 \subset \Pi _2 \otimes \chi _2 |_{G^0} \simeq
\Pi _1 \otimes \chi _1 |_{G^0} \subset i_{G,G^0}(\pi _1)|_{G^0} \simeq \sum _{g \in G/G^0} \pi
_1^g.$$  Thus $\pi _2 \simeq \pi _1^g$ for some $g \in G$.  Conversely, if
 $\pi _2 \simeq \pi _1^g$ for some $g \in G$, then clearly
$$i_{G,G^0}(\pi _2) \simeq i_{G,G^0}(\pi _1^g) \simeq i_{G,G^0}(\pi _1).$$   \qed

\proclaim Lemma 4.2.  Suppose that for some $g \in G$, both $\pi $ and $\pi ^g$ are irreducible
constituents of $i_{G^0,M}(\sigma )$.  Then there is $x_0 \in G^0$ such that $gx_0 \in N_G(\sigma
) = \{ x \in N_G(M): \sigma ^x \simeq \sigma \}$.  Conversely, if 
$\pi $ is an irreducible
constituent of $i_{G^0,M}(\sigma )$ and if $g \in N_G(\sigma
)G^0$, then $\pi ^g $ is also an irreducible constituent of $i_{G^0,M}(\sigma )$ and
$\pi $ and $\pi ^g$ occur with the same multiplicities.  \endproclaim

\proof  Suppose that $\pi , \pi ^g \subset i_{G^0,M}(\sigma )$. Then
since $\pi ^g \subset i_{G^0, M^g}(\sigma ^g)$, we see that $ i_{G^0,M}(\sigma )$ and
$ i_{G^0,M^g}(\sigma ^g)$ have an irreducible constituent in common. Thus there is $x_0 \in G^0$
such that $M^{gx_0} = M$ and $\sigma ^{gx_0} \simeq \sigma $, i.e. $gx_0 \in N_G(\sigma )$.  
Conversely, if $g \in N_G(\sigma )G^0$, then the multiplicity of $\pi $ in $i_{G^0,M}(\sigma )$ is
equal to  the multiplicity of $\pi ^g $ in $i_{G^0,M}(\sigma )^g \simeq  i_{G^0,M^g}(\sigma ^g)
\simeq i_{G^0,M}(\sigma ) $.  \qed

\proclaim Lemma 4.3.  Let $\pi $ be an irreducible constituent of $i_{G^0,M}(\sigma )$ and let
$G_{\pi } = \{ x \in G: \pi ^x \simeq \pi \}$.
 Then
$$[G_{\pi }/G^0] = [N_{G_{\pi }}(\sigma )/N_{G^0}(\sigma )]$$ where
$N_{G^0}(\sigma ) = N_G(\sigma ) \cap G^0$ and $N_{G_{\pi }}(\sigma ) = 
N_G(\sigma ) \cap G_{\pi }$.  \endproclaim

\proof  Consider the mapping from $N_{G_{\pi }}(\sigma )$ to $G_{\pi }/G^0$ given by
$g \mapsto gG^0$.  Its kernel is $N_{G_{\pi }}(\sigma ) \cap G^0 = N_{G^0}(\sigma )$.
Further, given $g \in G_{\pi }$, $\pi ^g \simeq \pi $ occurs in $i_{G^0,M}(\sigma )$, so by Lemma
4.1 there is $x_0 \in G_0$ such that $gx_0 \in N_G(\sigma )$.   But $\pi ^{gx_0} \simeq \pi ^g
\simeq
\pi $ so that $gx_0 \in G_{\pi } \cap  N_G(\sigma ) =  N_{G_{\pi }}(\sigma )$ and
$gx_0G^0 = gG^0$.  Thus the mapping is surjective.  \qed

Let $C(\sigma )$ denote the algebra of $G$-intertwining operators for
$i_{G,M}(\sigma )$ and let $C^0(\sigma )$ denote the algebra of $G^0$-intertwining
operators for $i_{G^0,M}(\sigma )$.

\proclaim Lemma 4.4.  $\dim C(\sigma ) = \dim C^0(\sigma ) [N_G(\sigma )/N_{G^0}(\sigma )]$.
\endproclaim 

\proof  Let $$i_{G^0,M}(\sigma ) = \sum _{\pi \in S(\sigma )}m_{\pi } \pi $$ be the decomposition of
$i_{G^0,M}(\sigma )$ into irreducible constituents.  For $\pi _1, \pi _2 \in S(\sigma )$, we will
say that $\pi _1 \sim \pi _2$ if there is $g \in G$ such that $\pi _2 \simeq \pi _1^g$.  Then using
Lemma 4.2, we have $\pi _1 \sim \pi _2$ if and only if $\pi _2 \simeq \pi _1^g$ for some $g \in
N_G(\sigma )$.   We can write
$$i_{G^0,M}(\sigma ) = \sum _{\pi \in S(\sigma )/\sim }m_{\pi }
 \sum _{ g \in N_G(\sigma )/N_{G_{\pi }}(\sigma )} \pi ^g .$$

Now $$i_{G,M}(\sigma ) \simeq i_{G,G^0}(i_{G^0,M}(\sigma )) \simeq
 \sum _{\pi \in S(\sigma )/\sim }m_{\pi }
 \sum _{ g \in N_G(\sigma )/N_{G_{\pi }}(\sigma )}i_{G,G^0}( \pi ^g )$$
$$ = \sum _{\pi \in S(\sigma )/\sim }m_{\pi }
 [ N_G(\sigma )/N_{G_{\pi }}(\sigma )] r_{\pi } \sum _{\chi \in X/X(\Pi _{\pi } )} \Pi
_{\pi }\otimes \chi $$ where $\Pi _{\pi }$ is an irreducible representation of $G$ such that
$\pi \subset \Pi _{\pi }|_{G^0}$.  Because of Lemma 4.1, the representations $ \Pi
_{\pi }\otimes \chi $ are pairwise inequivalent as $\pi $ ranges over $S(\sigma )/\sim$ and $\chi $
ranges over $ X/X(\Pi _{\pi } )$.  Thus
$$\dim C(\sigma ) = \sum _{\pi \in S(\sigma )/\sim }m^2_{\pi }
 [ N_G(\sigma )/N_{G_{\pi }}(\sigma )] ^2 r^2_{\pi } [ X/X(\Pi _{\pi } )].$$
But by Lemmas 2.13 and 4.3, $$r^2_{\pi } [ X/X(\Pi _{\pi } )] = [G_{\pi }/G^0] = [ N_{G_{\pi
}}(\sigma )/N_{G^0}(\sigma )].$$  Thus
$$\dim C(\sigma ) = \sum _{\pi \in S(\sigma )/\sim }m^2_{\pi }
 [ N_G(\sigma )/N_{G_{\pi }}(\sigma )] ^2 [ N_{G_{\pi }}(
\sigma )/N_{G^0}(\sigma )]$$
$$ =   [ N_G(
\sigma )/N_{G^0}(\sigma )]\sum _{\pi \in S(\sigma )/\sim }m^2_{\pi }
 [ N_G(\sigma )/N_{G_{\pi }}(\sigma )] $$ $$ =  [ N_G(
\sigma )/N_{G^0}(\sigma )] \dim C^0(\sigma ).$$  \qed

We now want to find a basis for $C(\sigma )$.  We proceed as in the connected case. 
Let $A$ be the split component of the center of $M$ and write $\fraka _{\bfC }^*$ for the
dual of its complex Lie algebra.  Each $\nu \in \fraka _{\bfC }^*$ determines a one-dimensional
character $\chi _{\nu }$ of $M$ which is defined by
$$\chi _{\nu }(m) = q^{<H_P(m), \nu >}, m \in M.$$
  We write $(I_P(\sigma : \nu ),\calH
_P(\sigma :\nu ))$ and $(I^0_P(\sigma : \nu ),\calH ^0_P(\sigma
:\nu ))$
for the induced representations of $G$ and $G^0$ as above corresponding to $\sigma _{\nu } = \sigma
\otimes \chi _{\nu } $.  Let $K$ be a good maximal compact subgroup of $G^0$ with respect to a
minimal parabolic pair $(P_0,A_0)$ of $G^0$ such that $P_0 \subset P, A \subseteq A_0$.  Then $G^0 =
KP$ and we also have the usual compact realization of $I_P^0(\sigma )$ on
$\calH _P^K(\sigma ) = $ $$
 \{ f_K\in C^{\infty }(K,V): f_K(kmn) = \sigma ^{-1}(m) f_K(k) {\rm \ for
\ all \  } k \in K, m \in M \cap K, n \in N \cap K  \}.$$
The intertwining operators between $\calH _P^0(\sigma , \nu )$ and $\calH _P^K(\sigma )$ are given by
$$F_P^K(\nu ): \calH _P^0(\sigma :\nu ) \rightarrow \calH _P^K(\sigma ),$$ 
$$F_P^K(\nu ) \phi (k) = \phi (k), \phi \in \calH _P^0(\sigma :\nu ) , k \in K.$$
For all $f_K \in \calH _P^K(\sigma ), x \in G^0$,
$$F_P^K(\nu )^{-1} f_K(x) = \delta _P ^{-{1\over 2}}(\mu (x))\sigma ^{-1}_{\nu }(\mu (x)) 
f_K(\kappa (x)).$$ Here for any $x \in G^0, \kappa (x) \in K, \mu (x) \in M$ are chosen so that $x
\in
\kappa (x) \mu (x) N$.  

Since we don't know whether there is a ``good'' maximal compact subgroup for $G$, we don't have a
single compact realization for $I_P(\sigma )$.  However we can proceed one coset at a time as
follows.  Write $G$ as a disjoint union of cosets, $G = \cup _{i=1}^k x_iG^0$.  For
any
$f
\in
\calH _P(\sigma ,
\nu ), 1
\leq i
\leq k$, we can define 
$$f _i(x) = \cases {f (x), & if $x \in x_iG^0$;\cr 0, & otherwise.\cr}$$
Then $f _i \in \calH _P(\sigma :\nu )$ for each $i$ and $f = \sum _{i=1}^k f _i$. 
For each $1 \leq i \leq k$, define 
$$F_P^i(\nu ): \calH _P(\sigma :\nu ) \rightarrow \calH _P^K(\sigma )$$ by
$$F_P^i(\nu ) f (k) = f (x_ik), f \in \calH _P(\sigma :\nu ) , k \in K.$$  Define 
$$F_P^i(\nu )^{-1}: \calH _P^K(\sigma ) \rightarrow \calH _P(\sigma :\nu )$$ by
$$F_P^i(\nu )^{-1} f_K(x) = \cases{ \delta _P ^{-{1\over 2}}(\mu (x_0))\sigma ^{-1}_{\nu }(\mu (x_0))
 f_K(\kappa (x_0)), & if $x = x_i x_0, x_0 \in G^0$;\cr 0, & otherwise.\cr}$$
Then $F_P^i(\nu ) F_P^i(\nu )^{-1}f_K = f_K$ for all $f_K \in \calH _P^K(\sigma )$ and
$ F_P^i(\nu )^{-1}F_P^i(\nu ) f = f _i$ for all $f \in \calH _P(\sigma :\nu )$.

If $P = MN$ and $P'=MN'$ are two parabolic subgroups of $G^0$ with Levi component $M$, then we have
the formal intertwining operators $$J^0(P':P:\sigma :\nu ) :\calH ^0 _P(\sigma : \nu ) \rightarrow
\calH ^0_{P'}(\sigma : \nu )$$ given by the standard integral formula
$$J^0(P':P:\sigma :\nu )\phi (x) = \int _{\barN \cap N'} \phi (x\barn )d \barn, \ \ x \in G^0.$$
Here $M \barN$ is the opposite
parabolic to $P$ and $d\barn $ is normalized Haar measure on $\barN \cap N'$.  We can define
 $$J(P':P:\sigma :\nu ) :\calH _P(\sigma : \nu ) \rightarrow
\calH _{P'}(\sigma : \nu )$$ by the same formula
$$J(P':P:\sigma :\nu )f(x) = \int _{\barN \cap N'} f(x\barn )d \barn, \ \ x \in G.$$  
In order to talk rigorously about holomorphicity and analytic continuation of these operators we
transfer them to the compact realizations.  Thus we have
$$J_K^0(P':P:\sigma :\nu ) :\calH _P^K(\sigma ) \rightarrow \calH _P^K(\sigma ) $$ defined for
$f_K \in \calH _P^K(\sigma ), k \in K$, by 
$$J_K^0(P':P:\sigma :\nu )f_K(k) = F_{P'}^K(\nu
)J^0(P':P:\sigma :\nu )F_P^K(\nu )^{-1}f_K(k) =$$ $$\int _{\barN \cap N'} \delta _P ^{-{1\over
2}}(\mu (k\barn ))\sigma ^{-1}_{\nu }(\mu (k\barn ))  f_K(\kappa (k\barn ))
d\barn .$$ 

\proclaim Theorem 4.5. (Harish-Chandra\Lspace \Lcitemark 7\Citecomma
13\Rcitemark \Rspace{})  Suppose $\sigma $ is an irreducible discrete series
representation of $M$.  There is a chamber
$\fraka ^* _{\bfC }(P':P)$ in $\fraka _{\bfC }^*$ such that for $\nu \in \fraka ^*_{\bfC }(P':P)$
the formal intertwining operator $J^0(P':P:\sigma :\nu )$ converges and defines a bounded operator.
For a fixed $f_K \in \calH _P^K(\sigma )$, the mapping
$$\nu \mapsto J_K^0(P':P:\sigma :\nu )f_K$$ from $\fraka ^* _{\bfC }(P':P)$ to $\calH _P^K(\sigma )$
is holomorphic.  Further, it extends to a meromorphic function on all of $\fraka _{\bfC }^*$. 
\endproclaim  

Because of Harish-Chandra's theorem we can define $J^0(P':P:\sigma :\nu )$ for all $\nu \in \fraka
^*_{\bfC }$ by
$$J^0(P':P:\sigma :\nu ) =  F_{P'}^K(\nu
)^{-1} J_K^0(P':P:\sigma :\nu )F_P^K(\nu ).$$

\proclaim Corollary 4.6. Suppose $\sigma $ is an irreducible discrete series
representation of $M$.  Then
the formal intertwining operator $J(P':P:\sigma :\nu )$ converges and defines a bounded operator
for $\nu \in \fraka ^*_{\bfC }(P':P)$.  Further, for all $\nu \in \fraka ^*_{\bfC }(P':P)$, we have
$$J(P':P:\sigma :\nu )= \sum _{i=1}^k  F_{P'}^i(\nu
)^{-1} J_K^0(P':P:\sigma :\nu )F_P^i(\nu ) .$$
\endproclaim 

\proof  It is clear from the definitions that for $\nu \in \fraka ^*_{\bfC }(P':P),
f \in \calH _P(\sigma : \nu ), x \in G, x_0 \in G^0$,
$$J(P':P:\sigma :\nu )f(xx_0) = J^0(P':P:\sigma :\nu )\phi (x_0)$$ where $\phi = 
l(x^{-1})f|_{G^0}$ is the restriction to $G^0$ of the left translate of $f$ by $x^{-1}$.  Thus the
integral converges.  Fix
$1
\leq i
\leq k$.  An elementary calculation shows that 
$$ F_{P'}^i(\nu
)J(P':P:\sigma :\nu )F_P^i(\nu )^{-1} = J_K^0(P':P:\sigma :\nu ).$$ 
Thus for any $f \in \calH _P(\sigma )$, since clearly $J(P':P:\sigma :\nu )(f _i) =
(J(P':P:\sigma :\nu )f )_i$ we have
 $$F_{P'}^i(\nu
)^{-1}J_K^0(P':P:\sigma :\nu )F_P^i(\nu ) f $$ 
$$= F_{P'}^i(\nu
)^{-1} F_{P'}^i(\nu
)J(P':P:\sigma :\nu )F_P^i(\nu )^{-1} F_P^i(\nu ) \phi = (J(P':P:\sigma :\nu )f )_i.$$ \qed

Because of Corollary 4.6 we can define $J(P':P:\sigma :\nu )$ for all $\nu \in \fraka
^*_{\bfC }$ by the formula
$$J(P':P:\sigma :\nu ) = \sum _{i=1}^k  F_{P'}^i(\nu
)^{-1} J_K^0(P':P:\sigma :\nu )F_P^i(\nu ) .$$  

\proclaim Corollary 4.7.  Let $\nu \in \fraka ^*_{\bfC }, f \in \calH _P(\sigma : \nu ), x \in G, x_0
\in G^0$.  Then $$J(P':P:\sigma :\nu )f(xx_0) = J^0(P':P:\sigma :\nu )\phi (x_0)$$ where $\phi = 
l(x^{-1})f|_{G^0}$.  \endproclaim

\proof  Suppose that $1 \leq j \leq k$.  
Because $J^0(P':P:\sigma :\nu )$ commutes with the left action of $G^0$, it is enough to prove
the result for $x = x_j$.  For any  $\nu \in \fraka ^*_{\bfC }, f \in \calH _P(\sigma : \nu ), x_0
\in G^0, J(P':P:\sigma :\nu )f(x_jx_0) =$ $$\sum _{i=1}^k  F_{P'}^i(\nu
)^{-1} J_K^0(P':P:\sigma :\nu )F_P^i(\nu )f(x_jx_0) =   F_{P'}^j(\nu
)^{-1} J_K^0(P':P:\sigma :\nu )F_P^j(\nu )f(x_jx_0) $$
$$ =  F_{P'}^K(\nu
)^{-1} J_K^0(P':P:\sigma :\nu )F_P^K(\nu )\phi _j(x_0)$$
 where
$\phi _j= 
L(x_j^{-1})f|_{G^0}$.   \qed

Fix scalar normalizing factors $r(P':P:\sigma :\nu )$ as in\Lspace \Lcitemark 2\Rcitemark \Rspace{}
used to define the normalized intertwining operators 
$$R^0(P':P:\sigma ) = r(P':P:\sigma :0)^{-1}J^0(P':P:\sigma :0).$$  Using Corollary 4.7, the fact
that 
$r(P':P:\sigma :\nu )^{-1}J^0(P':P:\sigma :\nu)$ is holomorphic and non-zero at $\nu = 0$ will imply
that $ r(P':P:\sigma :\nu )^{-1}J(P':P:\sigma :\nu)$ is also holomorphic and non-zero at $\nu = 0$. 
Thus we can use the same normalizing factors to define $$R(P':P:\sigma  ) = r(P':P:\sigma :0
)^{-1}J(P':P:\sigma :0 ).$$   For the normalized intertwining operators we will also have the
formula
$$R(P':P:\sigma )f(xx_0) = R^0(P':P:\sigma )\phi (x_0)$$ where notation is as in Corollary 4.7.

\proclaim Lemma 4.8.  Suppose $P_1, P_2$, and $P_3$ are parabolic subgroups of $G^0$ with Levi
component $M$.  Then $R(P_1:P_3:\sigma ) = R(P_1:P_2:\sigma )R(P_2:P_3:\sigma  )$.  \endproclaim

\proof  This follows easily using Corollary 4.7 from the corresponding formula for the connected
case.  \qed

\proclaim Lemma 4.9.  Suppose $P_1$ and $P_2$ are parabolic subgroups of $G^0$ with Levi component
$M$.  Then $R(P_2:P_1:\sigma )$ is an equivalence from $\calH _{P_1}(\sigma )$ onto
$\calH _{P_2}(\sigma )$.  \endproclaim

\proof  It follows from Lemma 4.8 that
$$R(P_1:P_2:\sigma )R(P_2:P_1:\sigma ) = R(P_1:P_1:\sigma ).$$  
But it follows from the integral formula that $J(P_1:P_1:\sigma :\nu )$ is the identity operator
for $\nu $ in the region of convergence, and hence for all $\nu $.  Thus 
$R(P_1:P_1:\sigma )$ is a non-zero constant times the identity operator and so
$R(P_2:P_1:\sigma )$ is invertible.  \qed

Let $x \in N_G(M)$.  Then if $P_1 = MN_1$ is a parabolic subgroup of $G^0$ with Levi component $M$,
so is $P_2 = xPx^{-1} = xMx^{-1}xN_1x^{-1} = MN_2$.  Let $dn_1$ and $dn_2$ denote normalized Haar
measure on $N_1$ and $N_2$ respectively.  
(Thus $dn_i$ assigns measure one to $K \cap N_i, i = 1,2$.)
  Define $\alpha _{P_1} (x) \in \bfR ^+$ by
$$\int _{N_2} \phi (n_2) dn_2 = \alpha _{P_1}(x)\int _{N_1} \phi (xn_1x^{-1})dn_1 , \phi \in
C_c^{\infty }(N_2).$$ For all $m \in M$, 
$$\alpha _{P_1}(m) = \delta _{P_1}(m).$$  Further, if $x \in N_G(K) \cap N_G(M)$, then
$\alpha _{P_1}(x) = 1$.  The following lemma is an easy consequence of the definition.

\proclaim Lemma 4.10.  Let $x,y \in N_G(M)$.  Then $$\alpha _P(yx) = \alpha
_{xPx^{-1}}(y)\alpha _P(x).$$  Moreover if $m \in M, y \in N_G(M)$, then
$$\alpha _P(ym) = \alpha _P(y)\delta
_P(m).$$ \endproclaim

Now let $N_G(\sigma ) = \{ g \in N_G(M): \sigma ^g \simeq \sigma \}$ and let $W_G(\sigma ) =
N_G(\sigma )/M$.   Let $W_{G^0}(\sigma ) = (N_G(\sigma ) \cap G^0)/M$.
If
$w \in W_G(\sigma )$, $\sigma $ can be extended to a representation of the group $M_w$ generated
by $M$ and any representative $n_w$ for $w$ in $N_G(\sigma )$.  Fix such an extension and denote
it by $\sigma _w$. 
Now we can define an intertwining operator
$$A_P(w): \calH _{w^{-1}Pw}(\sigma ) \rightarrow \calH _P(\sigma )$$ by 
$$(A_P(w)f)(x) =
\sigma _w(n_w)\alpha _{w^{-1}Pw}(n_w)^{1\over 2} f(xn_w).$$  
The $\alpha _{w^{-1}Pw}$ term is not used in the connected case because coset representatives $n_w$ can be
chosen in $K$ where $\alpha _{w^{-1}Pw} = 1$.  In the general case we don't know if there is a natural
choice of coset representatives with $\alpha _{w^{-1}Pw} = 1$.  Thus we add the $\alpha _{w^{-1}Pw}$
term so that $A_P(w)$ is independent of the coset representative $n_w$ for $w$.  $A_P(w)$ does
however depend on the choice of the extension $\sigma _w$. 

\proclaim Lemma 4.11.  The intertwining operator $A_P(w)$ is independent of the choice of coset
representative $n_w$ for $w \in W_G(\sigma )$.  For $w_1, w_2 \in W_G(\sigma )$ there is a
non-zero complex constant $c_P(w_1,w_2)$ so that
$$A_P(w_1w_2) = c_P(w_1,w_2)A_P(w_1)A_{w_1^{-1}Pw_1}(w_2).$$  \endproclaim

\proof Using Lemma 4.10 we have for any $m
\in M, f \in \calH _{w^{-1}Pw}(\sigma )$, 
$$\sigma _w(n_wm)\alpha _{w^{-1}Pw}(n_wm)^{1\over 2} f(xn_wm) $$ $$=
\sigma _w(n_w)\sigma (m) \alpha _{w^{-1}Pw}(n_w)^{1\over 2}\delta
_{w^{-1}Pw}(m)^{1\over 2}\sigma (m)^{-1}\delta _{w^{-1}Pw}(m)^{-{1\over 2}} f(xn_w) $$
$$= \sigma _w(n_w)\alpha _{w^{-1}Pw}(n_w)^{1\over 2} f(xn_w).$$  Thus 
the intertwining operator is independent of the choice of coset
representative. 
  
 Let $w_1, w_2 \in W_G(\sigma )$ and fix coset representatives $n_1, n_2$ for $w_1, w_2$
respectively.  Then we can use $n_1n_2$ as a coset representative for $w_1w_2$.   We will also write
$P_1 = w_1^{-1}Pw_1$ and $P_{12} =  w_2^{-1} w_1^{-1}Pw_1w_2$.
For any $w \in W_G(\sigma )$ and representative $n_w$ for $w$ we have
$$\sigma _w(n_w) \sigma (m) \sigma _w(n_w)^{-1} = \sigma (n_wmn_w^{-1})$$ for all $m \in M$.  Thus
$$ \sigma _{w_1w_2}( n_1n_2) \sigma _{w_2}(n_2)^{-1}\sigma _{w_1}(n_1)^{-1} $$
is a non-zero self-intertwining operator for $\sigma $ and hence there is a non-zero constant $c $
so that  $$ \sigma _{w_1w_2}( n_1n_2) = c\sigma _{w_1}(n_1) \sigma _{w_2}(n_2) .$$
Further, from Lemma 4.10 we have
$$ \alpha _{P_{12}}(n_1n_2) = \alpha _{P_1}(n_1)\alpha _{P_{12}}(n_2).$$ 
 
By definition, for any $x
\in G, f \in \calH _{P_{12}}(\sigma )$,
$$A_P(w_1w_2)f(x) = \sigma _{w_1w_2} (n_1n_2) \alpha _{P_{12}} (n_1n_2)^{1\over 2} f(xn_1n_2) $$
$$=c\sigma _{w_1}(n_1) \alpha _{P_1}(n_1)^{1\over 2}\sigma _{w_2}(n_2)
\alpha _{P_{12}}(n_2)^{1\over 2}
f(xn_1n_2)$$
$$ =  c A_P(w_1)A_{w_1^{-1}Pw_1}(w_2) f(x).$$  \qed

If $w \in W_{G^0}(\sigma )$ and $n_w$ is chosen to be in $K$, then we also have
$$A^0(w): \calH _{w^{-1}Pw}^0(\sigma ) \rightarrow \calH _P^0(\sigma )$$ given by $$(A^0(w)\phi
)(x) = \sigma _w(n_w) \phi (xn_w).$$ 
The compositions 
 $$R^0(w_0,\sigma ) = A^0(w_0)R^0(w_0^{-1}Pw_0:P:\sigma ), w_0 \in W_{G^0}(\sigma ), $$ 
$$R(w,\sigma ) = A_P(w)R(w^{-1}Pw:P:\sigma ), w \in W_G(\sigma ), $$ give
self-intertwining operators for  $I_P^0(\sigma )$ and $I_P(\sigma )$ respectively. 
 
\proclaim Lemma 4.12.  There is a cocycle $\eta $ so that 
 $$R(w_1w_2, \sigma ) = \eta (w_1,w_2) R(w_1,\sigma )R(w_2,\sigma )$$
for all $w_1,w_2 \in W_G(\sigma )$. 
 \endproclaim

\proof  Using Lemmas 4.8 and 4.11, we have $$R(w_1w_2, \sigma ) =
A_P(w_1w_2)R(w_2^{-1}w_1^{-1}Pw_1w_2:P:\sigma )$$ $$ = 
c_P(w_1,w_2)A_P(w_1)A_{w_1^{-1}Pw_1}(w_2)
R(w_2^{-1}w_1^{-1}Pw_1w_2:w_2^{-1}Pw_2:\sigma )R(w_2^{-1}Pw_2:P:\sigma ).$$
We will show that there is a non-zero constant $c_P'(w_1,w_2)$ so that
 $$A_{w_1^{-1}Pw_1}(w_2)R(w_2^{-1}w_1^{-1}Pw_1w_2:w_2^{-1}Pw_2: \sigma )A_P(w_2)^{-1} $$ 
$$= 
c_P'(w_1,w_2)R(w_1^{-1}Pw_1:P:\sigma ).$$   When this is established we will have
$R(w_1w_2, \sigma )= $ 
$$c_P(w_1,w_2)c_P'(w_1,w_2)A_P(w_1)R(w_1^{-1}Pw_1:P:\sigma )
A_P(w_2)R(w_2^{-1}Pw_2:P:\sigma ) $$
$$=c_P(w_1,w_2)c_P'(w_1,w_2) R(w_1,\sigma )R(w_2,\sigma ).$$  Thus $\eta $ can be
defined by $\eta (w_1,w_2) = c_P(w_1,w_2)c_P'(w_1,w_2)$. 
It is immediate from the formulas for the composition of
the operators $R(w_i,\sigma)$ that $\eta$ is a $2$--cocycle.

In order to prove the above identity we need to go back to the original
definition of  the standard intertwining operators.
$W_G(\sigma )$ acts on $\fraka _{\bfC }^*$ and $\chi _{\nu }(n_w^{-1}mn_w) = \chi _{w\nu }(m)$
if $n_w$ is any representative of $w \in W_G(\sigma )$.  
Now for any $\nu \in \fraka _{\bfC }^*
, f \in \calH _{w^{-1}Pw}(\sigma , \nu
), x \in G, m \in M, n \in N$, if $A_P(w)$ is defined exactly as above, we have
$$A_P(w)f(xmn) = \sigma _w(n_w)\alpha _{w^{-1}Pw}(n_w)^{1\over 2} f(xmnn_w)$$
$$=  \sigma _w(n_w)\alpha _{w^{-1}Pw}(n_w)^{1\over 2}\delta _{w^{-1}Pw}(n_w^{-1}m^{-1}n_w)^{-{1\over 2}} 
\sigma (n_w^{-1}m^{-1}n_w) \chi _{\nu }(n_w^{-1}m^{-1}n_w) f(xn_w)$$
$$=  \delta _P(m)^{-{1\over 2}}\sigma (m)^{-1} \chi _{w\nu }(m) ^{-1}
\sigma _w(n_w)\alpha _{w^{-1}Pw}(n_w)^{1\over 2} f(xn_w)$$
$$=  \delta _P(m)^{-{1\over 2}}\sigma (m)^{-1} \chi _{w\nu }(m) ^{-1}A_P(w)f(x).$$
Thus $A_P(w)$ maps $\calH _{w^{-1}Pw}(\sigma , \nu )$ to $\calH _P (\sigma , w\nu )$.  Now
write $P_1 = w_1^{-1}Pw_1, P_2 = w_2^{-1}Pw_2,\hfil \break P_{12} = w_2^{-1}w_1^{-1}Pw_1w_2$ and
 consider the composition
$$A_{P_1}(w_2)J(P_{12}:P_2: \sigma :w_2^{-1}\nu )A_P(w_2)^{-1}$$
It maps
$$\calH _P(\sigma , \nu ) \rightarrow \calH _{P_2}(\sigma , w_2^{-1}\nu ) \rightarrow
\calH _{P_{12}}(\sigma , w_2^{-1}\nu ) \rightarrow \calH _{P_1}(\sigma ,
\nu ).$$

If $\nu $ is in the region of convergence for the integral formula for 
$J(P_{12}:P_2: \sigma :w_2^{-1}\nu )$ and $d\barn _2$ denotes normalized Haar measure on 
$\barN _2 \cap N_{12}$, then we
have $$A_{P_1}(w_2)J(P_{12}:P_2: \sigma :w_2^{-1}\nu )A_P(w_2)^{-1}f(x) $$
$$=
\sigma _{w_2}(n_{w_2})\alpha _{P_{12}} (n_{w_2})^{1\over 2}J(P_{12}:P_2: \sigma :w_2^{-1}\nu
)A_P(w_2)^{-1}f(xn_{w_2}) $$
$$ = \sigma _{w_2}(n_{w_2})\alpha _{P_{12}} (n_{w_2})^{1\over 2}
\int _{ \barN _2 \cap N_{12}} A_P(w_2)^{-1}f(xn_{w_2}\barn _2)d \barn _2  $$
$$=  \sigma _{w_2}(n_{w_2})\alpha _{P_{12}} (n_{w_2})^{1\over 2}\int _{ \barN _2 \cap N_{12}}
\alpha _{P_2} (n_{w_2})^{-{1\over 2}}  \sigma _{w_2}(n_{w_2})^{-1}f(xn_{w_2}\barn _2 n_{w_2}^{-1})d
\barn _2 $$ $$=\alpha _{P_{12}} (n_{w_2})^{1\over 2}\alpha _{P_2} (n_{w_2})^{-{1\over 2}} \int
_{\barN _2 \cap N_{12}} f(xn_{w_2}\barn _2 n_{w_2}^{-1}) d\barn _2.$$
But $\barN _2 \cap N_{12} =w_2^{-1}(\barN \cap N_1)w_2$ so if $d\barn $ denotes normalized Haar
measure on $\barN \cap N_1$, there is a positive real number $r$ so that 
$$ \int _{\barN _2 \cap
N_{12}} \phi (n_{w_2}\barn _2 n_{w_2}^{-1}) d\barn _2= r  \int _{\barN \cap N_1}
\phi (\barn ) d\barn $$ for all $\phi \in C_c^{\infty }(\barN \cap N_1)$.  Thus  
 $$A_{P_1}(w_2)J(P_{12}:P_2: \sigma :w_2^{-1}\nu )A_P(w_2)^{-1}f(x) $$
$$= \alpha _{P_{12}} (n_{w_2})^{1\over 2}\alpha _{P_2} (n_{w_2})^{-{1\over 2}}
r  \int _{\barN \cap N_1} f(x\barn )
d\barn $$
$$ =  \alpha _{P_{12}} (n_{w_2})^{1\over 2}\alpha _{P_2} (n_{w_2})^{-{1\over 2}}
r  J(P_1:P:\sigma :\nu )f(x).$$
Setting $$c_P''(w_1,w_2) =   \alpha _{P_{12}} (n_{w_2})^{1\over 2}\alpha _{P_2} (n_{w_2})^{-{1\over
2}} r $$  we have
 $$A_{P_1}(w_2)J(P_{12}:P_2: \sigma :w_2^{-1}\nu )A_P(w_2)^{-1}
= c_P''(w_1,w_2) J(P_1:P:\sigma :\nu )$$ for all $\nu $ in the region of convergence for 
the integral formula for 
$J(P_{12}:P_2: \sigma :w_2^{-1}\nu )$.  
By analytic continuation, the identity is valid for all $ \nu  $.
Now divide both sides of the equation by $r(P_{12}:P_2:\sigma :w_2^{-1}\nu )$ and evaluate at $\nu =
0$. We obtain
$$ A_{P_1}(w_2)R(P_{12}:P_2: \sigma )A_P(w_2)^{-1} = c_P'(w_1,w_2) R(P_1:P:\sigma ) $$ where
$$c_P'(w_1,w_2) =  c_P''(w_1,w_2) r(P_1:P:\sigma :0)r(P_{12}:P_2:\sigma :0)^{-1}.$$
\qed

For $\phi \in \calH _P^0(\sigma )$ define $f = \Phi (\phi) \in \calH
_P(\sigma )$ such that $f(x) = 0$ if $x \not \in G^0$ and $f(x_0) = \phi (x_0)$ if $x_0 \in G^0$.

\proclaim Lemma 4.13.  Let $f = \Phi (\phi)$ as above.  Then for all $w, w_1 \in W_G(\sigma ),
x_0 \in G^0$, $$R(w,\sigma )f(x_0n_{w_1}^{-1}) =0$$ unless $w = w_1 w_0, w_0 \in W_{G^0}(\sigma )$.  
If $w = w_1 w_0, w_0 \in W_{G^0}(\sigma )$, then $R(w,\sigma )f(x_0n_{w_1}^{-1}) =$
$$\eta (w_1, w_0)\sigma _{w_1}(n_{w_1})\alpha _{w_1^{-1}Pw_1}(n_{w_1})^{1\over 2}
R^0(w_1^{-1}Pw_1:P:\sigma ) R^0(w_0:\sigma )\phi (x_0).$$  \endproclaim

\proof  For any $x \in G, w \in W_G(\sigma )$, using Corollary 4.7,
$R(w^{-1}Pw:P:\sigma )f(x) = R^0(w^{-1}Pw:P:\sigma )\phi ' (1)$ where $\phi ' $ is the restriction
to $G^0$ of $l(x^{-1})f$. Now since $f = \Phi (\phi )$ is supported on $G^0$, $\phi ' =0$ unless $x
\in G^0$.  If $x_0 \in G^0$, then 
$R(w^{-1}Pw:P:\sigma )f(x_0) = R^0(w^{-1}Pw:P:\sigma )\phi (x_0)$.  

Now by definition, $R(w,\sigma )f(x_0n_{w_1}^{-1}) =$ $$ \sigma
_{w}(n_{w})\alpha _{w^{-1}Pw}(n_w)^{1\over 2}R(w^{-1}Pw:P:\sigma ) f(x_0n_{w_1}^{-1}n_w).$$  By the
above this is zero unless $n_{w_1}^{-1}n_w \in G^0$, that is unless $w = w_1 w_0, w_0 \in
W_{G^0}(\sigma )$.  In this case, using Lemma 4.12, 
 $$R(w_1 w_0,\sigma )f(x_0n_{w_1}^{-1}) = \eta (w_1, w_0) R(w_1:\sigma )R(w_0:\sigma )f(x_0n_{w_1}^{-1})
$$ $$ = \eta (w_1, w_0)\sigma _{w_1}(n_{w_1})\alpha _{w_1^{-1}Pw_1}(n_{w_1})^{1\over 2}
 R(w_1Pw_1^{-1}:P:\sigma )R(w_0:\sigma )f(x_0) $$ $$ =
\eta (w_1, w_0)\sigma _{w_1}(n_{w_1})\alpha _{w_1^{-1}Pw_1}(n_{w_1})^{1\over 2}
 R^0(w_1Pw_1^{-1}:P:\sigma )R^0(w_0:\sigma )\phi (x_0).$$ \qed

 Recall that if we write $W_{G^0}^0(\sigma )$ for the subgroup of elements 
$w \in W_{G^0}(\sigma )$ such that $R^0(w,\sigma ) $ is scalar, then $W_{G^0}^0(\sigma )= W(\Phi
_1)$ is generated by reflections in a set $\Phi _1$ of reduced roots of $(G,A)$.   Let $\Phi ^+$ be
the positive system of reduced roots of $(G,A)$ determined by $P$ and let $\Phi _1^+ = \Phi _1 \cap
\Phi ^+$.  If we define
 $$R^0_{\sigma } = \{ w \in W_{G^0}(\sigma ): w\beta \in \Phi ^+ 
{\rm \ for \ all \ } \beta \in \Phi _1^+ \},$$  then $W_{G^0}(\sigma )$ is the semidirect product of
$R^0_{\sigma }$ and $W(\Phi _1)$.   $R^0_{\sigma }$ is called the $R$-group for $I_P^0(\sigma
)$ and the operators $$\{R^0(r,\sigma ), r \in R^0_{\sigma } \}$$ form a basis for the algebra of
intertwining operators of $I_P^0(\sigma )$.  We will define $$R_{\sigma } =  \{ w \in W_G(\sigma ):
w\beta \in \Phi ^+ {\rm \ for \ all \ } \beta \in \Phi _1^+ \}.$$  Clearly 
$R_{\sigma } \cap G^0 = R^0_{\sigma }$.

\proclaim Lemma 4.14.  $R(w, \sigma )$ is scalar if $w \in W(\Phi _1)$, and
$W_G(\sigma )$ is the semidirect product of $W(\Phi _1)$ and $R_{\sigma
}$.  \endproclaim

\proof   Fix $w \in W(\Phi _1), n_w\in K $ a representative for $w$, and $f \in \calH
_P(\sigma )$.   Then for all $x \in G$, $$R(w, \sigma
)f(x) = \sigma _w (n_w) R(w^{-1}Pw:P:\sigma )f(xn_w) = \sigma _w (n_w) R^0(w^{-1}Pw:P:\sigma
)\phi (n_w) = $$ where $\phi (x_0) = f(xx_0)$ for all $x_0 \in G^0$.  Thus $$R(w, \sigma
)f(x) = R^0(w,
\sigma )\phi (1).$$
 But since $w \in W(\Phi _1)$ there is a constant $c_w$ such that
$R^0(w, \sigma )\phi = c_w \phi $ for all $\phi \in \calH ^0_P(\sigma )$.
Thus
 $$R(w, \sigma )f(x) =  c_w \phi (1) = c_w f(x)$$ so $R(w, \sigma )$ is scalar.  

We must show that for any $w \in W_G(\sigma ), w\Phi _1 = \Phi _1$ so that 
 $$R_{\sigma } =  \{ w \in W_G(\sigma ):
w\Phi _1^+ = \Phi _1^+ \}.$$   Then as in the connected case it will be clear that $W_G(\sigma )$ is the
semidirect product of $R_{\sigma }$ and $W(\Phi _1)$.  
Let $w \in W_G(\sigma ), \alpha \in \Phi _1$, and let $s_{\alpha }\in W(\Phi _1)$
denote the reflection corresponding to $\alpha $.  Then $ws_{\alpha }w^{-1} = s_{w\alpha }\in
W_{G^0}(\sigma )$.  But $R(s_{\alpha }:\sigma )$ is scalar so that using Lemma 4.12, so is 
$R(s_{w\alpha }, \sigma ) = R(ws_{\alpha }w^{-1}:\sigma )$.  This implies as above that 
$R^0(s_{w\alpha }, \sigma )$ is scalar.  Hence $s_{w\alpha } \in W(\Phi _1 )$ and $w\alpha
\in \Phi _1$. \qed

\proclaim Lemma 4.15.  The dimension of $C(\sigma )$ is equal to $[R_ {\sigma }]$. \endproclaim

\proof  By Lemma 4.4, $$\dim C(\sigma ) = [N_G(\sigma )/N_{G^0}(\sigma )] \dim C^0(\sigma ) =
 [W_G(\sigma )/W_{G^0}(\sigma )] [R^0_{\sigma }] $$ $$= [W_G(\sigma )/W_{G^0}(\sigma )]
[W_{G^0}(\sigma ) /W(\Phi _1)] = [W_G(\sigma )/W(\Phi _1 )] = [R_{\sigma }].$$ 
\qed 

\proclaim Theorem 4.16.  The operators $\{R(r,\sigma ), r \in R_{\sigma } \}$ form a basis for the
algebra of intertwining operators of $I_P(\sigma )$. \endproclaim

\proof   By Lemma 4.15 it suffices to show that the operators are linearly independent.  Suppose
that $c_w, w \in R_{\sigma }$, are constants so that $$\sum _{w \in R_{\sigma }}c_w R(w, \sigma
)f (x) = 0$$ for all $f \in \calH _P(\sigma ), x \in G$.  Fix $w_1 \in R_{\sigma }$.  
Then for all $f = \Phi (\phi)$ with
$\phi \in \calH _P^0(\sigma )$ and all $ x_0 \in G^0$, we have
 $$\sum _{w \in R_{\sigma }}c_w R(w, \sigma
)f(x_0n_{w_1}^{-1}) = 0.$$
Now by Lemma 4.13, $ R(w, \sigma
)f(x_0n_{w_1}^{-1}) = 0$ unless $w = w_1w_0$ where $w_0 \in W_{G^0}(\sigma ) \cap R_{\sigma } =
R^0_{\sigma }$.  Now again using Lemma 4.13,
 $$0 = \sum _{w _0\in R^0_{\sigma }}c_{w_1w_0} R(w_1w_0, \sigma
)f(x_0n_{w_1}^{-1}) $$ $$= \sum _{w _0\in R^0_{\sigma }}c_{w_1w_0} \eta (w_1,w_0)
\sigma _{w_1}(n_{w_1})\alpha _{w_1^{-1}Pw_1}(n_{w_1})^{1\over 2} 
R^0(w_1^{-1}Pw_1:P:\sigma ) R^0(w_0:\sigma )\phi (x_0)$$
$$=\sigma _{w_1}(n_{w_1})\alpha _{w_1^{-1}Pw_1}(n_{w_1})^{1\over 2}
 R^0(w_1^{-1}Pw_1:P:\sigma )\sum _{w _0\in R^0_{\sigma }}c_{w_1w_0} \eta (w_1,w_0)
  R^0(w_0:\sigma )\phi (x_0).$$   Thus
$$\sum _{w _0\in R^0_{\sigma }}c_{w_1w_0} \eta (w_1,w_0)
  R^0(w_0:\sigma )\phi (x_0) =0$$ for all $\phi \in \calH _P^0(\sigma )$ and all $ x_0 \in G^0$.
 Now since we know that the operators $R^0(w_0:\sigma )$ are
linearly independent on $\calH _P^0(\sigma )$, we can conclude that the $c_{w_1w_0} \eta (w_1,w_0)$
and hence the $ c_{w_1w_0} $ are all zero.  \qed

As in Arthur\Lspace \Lcitemark 2\Rcitemark \Rspace{} we now have to deal with the cocycle $\eta $ of Lemma 4.12.  Fix a finite central
extension $$1 \rightarrow Z_{\sigma } \rightarrow \tilde R _{\sigma } \rightarrow R_{\sigma }
\rightarrow 1$$ over which $\eta $ splits.  Also define the functions $\xi _{\sigma }:\tilde R
_{\sigma }
\rightarrow \bfC ^*$ and the character $\chi _{\sigma }$ of $Z_{\sigma }$ as in Arthur\Lspace \Lcitemark 2\Rcitemark \Rspace{}.  Then we
obtain a homomorphism $$\tilde R (r,\sigma ) = \xi _{\sigma } ^{-1}(r)R(r, \sigma ), r \in \tilde R
_{\sigma },$$ of $\tilde R _{\sigma }$ into the group of unitary intertwining operators for
$I_P(\sigma )$ which transforms by
$$\tilde R (zr, \sigma ) = \chi _{\sigma }(z)^{-1} \tilde R (r,\sigma ), z \in Z_{\sigma }, r \in
\tilde R _{\sigma }.$$

Now we can define a representation $\calR $ of $\tilde R _{\sigma } \times G$
on $\calH _P(\sigma )$ given by $$\calR (r,x) = \tilde R (r,\sigma ) I_P(\sigma ,x), 
r \in \tilde R _{\sigma }, x \in G.$$  Let  $\Pi (\tilde R _{\sigma
}, \chi _{\sigma })$ denote the set of irreducible representations of $\tilde R _{\sigma }$ with
$Z_{\sigma }$ central character $\chi _{\sigma }$ and let $\Pi _{\sigma }(G)$ denote the set of
irreducible constituents of $I_P(\sigma )$.  

\proclaim Theorem 4.17.  There is a bijection $\rho \mapsto \pi _{\rho }$ of $\Pi (\tilde R _{\sigma
}, \chi _{\sigma })$ onto $\Pi _{\sigma }(G)$ such that 
$$\calR = \oplus _{\rho \in \Pi (\tilde R _{\sigma }, \chi _{\sigma })}\ \ \  (\rho \chk \otimes
\pi _{\rho }).$$  \endproclaim

\proof  Write the decomposition of $\calR $ into irreducibles as
$$\calR = \sum _{\rho , \pi } m_{\rho , \pi } \ (\rho \chk \otimes \pi ) $$ where
$\rho $ runs over $\Pi (\tilde R _{\sigma }, \chi _{\sigma }), \pi $ runs over
$\Pi _{\sigma }(G)$, and each $ m_{\rho , \pi } \geq 0$.  This corresponds to a decomposition
$$ \calH _P (\sigma ) \simeq \sum _{\rho , \pi } (V_{\rho \chk} \otimes W _{\pi })^
{m_{\rho , \pi }}$$ where $V_{\rho \chk}$ and $W _{\pi }$ denote the representation spaces
for the irreducible representations $\rho \chk $ and $\pi $.  For each $\rho $, we also write
$$\calH _{\rho } \simeq  V_{\rho \chk} \otimes \sum _{ \pi } W _{\pi }^
{m_{\rho , \pi }}$$ for the $\rho $-isotypic component of $\calH _P(\sigma )$.  Each subspace
$\calH _{\rho }$ is invariant under the action of $\calR $, in particular by all of the
intertwining operators $\tilde R (r,\sigma ), r \in \tilde R _{\sigma }$.  Since these
intertwining operators span the space $C(\sigma )$ of $G$-intertwining operators for $\calH
_P(\sigma )$, each $T \in C(\sigma )$ must satisfy $T (\calH _{\rho }) \subset \calH _{\rho }$ for
every $\rho $. 

We will first show that given $\pi $, there is at most one $\rho $ such that $m _{\rho , \pi } >
0$.  So fix $\pi $ and suppose that there are $\rho _1$ and $\rho _2 $ such that 
$m _{\rho _i, \pi } > 0, i = 1,2$.  Then $W_{\pi }$ occurs as a $G$-summand of  
both $\calH {\rho _1}$ and $\calH {\rho _2}$.  Thus there is $T_{12} \not = 0$ in $Hom _G(\calH
{\rho _1},\calH {\rho _2})$.  We can extend $T_{12}$ to an element $T$ of $C(\sigma )$ by
setting $T = T_{12}$ on $\calH _{\rho _1}$ and $T = 0$ on $\calH _{\rho }, \rho \not = \rho _1$.
Thus there is $T $ in $C(\sigma )$ such that 
$T (\calH _{\rho _1}) \subset \calH _{\rho _2}$.  But by the remark in the previous paragraph, 
$T (\calH _{\rho _1}) \subset \calH _{\rho _1}$.  Now since $T (\calH _{\rho _1}) \not = 0$, we
must have $\rho _1 \simeq \rho _2$.  

Now fix $\rho $ and look at $\calH _{\rho } \simeq  V_{\rho \chk} \otimes W$ where $W = \sum _{
\pi } W _{\pi }^ {m_{\rho , \pi }}$.  We will show that $W$ is irreducible as a
$G$-module. Thus suppose that $W = W_1 \oplus W_2$ where $W_1, W_2$ are $G$-submodules of $W$. 
Then we can define $T \in C(\sigma )$ by $T (v \otimes (w_1 + w_2)) = v \otimes w_1$ if $v \in 
V_{\rho \chk}, w_1 \in W_1, w_2 \in W_2$, and $T = 0$ on $\calH _{\rho '}$ if $\rho ' \not \simeq
\rho $.  Now since $T \in C(\sigma )$, we can write $T = \sum _r c_r \tilde {\calR }(r, \sigma )$
where
$r$ runs over $\tilde R _{\sigma }$. Thus for all $v \in  V_{\rho \chk}, w_1 \in W_1, w_2 \in W_2$,
we have
$$ v \otimes w_1 = T (v \otimes (w_1 + w_2)) = (\sum _r c_r \rho \chk (r)v ) \otimes w_1
+ (\sum _r c_r \rho \chk (r)v )\otimes w_2.$$
Suppose $W_2 \not = 0$. This implies that $\sum _r c_r \rho \chk (r)v =0$ for all $v$ so that
$v \otimes w_1 = 0$ for all $w_1$.  Thus $W_1 = \{ 0 \}$ and hence $W$ is irreducible.  But this
implies that $m_{\rho , \pi } \leq 1$ for all $\pi $ and that
there is at most one $\pi $ such that $m_{\rho , \pi } =1$.

Define $$\Pi _1 = \{ \rho \in \Pi (\tilde R _{\sigma }, \chi _{\sigma }) : m_{\rho , \pi } =1
{\rm \ for \ some \ } \pi \}.$$
For each $\rho \in \Pi _1$ we have shown that the representation $\pi $ such that
 $ m_{\rho , \pi
} =1$ is unique.  Thus we will call it $\pi _{\rho }$.  Further, we have shown that $\pi _{\rho 
_1} \simeq \pi _{\rho _2}$ just in case $\rho _1 \simeq \rho _2$.  Further, by definition of
$\Pi _{\sigma }(G)$, each $\pi \in \Pi _{\sigma }(G)$ occurs in $\calH_P(\sigma )$ and so must
be of the form $\pi _{\rho }$ for some $\rho \in \Pi _1$.  Thus to complete the proof of the
theorem we need only show that $\Pi _1 = \Pi (\tilde R _{\sigma
}, \chi _{\sigma })$.
 
Since $$\calH _P(\sigma ) \simeq \sum  _{\rho \in \Pi
_1}  (V_{\rho \chk} \otimes W _{\pi _{\rho }}),$$  we must have $\dim C(\sigma ) = \sum
_{\rho \in \Pi _1} (deg \ \rho )^2 $.  But since the $R(r, \sigma ), r \in R_{\sigma }$, form a
basis for $C(\sigma )$, we know that $$\dim C(\sigma ) = [R_{\sigma }] = \sum _{\rho \in \Pi
(\tilde R _{\sigma }, \chi _{\sigma })} (deg \ \rho )^2.$$  Thus $\Pi _1 = \Pi (\tilde R _{\sigma
}, \chi _{\sigma })$.  \qed

\medbreak
\noindent
{\bf Remark 4.18} Suppose now that $G/G^0$ is cyclic.  Then in\Lspace \Lcitemark 1\Rcitemark \Rspace{},
Arthur predicts a dual group construction of a group $R_{\psi,\sigma},$
in terms of the conjectural parameter $\psi$ for the $L$--packet of 
$\sigma,$ which should also describe the components of $\Ind_P^G(\sigma).$
In particular, $R_{\psi,\sigma}=W_{\psi,\sigma}/W_{\psi,\sigma}^0,$
where these groups are defined in terms of centralizers of the
image of $\psi.$
Furthermore, Arthur has conjecturally identified $W_{\psi,\sigma}^0$
with $W(\Phi_1)$ and $W_{\psi,\sigma}$ with $W_G(\s).$
Thus, if the conjectural parameterization exists in the connected case,
and Shelstad's Theorem\Lspace \Lcitemark 12\Rcitemark \Rspace{} extends to the $p$-adic  case,
then it must be the case that $R_{\sigma}\simeq R_{\psi,\sigma}.$
That is, we have shown that there is a group side construction
of Arthur's $R$--group, if such an object exists.
(For more details see\Lspace \Lcitemark 1\Rcitemark \Rspace{}
and\Lspace \Lcitemark 6\Rcitemark \Rspace{}, particularly Sections 1 and 4.)

\section {$R$-groups for $\protect\Ind_P^G(\sigma )$}

In this section we will study representations of $G$ which are induced from discrete series
representations of a parabolic subgroup $P$ of $G$.  Thus we revert to the notation that
parabolic subgroups of $G^0$ are denoted by $P^0$.

Let $P=MN$ be a cuspidal parabolic subgroup of $G$.  Let $\sigma $ be an irreducible discrete
series representation of $M$ and let $\sigma _0$ be an irreducible constituent of the
restriction of $\sigma $ to $M^0$.  We want to find a basis for the intertwining
algebra $C(\sigma )$ of the induced representation
$\Ind _P^G(\sigma )$.
Since $\sigma $ is
contained in $\Ind_{M^0}^M (\sigma _0)$ we know $\Ind _P^G(\sigma )$ is contained
in
$\Ind _P^G(\Ind _{M^0}^M(\sigma _0)) \simeq \Ind _{P^0}^G (\sigma _0)$.   In \S 4 we found a
basis for the intertwining algebra $C(\sigma _0)$ of $\Ind _{P^0}^G (\sigma _0)$.  We will see
how to obtain  a
basis for $C(\sigma )$
 by restricting the intertwining operators defined in \S 4. 

We first need
to embed $\sigma $ in a family
$\sigma _{\nu }, \nu \in \fraka _{\bfC }^*$, where $\fraka $ is the real Lie algebra of the
split component $A$ of $M$.
Write $X(M), X(A)$ for the groups of rational characters of $M,A$ respectively. 
 Let $$r: X(M) \otimes _{\bfZ } \bfR \rightarrow X(A) \otimes
_{\bfZ } 
\bfR$$ be the map given by restriction.   That is $r( \chi \otimes t) = \chi |_A \otimes t$ for
$\chi \in X(M), t \in \bfR$.  

\proclaim Lemma 5.1.  The homomorphism $r: X(M) \otimes _{\bfZ } \bfR \rightarrow X(A) 
\otimes _{\bfZ } \bfR$ is
surjective.  \endproclaim

\proof   Since $G$ is a linear group we have an embedding of $G$ in $L = GL(V),$  where $V$ is a finite dimensional $F$--vector space. Since $A$ is a
split torus, the action of $A$ on $V$ can be diagonalized.  For any $\chi \in X(A)$ let
$V(\chi ) = \{ v \in V : av = \chi  (a) v$ for all $ a \in A \} $.  Let $\chi _i , 1 \leq i
\leq k$, denote the distinct elements of $X(A)$ such that $V_i = V(\chi _i) \not = \{ 0 \}$.  
We can identify $a \in A$ with the block diagonal matrix with diagonal entries
$\chi _i(a) I_{d_i}$ where $d_i = \dim V_i $ and $I_{d_i}$ denotes the identity matrix
of size $d_i, 1 \leq i \leq k$.  

Since $M = C_G(A)$, we have $M \subset C_L(A) \simeq GL(V_1) \times GL(V_2) \times ... \times
GL(V_k)$.  For each $1 \leq i \leq k$ we can define $\det _i \in X(M)$ by $\det  _i
(m_1,...,m_k) = \det m_i$.  Now $\det _i \otimes d_i^{-1} \in X(M) \otimes _{\bfZ }\bfR $
 and for $a
\in A$,  $\det _i \otimes d_i^{-1} (a) = \chi _i(a)^{d_i} \otimes d_i ^{-1} = \chi _i
(a) \otimes 1$.  Thus $r(\det _i \otimes d_i^{-1}) = \chi _i \otimes 1$.  
The $\chi _i, 1 \leq i \leq k$, are generators of
$X(A)$, although they need not be independent.  
Thus $r$ is surjective.  \qed
 
Let $X_0(M) = \{ \chi \in X(M) : \chi |_{M^0} = 1 \}$.

\proclaim Lemma 5.2.  The kernel of $r$ is $X_0(M) \otimes _{\bfZ } \bfR$. \endproclaim

\proof  Suppose $\chi \otimes t $ is in the kernel of $r$ where $\chi \in X(M)$ and $t \in
\bfR$.  If $t = 0$ then $\chi \otimes t$ is the identity element.  Assume $t \not = 0$. Then
$|\chi (a)|_F^t = 1$ for all $a \in A$ implies that $|\chi (a)|_F =
1$ for all $a \in A$.  Since $\chi |_A$ is a rational character of a split torus this implies
that $\chi |_A = 1$.  But restriction from $X(M^0)$ to $X(A)$ is injective 
\Lcitemark 13\LIcitemark{}, Lemma 0.4.1\RIcitemark \Rcitemark \Rspace{}, so that
$\chi |_{M^0} = 1$.  Thus $\chi \in X_0(M)$.  \qed

Recall the homomorphism $H_{M^0}:M^0 \rightarrow Hom (X(M^0), \bfZ)$ defined by
$$<H_{M^0}(m), \chi > = log _q |\chi (m)|_F, m \in M^0, \chi \in X(M^0).$$  
Define an analogous homomorphism $H_M:M \rightarrow Hom (X(M), \bfZ )$ by
 $$<H_M(m), \chi > = log _q |\chi (m)|_F, m \in M, \chi \in X(M).$$  

\proclaim Lemma 5.3.  Suppose that $\chi \in X_0(M)$.  Then $<H_M(m), \chi > = 0$ for all $m \in
M$.  \endproclaim.

\proof  Let $\chi \in X_0(M)$.  Thus $\chi (m_0) = 1 $ for all $m_0 \in M^0$.  Let $d $ be the
index of $M^0$ in $M$.  Thus $m^d \in M^0$ for all $m \in M$ so that $\chi (m^d) = 1$ for all
$m \in M$.  Thus $\chi (m)$ is a $d^{th}$ root of unity and $|\chi (m)|_F = 1$ for all $m \in
M$.  Thus $<H_M(m), \chi > = log _q |\chi (m)|_F = 0$ for all $m \in M$. \qed

Recall that $Hom (X(M^0), \bfZ) \otimes _{\bfZ } \bfR \simeq 
Hom (X(A),\bfZ) \otimes _{\bfZ } \bfR = \fraka $ is the real Lie algebra of $A$,
$\fraka ^*= X(A) \otimes _{\bfZ } \bfR $ is its real dual, and
$\fraka ^*_{\bfC } = \fraka ^* \otimes _{\bfR } \bfC $ is its complex dual.  For each $\nu
\in \fraka ^*_{\bfC } $ we have a character $\chi ^0_{\nu }$ of $M^0$ defined by $\chi ^0_{\nu }
(m) = q^{<H_{M^0}(m), \nu >}, m \in M^0$.  
By Lemmas 5.1 and 5.2 the mapping $r$ above induces an
isomomorphism 
$$r_*: { X(M) \otimes _{\bfZ } \bfC \over X_0(M) \otimes _{\bfZ } \bfC}   \simeq X(A)
\otimes _{\bfZ } \bfC \simeq \fraka ^*_{\bfC }.$$
By Lemma 5.3, for each $m \in M$, $H_M(m)$ is an element of the complex dual of
${X(M) \otimes _{\bfZ } \bfC \over X_0(M) \otimes _{\bfZ } \bfC}$.  Thus for each $\nu \in
\fraka _{\bfC }^*$, we can define a character $\chi _{\nu } $ of $M$  by $$\chi _{\nu
} (m) = q^{<H_M(m), r_*^{-1}(\nu )>}, m \in M.$$ 

\proclaim Lemma 5.4.  For all $\nu \in \fraka _{\bfC }^*$, the restriction of $\chi _{\nu }$ to
$M^0$ is $\chi ^0_{\nu }$.  \endproclaim 

\proof  For $m_0 \in M^0, \chi \in X(M)$ we have $$<H_M(m_0), \chi > = log _q |\chi (m_0)|_F = 
log _q |\chi _0 (m_0)|_F = <H_{M_0}(m_0), \chi _0>$$ where $\chi _0$ denotes the restriction
of $\chi $ to $M^0$.  Since the isomorphism $r_*$ comes from the restriction map it is easy to
see that $<H_M(m_0), r_*^{-1}(\nu )> = <H_{M^0}(m_0), \nu >$ for all $m_0 \in M^0, \nu \in
\fraka ^*_{\bfC }$. \qed

As above, let $\sigma $ be an irreducible discrete series representation of $M$ and let $\sigma
_0$ be an irreducible constituent of the restriction of $\sigma $ to $M^0$.  
Let $V_0$ be the representation space for $\sigma _0$ and let 
$$W = \{ f:M \rightarrow V_0: f(mm_0) = \sigma _0(m_0)^{-1} f(m) {\rm \ for  \ all \ }m \in M,
m_0 \in M^0\}.$$
Then $M$ acts on $W$ by left translation and we will call this induced representation 
$(I_M, W)$.  
Let $V$ denote the representation space for $\sigma $ and fix a non-zero intertwining operator
$S:V \rightarrow W$ so that $S\sigma (m) = I_M(m) S$ for all $m \in M$.  Since $I_M$ is
unitary, we can also define a projection operator $P:W \rightarrow V$ so that
$PI_M(m) = \sigma (m)P$ for all $m \in M$ and $PSv = v$ for all $v \in V$.   We also
define representation spaces $$\thickmuskip=
.5\thickmuskip \calH _P(\sigma ) = \{ \phi\in C^\infty(G,V): \phi (xmn) = \delta
_P^{-{1\over 2}} (m) \sigma (m) ^{-1} \phi (x) {\rm \ for  \ all \ } x \in G, m \in M, n \in
N \};$$ 
$\calH _{P^0}(\sigma _0) = $
$$\{ \psi \in C^\infty(G, V_0): \psi (xm_0n) = \delta
_{P^0}^{-{1\over 2}} (m_0) \sigma _0(m_0) ^{-1} \psi  (x) {\rm \ for  \ all \ } x \in G, m_0
\in M^0, n \in N \};$$ 
$\calH _P(I_M) = $
$$\{ \psi\in C^\infty(G, W): \psi (xmn) = \delta
_P^{-{1\over 2}} (m) I_M(m) ^{-1} \psi (x) {\rm \ for  \ all \ } x \in G, m \in M, n \in N
\}.$$ In each case $G$ acts on the representation space by left translations and we call the
induced representations $I_P(\sigma ), I_{P^0}(\sigma _0)$, and $I_P(I_M)$ respectively.  They
are the representations $\Ind_P^G(\sigma ), \Ind _{P^0}^G(\sigma _0)$, and $\Ind _P^G(\Ind
_{M^0}^M(\sigma _0))$ respectively. 

The intertwining operators $S:V \rightarrow W$ and $P:W \rightarrow V$ induce intertwining
operators $S^*$ from $(I_P(\sigma ),\calH _P(\sigma ))$ to $(I_P(I_M),\calH _P(I_M))$ 
and $P^*$ from $(I_P(I_M),\calH _P(I_M))$ to $(I_P(\sigma ),\calH _P(\sigma ))$
given by
$$(S^*\phi )(x) = S\phi (x) {\rm \ for  \ all \ } \phi \in \calH _P(\sigma ), x \in G;$$
$$(P^*\psi )(x) = P\psi (x) {\rm \ for  \ all \ } \psi \in \calH _P(I_M ), x \in G.$$
There is also an equivalence $T$ between $(I_P(I_M),\calH _P(I_M))$ and
$( I_{P^0}(\sigma _0),\calH _{P^0}(\sigma _0) )$ given by
$$(T\psi )(x) = \psi (x)(1) {\rm \ for  \ all \ } \psi \in \calH _P(I_M), x \in G .$$
Its inverse is given by
$$(T^{-1}\psi ')(x)(m) = \delta _P^{1\over 2}(m) \psi ' (xm) {\rm \ for  \ all \ }
\psi ' \in \calH _{P^0}(\sigma _0) , x \in G, m \in M .$$

Recall for each $\nu \in \fraka ^*_{\bfC }$ we have defined characters $\chi _{\nu }^0$ of $M^0$
and $\chi _{\nu }$ of $M$ such that $\chi _{\nu }^0$ is the restriction of $\chi _{\nu }$ to
$M^0$.  We use these characters to define representations $\sigma (\nu ) = \sigma \otimes \chi
_{\nu }$ and $I_M(\nu ) = I_M \otimes \chi _{\nu }$ of $M$ and $\sigma _0 (\nu ) = \sigma _0
\otimes \chi _0(\nu )$ of $M^0$.  As above we use these to form induced representation spaces
$\calH _P(\sigma , \nu ) = \calH _P(\sigma (\nu )), \calH _P(I_M, \nu ) = \calH _P( I_M(\nu ))$,
and $\calH _{P^0}(\sigma , \nu ) = \calH _{P^0}(\sigma _0(\nu ))$. The intertwining operators
$S:V \rightarrow W$ and $P:W \rightarrow V$ also intertwine $\sigma (\nu )$ and $I_M(\nu )$ and
so as above define induced intertwining operators $S^*_{\nu }: \calH _P(\sigma ,\nu )
\rightarrow \calH _P( I_M,\nu )$ and $P^*_{\nu }:\calH _P( I_M,\nu ) \rightarrow
\calH _P(\sigma ,\nu )$.
   There are also equivalences $T_{\nu }:
\calH _P(I_M,\nu ) \rightarrow \calH _{P^0}(\sigma _0,\nu ) $ given by
$$(T_{\nu }\psi )(x) = \psi (x)(1) {\rm \ for  \ all \ } \psi \in \calH _P(I_M,\nu ), x \in G
.$$ The inverses are given by
$$(T^{-1}_{\nu }\psi ')(x)(m) = \delta _P^{1\over 2}(m) \chi _{\nu }(m) \psi ' (xm) {\rm \ for 
\ all \ } \psi ' \in \calH _{P^0}(\sigma _0, \nu ) , x \in G, m \in M .$$

Suppose $P_1=MN_1$ and $P_2=MN_2$ are two cuspidal parabolic subgroups of $G$ with Levi
component
$M$.  In \S 4 we defined a meromorphic family of intertwining operators 
$$J(P_2^0:P_1^0:\sigma _0:\nu ): \calH _{P_1^0}(\sigma _0, \nu ) \rightarrow
 \calH _{P_2^0}(\sigma _0, \nu ).$$
We can transfer these intertwining operators to the equivalent spaces $\calH _{P_i}(I_M,\nu
), i = 1,2$, by means of the equivalences $T_{\nu , P_i}$. Thus we define
$$J(P_2:P_1:I_M:\nu ):\calH _{P_1}(I_M, \nu ) \rightarrow \calH _{P_2}(I_M, \nu )$$ by
$$J(P_2:P_1:I_M:\nu ) = T_{\nu , P_2}^{-1} J(P_2^0:P_1^0:\sigma _0:\nu ) T_{\nu , P_1}.$$
We can also define
$$J(P_2:P_1:\sigma :\nu ):\calH _{P_1}(\sigma , \nu ) \rightarrow \calH _{P_2}(\sigma , \nu )$$
by $$J(P_2:P_1:\sigma :\nu ) = P^*_{\nu , P_2 } J(P_2:P_1:I_M:\nu ) S^*_{\nu ,
P_1}.$$

\proclaim Lemma 5.5.  Suppose that $\nu \in \fraka _{\bfC}^* (P_2^0:P_1^0)$ so that
$J(P_2^0:P_1^0:\sigma _0:\nu )$ is given by the convergent integral
$$J(P_2^0:P_1^0:\sigma _0:\nu )\psi ' (x) = \int _{\barN _1 \cap N_2} \psi '(x \barn) d \barn,
x \in G, \psi ' \in \calH _{P_1^0}(\sigma _0, \nu ).$$
Then $J(P_2:P_1:I_M:\nu )$ is given by the convergent integral
$$J(P_2:P_1:I_M:\nu ) \psi (x) =  \int _{\barN _1 \cap N_2} \psi (x \barn) d \barn,
x \in G, \psi \in \calH _{P_1}(I_M, \nu )$$ and
$J(P_2:P_1:\sigma :\nu )$ is also given by the convergent integral
$$J(P_2:P_1:\sigma:\nu ) \phi (x) =  \int _{\barN _1 \cap N_2} \phi (x \barn) d \barn,
x \in G, \phi \in \calH _{P_1}(\sigma , \nu ).$$
 \endproclaim

\proof  Using the definitions of the operators and the transformation property of the
representation space $\calH _{P_1}(I_M, \nu)$, we have for all $x \in G, m \in M, \psi \in 
 \calH _{P_1}(I_M, \nu)$, 
 $$J(P_2:P_1:I_M:\nu ) \psi (x) (m) = 
 T_{\nu , P_2}^{-1} J(P_2^0:P_1^0:\sigma _0:\nu ) T_{\nu , P_1} \psi (x) (m) $$
$$= \chi _{\nu }(m) \delta _{P_2}^{1\over 2}(m) 
(J(P_2^0:P_1^0:\sigma _0:\nu ) T_{\nu , P_1} \psi )(xm) $$
$$= \chi _{\nu }(m) \delta _{P_2}^{1\over 2}(m) 
\int _{\barN _1 \cap N_2} ( T_{\nu , P_1} \psi )(xm \barn )  d \barn $$
$$= \chi _{\nu }(m) \delta _{P_2}^{1\over 2}(m) 
\int _{\barN _1 \cap N_2}  \psi (xm \barn ) (1) d \barn $$
$$=  \delta _{P_2}^{1\over 2}(m) 
\int _{\barN _1 \cap N_2} \delta _{P_1}^{-{1\over 2}}(m) \psi (xm \barn m^{-1})(m) d \barn .$$
Now there is a homomorphism $\beta :M \rightarrow \bfR ^+$ so that for all $m \in M, f \in
C_c^{\infty }(\barN _1 \cap N_2) $, 
$$\int _{\barN _1 \cap N_2} f(m \barn m^{-1}) d \barn = 
\beta (m)\int _{\barN _1 \cap N_2} f( \barn ) d \barn .$$
For $m \in M^0$ we know that $$\beta (m) = 
\delta _{P_1}^{1\over 2}(m) 
\delta _{P_2}^{-{1\over 2}}(m).$$
Since $\beta \delta _{P_1}^{-{1\over 2}} 
\delta _{P_2}^{1\over 2}$ is a homomorphism from the finite group $M/M^0$ into $\bfR ^+$,
it must be identically one.  Hence for all $m \in M$ we have
$$\int _{\barN _1 \cap N_2} \psi (xm \barn m^{-1})(m) d \barn = 
\delta _{P_1}^{1\over 2}(m) 
\delta _{P_2}^{-{1\over 2}}(m)\int _{\barN _1 \cap N_2} \psi (x \barn )(m) d \barn .$$
Thus
$$J(P_2:P_1:I_M:\nu ) \psi (x) (m) = 
\int _{\barN _1 \cap N_2}\psi (x\barn )(m) d \barn .$$ 

Now for $x \in G, m \in M, \phi \in 
 \calH _{P_1}(\sigma, \nu)$, we have
 $$J(P_2:P_1:\sigma :\nu ) \phi (x)  = 
[ P^*_{\nu , P_2 } J(P_2:P_1:I_M:\nu ) S^*_{\nu ,
P_1}\phi ] (x) $$
$$ = P \cdot [ J(P_2:P_1:I_M:\nu ) S^*_{\nu ,
P_1}\phi ] (x) $$
$$ = P \cdot \int _{\barN _1 \cap N_2} [ S^*_{\nu , P_1}\phi ] (x\barn ) d \barn $$
$$ = P S \cdot \int _{\barN _1 \cap N_2} \phi (x\barn ) d \barn $$
$$ = \int _{\barN _1 \cap N_2} \phi (x\barn ) d \barn .$$
 \qed

 Let $r(P_2^0:P_1^0:\sigma _0:\nu )$ be the scalar normalizing factors used in \S 4 to define
the normalized intertwining operators 
$$R(P_2^0:P_1^0:\sigma _0) = r(P_2^0:P_1^0:\sigma _0: 0 )^{-1}J(P_2^0:P_1^0:\sigma _0:0 ).$$
The fact that $$ r(P_2^0:P_1^0:\sigma _0: \nu )^{-1}J(P_2^0:P_1^0:\sigma _0:\nu )$$ is
holomorphic and non-zero at $\nu = 0$ and that $T_{\nu , P_i}, i = 1,2$, are equivalences,
will imply that $$ r(P_2^0:P_1^0:\sigma _0: \nu )^{-1}J(P_2:P_1:I_M:\nu )$$ is also holomorphic
and non-zero at $\nu = 0$.  Thus we can define
$$R(P_2:P_1:I_M) = r(P_2^0:P_1^0:\sigma _0: 0 )^{-1}J(P_2:P_1:I_M:0 ).$$
We also define
$$R(P_2:P_1:\sigma ) = r(P_2^0:P_1^0:\sigma _0: 0 )^{-1}J(P_2:P_1:\sigma :0 ) =
P^*R(P_2:P_1:I_M)S^*.$$

\proclaim Lemma 5.6. Let $\phi \in \calH _{P_1}(\sigma )$.  Then $R(P_2:P_1:I_M)S^* \phi (x) \in
S(V)$ for all $x \in G$.  \endproclaim

{\hfuzz=13pt
\proof  For every $\nu $ we have an intertwining operator
$J(P_2:P_1:I_M: \nu )S^*_{\nu }: \calH _{P_1}(\sigma , \nu) \rightarrow \calH _{P_2}(I_M:\nu)$. 
In order to carry out arguments using the integral
formula and meromorphic extension of the intertwining operator we want a compact realization of
the representation.  Since we do not know if there is a maximal compact subgroup of $G$ which
meets every connected component, we proceed one coset at a time.  Let $G_M = G^0M$ and write
$G = \cup _{i=1}^k x_iG_M$.  Then $P \subset G_M$ for any parabolic subgroup $P$ of $G$ with
Levi component $M$.  Let $K^0$ be a good maximal compact subgroup of $G^0$ so that $G^0 =
K^0P^0$.  Thus $G_M = K^0P$.  Let 
$\calH _{K^0}(\sigma ) = $
$$\{ f_K\in C^\infty(K^0, V): f_K(kmn) = \sigma ^{-1}(m) f_K(k) {\rm \ for
\ all \  } m \in K^0 \cap M, n \in K^0 \cap N, k \in K^0 \}.$$
For any $\phi \in \calH _P(\sigma :\nu)$ and $1 \leq i \leq k$ we can
define  $$\phi _i(x) = \cases {\phi (x), & if $x \in x_iG_M$;\cr 0, & otherwise.\cr}$$
Then $\phi _i \in \calH _P(\sigma :\nu )$ for each $i$ and $\phi = \sum _{i=1}^k \phi _i$.
Thus every element of $\calH _P(\sigma :\nu )$ is a sum of elements supported on a single
coset of $G_M$ in $G$ and so it is enough to prove the lemma for $\phi \in \calH _{P_1}(\sigma )$
supported on a single coset of $G_M$ in $G$.  

Fix $1 \leq i \leq k$ and define 
$$F_i(\nu ): \calH _{P_1}(\sigma :\nu ) \rightarrow \calH _{K^0}(\sigma )$$ by
$F_i(\nu ) \phi (k) = \phi (x_ik)$ for all $k \in K^0$.  Define 
$$F_i^{-1}(\nu ): \calH _{K^0}(\sigma ) \rightarrow \calH _{P_1}(\sigma :\nu )$$ by
$$F_i^{-1}(\nu ) f_K(x) = \cases{ \delta _{P_1} ^{-{1\over 2}}(m)\sigma ^{-1}(m) \chi _{\nu
}^{-1}(m) f_K(k), & if $x = x_i kmn, k \in K^0, m \in M, n \in N$;\cr 0, & otherwise.\cr}$$
Then $F_i(\nu ) F_i^{-1}(\nu )f_K = f_K$ for all $f_K \in \calH _{K^0}(\sigma )$ and
$ F_i^{-1}(\nu )F_i(\nu ) \phi = \phi _i$ for all $\phi \in \calH _{P_1}(\sigma :\nu )$.}

Fix $1 \leq i \leq k$ and $\phi \in \calH _{P_1}(\sigma ) = \calH _{P_1}(\sigma : 0)$ such
that $\phi = \phi _i$ is supported on $x_iG_M$.  Let $f_K = F_i(0) \phi \in \calH _{K^0}(\sigma
)$ and for each $\nu $ define $\phi _i(\nu ) \in \calH _{P_1}(\sigma :\nu
)$ by $\phi _i(\nu :x) = F_i^{-1}(\nu )f_K (x), x \in G$.
 We have $\phi _i(0) =
 F_i^{-1}(0 )F_i(0) \phi = \phi _i = \phi $.

 Fix $w ^* \in S(V)^{\perp
}$ and define $\Phi (\nu :x) = <J(P_2:P_1:I_M:\nu )S^*_{\nu }\phi _i(\nu :x), w^*>$.  Then
$\nu \mapsto \Phi (\nu :x)$ is a meromorphic function of $\nu \in \fraka _{\bfC }^*$ for each
$x \in G$.   
If $\nu \in \fraka _{\bfC}^* (P_2^0:P_1^0)$, by Lemma
5.5 we have
$$J(P_2:P_1:I_M:\nu )S^*_{\nu }\phi _i(\nu :x) = \int _{\barN _1 \cap N_2} S \phi _i(\nu :
x\barn ) d \barn \in S(V).$$  Thus $\Phi (x:\nu ) = 0$ for all 
$\nu \in \fraka _{\bfC}^* (P_2^0:P_1^0)$ and
hence for all $\nu $.  Thus $$J(P_2:P_1:I_M:\nu )S^*_{\nu }\phi _i(\nu :x)  \in S(V)$$ for all
$\nu $ and so $$ r(P_2^0:P_1^0:\sigma _0: \nu )^{-1}J(P_2:P_1:I_M:\nu )
S^*_{\nu }\phi _i(\nu :x)  \in S(V)$$ for all
$\nu \in \fraka _{\bfC }^*, x \in G$.  In particular for $\nu = 0$ we have
$R(P_2:P_1:I_M)S^* \phi (x) \in
S(V)$ for all $x \in G$.  \qed

\proclaim Corollary 5.7.  Let $\phi \in \calH _{P_1}(\sigma )$.  Then $R(P_2:P_1:I_M)S^*\phi =
S^*R(P_2:P_1:\sigma )\phi $. \endproclaim

\proof  This follows from Lemma 5.6 since $SP$ is the identity on $S(V)$.  \qed

\proclaim Lemma 5.8.  Suppose $P_1,P_2$, and $P_3$ are cuspidal parabolic subgroups of $G$ with
Levi component $M$.  Then $$R(P_1:P_3:I_M) = R(P_1:P_2:I_M)R(P_2:P_3:I_M)$$ and
$$R(P_1:P_3:\sigma ) = R(P_1:P_2:\sigma )R(P_2:P_3:\sigma ).$$  \endproclaim

\proof  The statement for $I_M$ follows easily from Lemma 4.8 since the
intertwining operators $T_{\nu }$ are equivalences.  Now using Corollary 5.7 we have
$$R(P_1:P_3:\sigma )  = P^*R(P_1:P_3:I_M)S^* =  P^*R(P_1:P_2:I_M)R(P_2:P_3:I_M)S^* $$
$$P^*R(P_1:P_2:I_M)S^*R(P_2:P_3:\sigma ) = R(P_1:P_2:\sigma )R(P_2:P_3:\sigma ).$$  \qed 

\proclaim Lemma 5.9.  Suppose $P_1$ and $P_2$ are parabolic subgroups of $G$ with Levi
component
$M$.  Then $R(P_2:P_1:\sigma )$ is an equivalence from $\calH _{P_1}(\sigma )$ onto
$\calH _{P_2}(\sigma )$.  \endproclaim

\proof  The proof is the same as that of Lemma 4.9.  \qed

Now as above we define $N_G(\sigma _0) = \{ g \in N_G(A): \sigma _0^g \simeq \sigma
_0 \}$.  If $w \in W_G(\sigma _0) = N_G( \sigma _0)/M^0$, $\sigma _0$ can be extended to a
representation $\sigma _{0,w}$ of the group $M^0_w $ generated by $M^0$ and any representative $n_w$ for
$w$.  Define $T(n_w): W \rightarrow W$ by
$$T(n_w) f(m) = \sigma _{0,w}(n_w)f(n_w^{-1}mn_w), m \in
M.$$  It is easy to check that $$T(n_w)I_M(m) = I_M(n_wmn_w^{-1})T(n_w), m \in M.$$
 Next we define intertwining operators
$B_P(w): \calH _{w^{-1}Pw}(I_M) \rightarrow \calH _P(I_M)$ by $$B_P(w)\psi (g) = 
\alpha _{w^{-1}Pw}(n_w)^{1\over 2} T(n_w) \psi
(gn_w), g \in G .$$
Finally we define
self-intertwining operators $R(w, I_M): \calH _P(I_M) \rightarrow \calH _P(I_M)$ by
$$R(w,I_M) = B_P(w) R(w^{-1}Pw:P:I_M).$$ 
Recall that in \S 4 we defined intertwining operators $$A_P(w): \calH _{w^{-1}P^0w}(\sigma _0)
\rightarrow \calH _{P^0}(\sigma _0) {\rm \ and \ } R(w, \sigma _0):\calH _{P^0}(\sigma _0)
\rightarrow
\calH _{P^0}(\sigma _0).$$  It is easy to check that
$$B_P(w) = T_P^{-1} A_P(w) T_{w^{-1}Pw} {\rm \ and \ } R(w, I_M) = T_P^{-1} R(w, \sigma _0)
T_P$$ where $T_P: \calH _P(I_M) \rightarrow \calH _{P^0}(\sigma _0)$ is the equivalence defined
above.  In particular this implies that $B_P(w)$ and $R(w, I_M)$ are independent of the coset
representatives chosen.

Define $N_G(\sigma )  = \{ g \in N_G(M): \sigma ^g \simeq \sigma \}$.  If $w \in W_G(\sigma )
= N_G( \sigma )/M$, $\sigma $ can be extended to a representation of the group $M_w
$ generated by $M$ and any representative $n_w$ for $w$.  Denote such an extension by $\sigma
_w$ and define $A_P'(w): \calH _{w^{-1}Pw}(\sigma ) \rightarrow \calH _P(\sigma )$ by
$$(A_P'(w)\phi )(x) = \sigma _w(n_w) \alpha _{w^{-1}Pw}(n_w)^{1\over 2} \phi (xn_w).$$ 

\proclaim Lemma 5.10.  The intertwining operator $A_P'(w)$ is independent of the choice of coset
representative $n_w$ for $w \in W_G(\sigma )$.  For $w_1, w_2 \in W_G(\sigma )$ there is a
non-zero constant $c_P(w_1,w_2)$ so that
$$A_P'(w_1,w_2) = c_P(w_1,w_2) A_P'(w_1)A'_{w_1^{-1}Pw_1}(w_2).$$  \endproclaim

\proof  The proof is exactly the same as that of Lemma 4.11.  \qed

Finally, for $w \in W_G(\sigma )$, we define
$R(w, \sigma ): \calH _P(\sigma ) \rightarrow \calH _P(\sigma )$ by
$$R(w,\sigma ) = A_P'(w) R(w^{-1}Pw:P:\sigma ).$$ 
Note that for $u \in W_G(\sigma _0)$ we could also have defined an intertwining operator
$$R'(u, \sigma ): \calH _P(\sigma ) \rightarrow \calH _P(\sigma )$$ by
$$R'(u,\sigma ) = P^*R(u,I_M)S^*.$$ 

We want to relate these two definitions.  Let $$W_G(\sigma _0,\sigma ) = 
[N_G(\sigma _0) \cap N_G( \sigma )]/M^0 \subset W_G(\sigma _0).$$
Suppose that $x \in N_G(\sigma _0) \cap N_G( \sigma )$.  Then $x$ represents an element
$xM^0 \in W_G(\sigma _0, \sigma )$ and an element $ xM \in W_G(\sigma )$. 
Let $$p: W_G(\sigma _0, \sigma ) \rightarrow W_G(\sigma )$$ be given by  $p(xM^0)= xM, x \in 
 N_G(\sigma _0) \cap N_G( \sigma )$. 

\proclaim Lemma 5.11.  The mapping $p$ is surjective.  Its kernel is $W_M(\sigma _0) =
N_M(\sigma _0)/M^0$.  \endproclaim

\proof  Let $w \in W_G(\sigma )$ and
 let $x \in N_G(\sigma )$ be a representative for $w$.  Then $\sigma _0^x$ is contained
in the restriction of $\sigma ^x \simeq \sigma $ to $M^0$ so that there is $m \in M$ such that
$\sigma _0^x \simeq \sigma _0^m$.  Hence $w $ has a representative $xm^{-1} \in N_G(\sigma )
\cap N_G(\sigma _0)$.  Thus $p$ is surjective.  Clearly $p(xM^0) = M$ just in case $x \in
N_G(\sigma _0) \cap N_G(\sigma ) \cap M = N_M(\sigma _0)$.  \qed

\proclaim Lemma 5.12.  Suppose that $u \in W_G(\sigma _0)$ is in the complement of
$W_G(\sigma _0, \sigma )$.  Then
$$R'(u, \sigma )  = 0.$$  If $u \in W_G(\sigma _0, \sigma )$,
then there is a
complex constant $c$ so that $$R'(u,\sigma ) = c R(p(u) ,\sigma ).$$   
\endproclaim

\proof  Using Corollary 5.7, for any $u \in W_G(\sigma _0)$,
$$R'(u,\sigma ) =  P^*R(u,I_M)S^* = P^* B_P(u)R(u^{-1}Pu:P:I_M)S^* $$
$$ = P^* B_P(u)S^*R(u^{-1}Pu:P:\sigma ).$$
But for any $g \in G, \phi \in \calH _{u^{-1}Pu} (\sigma )$, if $x \in N_G(\sigma _0)$ is a
representative for $u$,
$$P^*B_P(u)S^*\phi (g)  = \alpha _{u^{-1}Pu}(x)^{1\over 2}(PT(x)S) \phi (gx) .$$  Since
$T(x)$ intertwines $I_M$ and $I_M^x$, we see that $PT(x)S$ intertwines $\sigma $ and $\sigma
^x \simeq \sigma ^u$.  Thus $PT(x)S = 0$ and hence $R'(u,\sigma ) =0$  unless $\sigma ^u
\simeq \sigma $.  

Suppose $u \in  W_G(\sigma _0, \sigma )$. Write $u = u_x, w = w_x = p(u), x \in N_G(\sigma _0)
\cap N_G( \sigma )$.  Then $PT(x)S$ and $\sigma _w(x)$ both intertwine $\sigma $ and $\sigma
^x$, and $\sigma _w(x) \not = 0$.  Thus there is a complex constant $c'$ so that
$PT(x)S = c'\sigma _w(x)$.  Thus for any $g \in G$,  
$$R'(u,\sigma ) \phi (g) = \alpha _{u^{-1}Pu}(x)^{1\over 2}(PT(x)S)[R(u^{-1}Pu:P:\sigma )
\phi ] (gx) $$ $$= c' \alpha _{u^{-1}Pu}(x)^{1\over 2}
 \sigma _w(x)[R(u^{-1}Pu:P:\sigma ) \phi ] (gx)$$ $$ =
c'A'_P(w)R(w^{-1}Pw:P:\sigma ) \phi (g)  
=c'R(w, \sigma )\phi (g).$$ \qed

\proclaim Lemma 5.13.  The $R(w, \sigma ), w \in  W_G( \sigma )$, span the algebra
$C(\sigma )$ of self-intertwining operators on $\calH _P(\sigma )$. \endproclaim

\proof  Let $R$ be a self-intertwining operator for $\calH _P(\sigma )$.  Then $S^*RP^*$ is a
self-intertwining operator for $\calH _P(I_M)$, hence in the span of the $R(u, I_M), u \in
W_G(\sigma _0)$. But then $R = P^*S^*RP^*S^*$ is in the span of the $P^*R(u,I_M)S^* = R'(u,
\sigma ), u \in W_G(\sigma _0)$.  But by Lemma 5.12, each $R'(u,\sigma )$ is either zero or a
multiple of one of the operators $R(w,\sigma ), w \in W_G(\sigma )$.  \qed

\proclaim Lemma 5.14.  Let $u \in W_G(\sigma _0)$ and 
suppose that $R(u, I_M)$ is scalar.  Then $u \in W_G(\sigma _0, \sigma
)$ and $R(p(u), \sigma ) $ is scalar.  \endproclaim

\proof  Suppose that there is a constant $s \in \bfC $ such $R(u, I_M)\psi = s\psi $ for all
$\psi \in \calH _P(I_M)$.  Since $R(u,I_M) \not = 0, s \not = 0$.  Now for all $\phi \in
\calH _P(\sigma )$,  $$R'(u, \sigma )\phi = P^* R(u, I_M) S^* \phi = s P^*S^*\phi = s \phi
.$$  Thus $R'(u,\sigma )$ is scalar and non-zero.  Thus by Lemma 5.12 we have
$u \in W_G(\sigma _0,\sigma
)$. 
Further, by Lemma 5.12, there is a constant
$c$ so that $R'(u, \sigma ) = cR(p(u),\sigma )$.  Since $R'(u,\sigma ) \not = 0, c \not =
0$.  Thus $R(p(u),\sigma ) = c^{-1}R'(u,\sigma )$ is scalar.  \qed 

\proclaim Lemma 5.15. There is a cocycle $\eta $ so that
$$R(w_1w_2, \sigma ) = \eta (w_1,w_2)R(w_1, \sigma ) R(w_2, \sigma )$$
for all $w_1, w_2 \in W_G(\sigma )$.
  \endproclaim

\proof  The proof is similar to Lemma 4.12. \qed

As in \S 4, if $W^0_{G^0}(\sigma _0)$ is the subgroup of elements $u \in W_{G^0}(\sigma _0)$
such that $R^0(u,\sigma _0)$ is scalar, then $W^0_{G^0}(\sigma _0) = W(\Phi _1)$ is generated by
reflections in a set  $\Phi _1$ of reduced roots of $(G,A)$.    
Let $\Phi ^+, \Phi _1^+$ be defined as in \S 4.  Since $M$ centralizes $A$,
$W_G(\sigma ) \subset N_G(A)/M$ acts on roots of $A$ and we can define
 $$R_{\sigma } = \{ w \in W_G(\sigma ): w\beta \in \Phi ^+ { \rm \ for \ all \  } \beta
\in \Phi _1^+ \}.$$  We want to prove the following.

\proclaim Theorem 5.16.   The $R(w,\sigma ), w \in R_{\sigma }$,
form a basis for the algebra of intertwining operators of $I_P(\sigma )$.  \endproclaim

In order to prove Theorem 5.16, we will first compute the dimension of
$C(\sigma )$ using our knowledge of the dimension of $C(\sigma _0)$.   We denote the
equivalence classes of $\Ind_P^G(\sigma )$ and $\Ind_{P^0}^G(\sigma
_0)$ by
$i_{G,M}(\sigma )$ and $i_{G,M^0}(\sigma _0)$ respectively.
 Let $X$ and $Y$ denote the groups of unitary characters of $G/G^0$ and $M/M^0$ respectively.
For $\chi \in X$, let $\chi _M \in Y$ denote the restriction of $\chi $ to $M$.  Define
$$X(\sigma ) = \{ \chi \in X: \chi _M \otimes \sigma \simeq \sigma \};$$
$$X_1(\sigma ) = \{ \chi \in X: \chi \otimes i_{G,M}(\sigma ) \simeq i_{G,M}(\sigma )\} ;$$
$$Y(\sigma ) = \{ \eta \in Y: \sigma \otimes \eta \simeq \sigma \}.$$
Let $s$ denote the multiplicity of $\sigma _0$ in the restriction of $\sigma $ to $M^0$.

\proclaim Lemma 5.17.  $$\dim C(\sigma _0) = 
 \dim C(\sigma ) s^2 [X/X(\sigma )] [X_1(\sigma )/X(\sigma )].$$  \endproclaim

\proof  Using Lemma 2.13,
$$i_{M,M^0}(\sigma _0) = s \sum _{\eta \in
Y/Y(\sigma )} \sigma \otimes \eta .$$
This implies that
$$i_{G,M}(i_{M,M^0}(\sigma _0)) =  s \sum _{\eta \in
Y/Y(\sigma )} i_{G,M}(\sigma \otimes \eta ).$$
Since both $G/G^0$ and $M/M^0$ are finite abelian, it is clear that the  map $\chi \mapsto \chi _M$ induces
an isomorphism between
$X/X(\sigma )$ and $Y/Y(\sigma )$.  
Thus we can rewrite
$$i_{G,M}(i_{M,M^0}(\sigma _0)) =    s \sum _{\chi \in
X/X(\sigma )} i_{G,M}(\sigma \otimes \chi _M).$$  But by Corollary 3.3 the induced representations
$i_{G,M}(\sigma \otimes \chi _M)$ are either disjoint or equal.  Further,
$$i_{G,M}(\sigma \otimes \chi _M) = i_{G,M}(\sigma ) \otimes \chi = i_{G,M}(\sigma )$$ just in
case $\chi \in X_1(\sigma )$.  Thus we have
$$i_{G,M}(i_{M,M^0}(\sigma _0)) =  s 
[X_1(\sigma )/X(\sigma )]\sum _{\chi \in
X/X_1(\sigma )} i_{G,M}(\sigma ) \otimes \chi $$ where the representations 
$ i_{G,M}(\sigma ) \otimes \chi $ are disjoint for $\chi \in X/X_1(\sigma )$. 
Thus
$$\dim C(\sigma _0) = \dim C(\sigma ) s^2 [X_1(\sigma )/X(\sigma )]^2
[X/X_1(\sigma )]$$ 
$$ =  \dim C(\sigma ) s^2 [X/X(\sigma )] [X_1(\sigma )/X(\sigma )].$$
\qed

\proclaim Lemma 5.18.  $$s^2 [X/X(\sigma )] [X_1(\sigma )/X(\sigma )] =
[W_G(\sigma _0)]/[W_G(\sigma )].$$ \endproclaim

\proof  First, using Lemma 2.13 we have 
$$s^2 [X/X(\sigma )] = s^2 [Y/Y(\sigma )] = [N_M(\sigma _0)/M^0] = [W_M(\sigma _0)].$$
We claim that $$[X_1(\sigma )/X(\sigma )]  = [W_G(\sigma _0)/W_G(\sigma _0, \sigma )].$$
This would establish the lemma since by Lemma 5.11 we have
$$[W_G(\sigma )] = [W_G(\sigma _0, \sigma )/W_M(\sigma _0)].$$

We will define a bijection between
$X_1(\sigma )/X(\sigma )$ and $$(N_G( \sigma _0) \cap N_G(\sigma )) \backslash
N_G(\sigma _0) \simeq W_G(\sigma _0, \sigma ) \backslash W_G(\sigma _0) .$$  
Let $\chi \in X_1(\sigma )$.  The equivalence class of $\sigma \otimes \chi _M$ depends only
on the coset $\barchi $ of $\chi $ in $X_1(\sigma )/X(\sigma )$.
Further, by definition of $X_1(\sigma )$, we have
 $$i_{G,M}(\sigma
\otimes \chi _M) = i_{G,M}(\sigma ) .$$  By Corollary 3.2 there is $x \in N_G(A)$ such
that
$\sigma \otimes \chi _M \simeq \sigma ^x$. Thus  $\sigma |_{M^0} \simeq  \sigma ^x |_{M^0}$. 
Thus $\sigma _0^x$ occurs in $\sigma |_{M^0}$ and so there is $m \in M$ such that $\sigma _0^x
\simeq \sigma _0^m$.  Then $y = xm^{-1} \in N_G(\sigma _0)$.  Although $y\in N_G(\sigma
_0)$ is not uniquely  determined by $\chi $, 
$$\sigma \otimes \chi _M \simeq \sigma ^{y_1} \simeq \sigma ^{y_2}$$ if and only if
$y_1y_2^{-1} \in N_G(\sigma _0) \cap N_G(\sigma )$.  Thus for each $\barchi \in
X_1(\sigma )/X(\sigma )$ there is a unique coset $\barx (\barchi ) = (N_G( \sigma _0) \cap
N_G(\sigma )) x
 $ in  $(N_G( \sigma _0) \cap N_G(\sigma ))
\backslash N_G(\sigma _0)$ such that $\chi _M \otimes \sigma \simeq \sigma ^x$.    Finally,
given $x \in N_G( \sigma _0), \sigma ^x$ is a
constituent of $i_{M,M^0}(\sigma _0^x) \simeq i_{M,M^0}(\sigma _0)$ so that there is $\eta \in
Y$ such that $\sigma ^x \simeq \sigma \otimes \eta $.  Now let $\chi \in X$ such that $\chi _M
= \eta $.  Then $\sigma \otimes \chi _M \simeq \sigma ^x$ so that $\barx = \barx (\barchi )$.  
\qed

Recall from \S 4 that $W_G(\sigma _0)$ is the semidirect product of subgroups $R_{\sigma _0}$
and
$W(\Phi _1)$ where $R(w, \sigma _0)$ is scalar for $w \in W(\Phi _1)$ and the $R(r,\sigma_0), r
\in R_{\sigma _0}$, give a basis for $C(\sigma _0)$. 

\proclaim Lemma 5.19.  $$\dim C(\sigma ) = [W_G(\sigma )]/[W(\Phi _1)].$$ \endproclaim

\proof  Combining Lemmas 5.17 and 5.18 we have
$$\dim C(\sigma _0) = \dim C(\sigma ) \cdot {[W_G(\sigma _0)]\over [W_G(\sigma )]}.$$
But from Lemma 4.15, $\dim C(\sigma _0) = [R_{\sigma _0}] = [W_G(\sigma _0)]/[W(\Phi
_1)].$   \qed

Since $W(\Phi _1) \subset N_{G^0}(A)/M^0$, it can be naturally embedded in $W_G(A) =
N_G(A)/M$.  
 
\proclaim Lemma 5.20.  $W_G(\sigma )$ is the semidirect product of $W(\Phi _1)$ and $R_{\sigma
}$.  For $w \in W_G(\sigma )$, $R(w,\sigma )$ is scalar if and only if $w \in W(\Phi _1)$. 
\endproclaim

\proof  If $x \in N_{G^0}(\sigma _0)$ represents an element of $W(\Phi _1)$, then by
Lemma 4.14 $R(u_x,I_M)$ is scalar.  Thus by Lemma 5.14, $w_x \in W_G(\sigma )$ and
$R(w_x,\sigma )$ is scalar.  
 Let $W^0_G(\sigma )$ denote the set of all $w \in W_G(\sigma )$
such that $R(w,\sigma )$ is scalar.  By the above $W(\Phi _1) \subset W^0_G(\sigma )$.
Using Lemmas 5.13 and 5.19 we see that
$$[W_G(\sigma )]/[W(\Phi _1)] = \dim C(\sigma ) \leq [W_G(\sigma )]/[W_G^0(\sigma )].$$  
Thus $W(\Phi _1) = W^0_G(\sigma )$.

Now as in the proof Lemma 4.14, $W(\Phi _1)$ is a normal subgroup of $W_G(\sigma )$ and
so $w\Phi _1 = \Phi _1$ for all $w \in W_G(\sigma )$.  This implies that
$$R_{\sigma } = \{ w \in W_G(\sigma ): w\Phi _1^+ = \Phi _1^+\}$$ which yields the
semidirect product decomposition.  \qed

\medbreak
\noindent {\bf Proof of Theorem 5.16.}  It follows from Lemmas 5.13,5.15, and 5.20 that the
$R(w,\sigma ), \hfil \break w
\in R_{\sigma }$, span the algebra $C(\sigma )$.  Further, by Lemmas 5.19 and 5.20, $$\dim
C(\sigma ) = [W_G(\sigma )]/W(\Phi _1)] = [R_{\sigma }].$$  \qed

Let $\eta $ be the cocycle of Lemma 5.15.  Exactly as in \S 4 we can fix a finite central
extension $$1 \rightarrow Z_{\sigma } \rightarrow \tilde R _{\sigma } \rightarrow R_{\sigma }
\rightarrow 1$$ over which $\eta $ splits, a character $\chi _{\sigma } $ of $Z_{\sigma }$, and
a representation
$\calR $ of 
$\tilde R _{\sigma } \times G$
on $\calH _P(\sigma )$.  Let  $\Pi (\tilde R _{\sigma
}, \chi _{\sigma })$ denote the set of irreducible representations of $\tilde R _{\sigma }$ with
$Z_{\sigma }$ central character $\chi _{\sigma }$, and let $\Pi _{\sigma
}(G)$ denote the set of irreducible constituents of $I_P(\sigma )$.  

\proclaim Theorem 5.21.  There is a bijection $\rho \mapsto \pi _{\rho }$ of $\Pi (\tilde R
_{\sigma }, \chi _{\sigma })$ onto $\Pi _{\sigma }(G)$ such that 
$$\calR = \oplus _{\rho \in \Pi (\tilde R _{\sigma }, \chi _{\sigma })}\ \ \  (\rho \chk \otimes
\pi _{\rho }).$$  \endproclaim

\proof  The proof is exactly the same as that of Theorem 4.17.  \qed

\section {Examples}

\noindent
For applications involving comparisons of representations between groups and
twisted trace formulas it is customary to use the following definition of parabolic
subgroup.  Let
$P^0$ be a parabolic subgroup of
$G^0$.  Then
$P = N_G(P^0)$ is a parabolic subgroup of $G$.  Thus, using this definition, parabolic
subgroups of
$G$ are in one to one correspondence with parabolic subgroups of $G^0$.   We will show that
the parabolic subgroups obtained using this definition are also parabolic subgroups using the
definition of \S 2.  However they are not cuspidal in general.  Indeed, recall from
Proposition 2.10 that if
$P^0$ is a parabolic subgroup of $G$, then the corresponding cuspidal parabolic subgroup is
the smallest parabolic subgroup of $G$ lying over $P^0$.  On the other hand, if
$P$ is any subgroup of $G$ with $P \cap G^0 = P^0$, then $P \subset N_G(P^0)$.  Thus
$N_G(P^0)$ will be the largest parabolic subgroup of $G$ lying over $P^0$.  We will give
examples to show that this class of parabolic subgroups, which we call N-parabolic
subgroups (N for normalizer), do not yield a nice theory of parabolically induced
representations of $G$.

\proclaim Lemma 6.1.  Let $P^0$ be a parabolic subgroup of $G^0$.  Then $P=N_G(P^0)$ is a
parabolic subgroup of $G$.  It is the largest parabolic subgroup lying over
$P^0$.  If $M^0$ is a Levi component for $P^0$, then $M = N_G(M^0) \cap P$ is a
Levi component for $P$.  \endproclaim

\proof  Let $P^0 = M^0N$ be a Levi decomposition of $P^0$ and define $M= N_G(M^0) \cap
P$.  Then $M \cap G^0 = M^0$ and $MN \subset P$.   Let $x \in P=N_G(P^0)$.  Then $xM^0x^{-1}$
is a Levi subgroup of $P^0$ and so there is $n \in N$ such that $xM^0x^{-1} = nM^0n^{-1}$. 
Now
$n^{-1}x \in N_G(M^0) \cap N_G(P^0) = M$ and so $x \in NM = MN$.  Thus $P = MN$.

Let $A$ be the split component of $M^0$.  Then
$M$ normalizes $A$ and we define a Weyl group $W = M/C_M(A)$ where $C_M(A)$ denotes the
centralizer of $A$ in $M$.  Since $M^0 \subset C_M(A)$, we know that $W$ is a finite
abelian group.  The split component of $M$ is $A' = \{ a \in A: xax^{-1} = a$ for all $x\in
M \} =
 \{ a\in A: wa = a$ for all $a \in W\}$. 
Let $M' = C_G(A')$.  If we can show that $A'$ is the split component of $M'$, then $A'$ is
a special vector subgroup.

 Let $\Phi ^+ = \Phi (P^0,A)$ denote the set of roots of $A$ in $P^0$, $\fraka
$ the real Lie algebra of $A$, and $\fraka ^+$ the positive chamber of $\fraka $ with respect
to
$\Phi ^+$.  Fix $w \in W$ and define $\fraka _w = \{ H \in \fraka : wH = H\}$. 
Since $M$ normalizes $P^0$, we have $w\Phi ^+ =\Phi ^+$ and
$w\fraka ^+ =\fraka ^+$.  Let $k$ be the order of $w$.
Then since $\fraka ^+$ is convex, for any $H \in \fraka ^+$ we have $H_w = H + wH + w^2H
+...+w^{k-1}H \in \fraka ^+$ with $wH_w = H_w$.  Thus 
$\fraka _w^+ = \fraka _w \cap \fraka ^+ \not = \emptyset $ and so for any $\alpha \in \Phi
^+$, the restriction of $\alpha $ to $\fraka _w $ is non-zero.  Since $W$ is a finite abelian
group, an easy induction argument shows in fact that the restriction of $\alpha $ to
$\fraka '$, the real Lie algebra of $A'$, is non-zero for every $\alpha \in \Phi ^+$.
Thus
$M' \cap G^0 = C_{G^0}(A') = C_{G^0}(A) = M^0$.

Let $A''$ be the split component of $M'$.  Thus $A' \subset A''$.  But since $M'\cap
G^0 = M^0$ and $M \subset M'$, we have $A''=
\{ a \in A: xax^{-1} = a$ for all $x\in M'\} \subset A'$ .   Thus $A'=A''$ is the split
component of $M'$.  This implies that $A'$ is a special vector subgroup and that
$M'=C_G(A')$ is a Levi subgroup of $G$.  But since 
the restriction of $\alpha $ to
$\fraka '$, is non-zero for every $\alpha \in \Phi ^+$, we can choose a set of positive
roots $(\Phi ')^+$ of $L(G)$ with respect to $L(A')$ by restricting the roots in $\Phi
^+$.  With this choice of positive roots, we obtain a parabolic subgroup $P'=M'N'$ of $G$
with $N'=N$.  Thus $M'$ normalizes $N$.  It also normalizes $M^0$ since $M' \cap G^0 =
M^0$.  Thus $M' \subset N_G(M^0) \cap N_G(P^0) = M$.   Now $M'=M$ so $P'=MN = P$.   \qed

Let $(P_0^0,A_0)$ be a minimal p-pair in $G^0$ and let $\Delta $ denote the set
of simple roots of $A_0$ in $P_0^0$.   Then as usual the standard parabolic
subgroups of $G^0$ are indexed by subsets $\Theta $ of $\Delta $.  Write
$(P_{\Theta }^0, A_{\Theta })$ for the standard parabolic pair of $G^0$
corresponding to
$\Theta \subset \Delta $ and write $P_{\Theta } = N_G(P_{\Theta }^0)$, the
N-parabolic subgroup of $G$ lying over $P_{\Theta }^0$.   Let $N_G(P_0^0,A_0)$
be the set of elements in
$G$ that normalize both
$A_0$ and $P_0^0$.  Clearly $N_G(P_0^0,A_0) \cap G^0 = P_0^0 \cap N_{G^0}(A_0) =
C_{G^0}(A_0) = M^0_0$.  Write $W_G(P_0^0,A_0) = N_G(P_0^0,A_0)/M^0_0$.  Then
$W_G(P_0^0,A_0)$ acts on $\Delta $ and for each $\Theta \subset \Delta $ we
write $$W(\Theta ) = \{ w \in W_G(P_0^0,A_0): w\Theta = \Theta \}.$$

\proclaim Lemma 6.2.  For all $\Theta \subset \Delta $, 
$$ P_{\Theta } = \cup _{w \in W(\Theta )}wP_{\Theta }^0.$$ \endproclaim

\proof  This follows from Lemma 3.8. \qed

Now that we have a simple method of computing the groups $P_{\Theta }$, we will
give examples to show the following unpleasant facts.

\medbreak
\noindent {\bf Fact 6.3.}  Let $P_1^0$ and $P_2^0$ be parabolic subgroups of
$G^0$ and let $P_i = N_G(P_i^0)$, i = 1,2.  Then $P_1^0 \subset P_2^0$ does not imply that
$P_1
\subset P_2$. 

\medbreak
\noindent {\bf Example 6.4.}  Let $G = O(2n), n \geq 2$.  Then $G^0 = SO(2n)$
and  the minimal
parabolic subgroup
$B^0$ of $G^0$ is the group of upper triangular matrices in $SO(2n)$ with
$A_0$ the subgroup of diagonal matrices.   The simple roots are $$\Delta = \{ e_1
- e_2, e_2 - e_3,...,e_{n-1}-e_n, e_{n-1} + e_n\}.$$  The Weyl group
$W_G(B^0,A_0)$ has order two and is generated by the sign change $c_n:e_n
\mapsto -e_n$ that interchanges $e_{n-1}-e_n$ and $ e_{n-1} + e_n$.  Thus a
subset $\Theta $ of $\Delta $ is stable under $c_n$ just in case neither or
both of $e_{n-1}\pm e_n$ belong to $\Theta $.  So for example if $\Theta = \{
e_{n-1} - e_n \}$, then $B^0 \subset P^0_{\Theta }$, but $B =
B^0 \cup c_nB^0$ is not contained in $P_{\Theta }=P^0_{\Theta }$.

We call $P_0$ a minimal N-parabolic subgroup of $G$ if given any
N-parabolic subgroup $P$ of $G$ there is $x \in G$ such that $P_0 \subset
xPx^{-1}$.

\medbreak
\noindent {\bf Fact 6.5.}  $G$ need not have a minimal N-parabolic subgroup.

Suppose $P_0$ is a minimal N-parabolic subgroup of $G$.  Then it is easy to
see that $P_0^0 = P_0 \cap G^0$ is a minimal parabolic subgroup of $G^0$.  Now by
Lemma 6.2, $P_0$ meets every connected component of $G$.  But as Example 6.4
shows, there are N-parabolic subgroups of $O(2n)$ which are contained in the
identity component $SO(2n)$.  Thus no conjugate could contain the minimal
N-parabolic subgroup.

If $P^0 = M^0N$ is a Levi decomposition for $P^0$, then
we obtained a Levi decomposition $P=MN$ for $P= N_G(P^0)$ by defining $M = N_G(M^0) \cap P$. 
Thus $M$ depends on both $M^0$ and $P$.  

\medbreak
\noindent {\bf Fact 6.6.}   
Suppose that $P_1 = M_1N_1$ and $P_2=M_2N_2$ are N-parabolic
subgroups of $G$ such that $M_1^0 = M_2^0$, ie. $P_1^0$ and $P_2^0$ have the
same Levi subgroup.  Then it need not be true that $M_1 = M_2$, or even that
$M_1$ and $M_2$ have the same number of connected components.

\medbreak
\noindent {\bf Example 6.7.}  Let $G = O(8)$ as in Example 6.4.  Let $P_1 =
P_{\Theta _1}$ where $\Theta _1 = \{e_3 - e_4 \}$ and let $P_2' = P_{\Theta
_2}$ where $\Theta _2 = \{ e_1 - e_2\}$.  Then $P_1 = P_1^0$ is connected and
$P'_2 = (P'_2)^0 \cup c_4 (P'_2)^0$ meets both components of $G$.  Let $w =
(13)(24)\in N_G(A_0)/A_0$ be the Weyl group element that permutes the pairs
$(e_1,e_3)$ and $(e_2,e_4)$ and define $P_2 = wP_2'w^{-1}$.  Then $wA_{\Theta
_2}w^{-1} = A_{\Theta _1}$ and so
$P^0_2 = w(P_2')^0w^{-1}$ and $P_1^0$ both have Levi component $M_1^0 =
C_{G^0}(A_{\Theta _1})$.  However $M_1 = M_1^0$ is connected and $M_2 = M_2^0
\cup c_2 M_2^0$ meets both components of $G$.  

In addition to structural problems, the class of N-parabolic subgroups does not yield a
nice theory of parabolic induction.  One of the basic
cornerstones of representation theory in the connected case is that every irreducible
admissible representation is contained in a representation which is parabolically induced
from a supercuspidal representation and every tempered representation is a
subrepresentation of a representation which is parabolically induced
from a discrete series representation.  But if supercuspidal and discrete series
representations are defined as in the connected case and in \S 2, then the Levi component
$M$ of a parabolic subgroup $P$ has no supercuspidal or discrete series representations
unless $P$ is cuspidal.   Thus we will not in general be able to obtain all irreducible
admissible or tempered representations of $G$ via induction from supercuspidal or
discrete series representations of N-parabolic subgroups.

\medbreak
\noindent {\bf Example 6.8.}  Define $G = O(2) = SO(2) \cup w SO(2)$ as in
Remark 2.2.  The only parabolic subgroup of $G^0 = SO(2) \simeq F^{\times
}$ is itself.  Thus the only N-parabolic subgroup
of $G$ is $N_G(G^0) = G$ itself which has split component $Z = \{ 1 \}$. But $G$ has no
representations with compactly supported matrix coefficients, hence no supercuspidal
representations.  It also has no representations with square-integrable matrix
coefficients, hence no discrete series representations.  

As can be seen in Example 6.8, the problem with using the standard definitions
for supercuspidal and discrete series representations with N-parabolic subgroups $P=MN$ is
that the split component of $M$ may be smaller than the split component of $M^0$. 
In order to guarantee the existence of enough
supercuspidal and discrete series representations we could define a
representation of $M$ to be supercuspidal (respectively discrete series) just in
case its restriction to $M^0$ is supercuspidal (respectively discrete
series).   Then it would be easy to prove as in Theorem 2.18 that
every irreducible admissible (respectively tempered)
representation of $G$ is contained in a representation which is induced
from a supercuspidal (respectively discrete series) representation of an N-parabolic
subgroup. 

Another basic property of parabolic induction in the connected case is the
following.  Suppose that
$P_1 = M_1N_1$ and $P_2=M_2N_2$ are parabolic subgroups and $\sigma _i, i =
1,2$, are irreducible representations of $M_i$ which are both either
supercuspidal or discrete series.  Then if the induced representations $\Ind
_{P_i}^G(\sigma _i)$ are not disjoint, then the pairs $(M_1,\sigma _1)$ and
$(M_2, \sigma _2)$ are conjugate.  Further, in the discrete series case, the
induced representations are equivalent.   These properties fail in the disconnected case
when the
$P_i$ are N-parabolic subgroups of $G$ and supercuspidal and discrete series
representations are defined as above.  

\medbreak
\noindent {\bf Example 6.9.}  Let $G = O(8)$ and define $P_1$ and $P_2$ as in
Example 6.7.  Recall in this case that $M_1 = M_1^0\simeq GL(2) \times GL(1)^2$
while
$M_2 = M_2^0
\cup c_2M_2^0 \simeq GL(2) \times GL(1) \times O(2)$ with $M_2^0 = GL(2)
\times GL(1) \times SO(2) = M_1^0$.  Let $\sigma _0 = \rho \otimes \chi _1
\otimes \chi _2$ be an irreducible unitary supercuspidal representation of $M_1^0
= M_2^0$ where $\rho $ is an irreducible unitary supercuspidal representation of
$GL(2)$, $\chi _1$ is a unitary character of $GL(1)$, and $\chi _2$ is a
non-trivial unitary character of $GL(1)$ with $\chi _2^2 = 1$.  Then $\sigma
_0^{c_2} = \rho \otimes \chi _1 \otimes \chi _2^{-1} =
\sigma _0$ so there is an irreducible representation $\sigma
_2$ of $M_2$ which extends $\sigma _0$.  Further, 
$$\Ind _{M_2^0}^{M_2}(\sigma _0) = \sigma _2 \oplus (\sigma _2 \otimes \eta )$$
where $\eta $ is the non-trivial character of $M_2/M_2^0$.   The
representations $\sigma _2$ and $\sigma _2 \otimes \eta $ of $M_2$ would both be
supercuspidal (and discrete series) since they both restrict to $\sigma _0$ on
$M_2^0$.   Now
$$\Ind _{P_2^0}^G(\sigma _0)
\simeq \Ind _{P_2}^G \Ind _{M_2^0}^{M_2}(\sigma _0) = \Ind _{P_2}^G (\sigma _2)
\oplus \Ind _{P_2}^G (\sigma _2 \otimes \eta ).$$

Since $P_1^0$ and $P_2^0$ are parabolic subgroups of $G^0$ with the same Levi
component $M_1^0 = M_2^0$ we have
$$\Ind _{P_1^0}^{G^0}(\sigma _0) \simeq \Ind _{P_2^0}^{G^0}(\sigma _0).$$
But $\sigma _1 =
\sigma _0$ is an irreducible supercuspidal (and discrete series) representation
of
$M_1 = M_1^0$ and 
$$\Ind _{P_1}^G(\sigma _1 ) = \Ind _{P_1^0}^G(\sigma _0) \simeq \Ind _{G^0}^G \Ind
_{P_1^0}^{G^0}(\sigma _0) $$ $$\simeq \Ind _{G^0}^G \Ind
_{P_2^0}^{G^0}(\sigma _0) \simeq  \Ind _{P_2^0}^G(\sigma _0) \simeq
\Ind _{P_2}^G (\sigma _2) \oplus \Ind _{P_2}^G (\sigma _2 \otimes \eta ).$$
Thus we have irreducible supercuspidal (and discrete series) representations
$\sigma _1$ of
$M_1$ and $\sigma _2$ of $M_2$ so that
$$\Ind _{P_1}^G(\sigma _1 ) \simeq \Ind _{P_2}^G (\sigma _2) \oplus \Ind _{P_2}^G
(\sigma _2 \otimes \eta ).$$  Clearly $M_1$ and $M_2$ cannot be conjugate in $G$
since $M_1$ is connected while $M_2$ is not.  Further, the representations $\Ind
_{P_i}^G(\sigma _i)$ have a nontrivial intertwining, but are not equivalent.

A final nice property of parabolic induction in the connected case is the
theory of $R$-groups.  If $P=MN$ is a parabolic subgroup and $\sigma $ is
an irreducible discrete series representation of $M$, then $R$ is a subgroup of
$W_G(\sigma )$, the group of Weyl group elements fixing $\sigma $, which
determines the reducibility of
$\Ind_P^G(\sigma )$.  Its most basic property is that the dimension of the
algebra of self-intertwining operators of $\Ind_P^G(\sigma )$ is equal to the
number of elements in $R$.
 The following example shows that there could not be such a simple
$R$-group theory in the disconnected case for N-parabolic subgroups.

\medbreak
\noindent {\bf Example 6.10.}  Let $G^0 = SL_2(F) \times SL_2(F) = \{ (x,y) : 
x,y \in SL_2(F) \}$.  Let $G = G^0 \cup \gamma G_0$ where $\gamma (x,y) \gamma
^{-1} = (y,x)$.  Let $B_1 = A_1N_1$ denote the usual Borel subgroup of
$G_1 = SL_2(F)$ where $A_1$ is the subgroup of diagonal elements and $N_1$ is the
subgroup of upper triangular matrices with ones on the diagonal.  Then $B^0 = B_1
\times B_1$ is a Borel subgroup of $G^0$ and $B = N_G(B^0) = B^0 \cup \gamma
B^0$. $B^0$ and $B$ have Levi decompositions $B^0 = M^0N$ and $B = MN$ where
$M^0 = A^0 = A_1 \times A_1, M = M^0 \cup \gamma M^0$, and $N = N_1\times N_1$.
We have Weyl groups $W(G^0,A^0) = W(G_1,A_1) \times W(G_1,A_1) \simeq Z_2
\times Z_2$ and $W(G, A^0) = W(G^0,A^0) \cup \gamma W(G^0,A^0) \simeq D_8$,
the dihedral group of order 8.

Let $\chi _1$ be a non-trivial character of $A_1 \simeq F^{\times }$ of order two
so that $\Ind ^{G_1}_{ B_1} (\chi _1) = \pi _1 \oplus \pi _s$ is a reducible
principal series representation.  Note that the $R$-group  here is $R_1 =
W(G_1,A_1) \simeq Z_2$ and we denote the irreducible constituents of the
induced representation by $\pi _1, \pi _s$ to indicate that they are
parameterized by the trivial and sign characters $\rho _1$ and $\rho _s$ of $R_1$
respectively.  Now let $\chi = \chi _1 \otimes \chi _1$.  Then 
$$\Ind ^{G^0}_{B^0}(\chi ) \simeq (\pi _1 \oplus \pi _s) \otimes (\pi _1 \oplus \pi
_s) = \pi _{11} \oplus \pi _{1s} \oplus \pi _{s1} \oplus \pi _{ss}$$ where for
$i,j \in \{ 1,s \}$, we write $\pi _{ij} = \pi _i \otimes \pi _j$.  Note that
$\pi _{ij} ^{\gamma } \simeq \pi _{ji}$, so that $\pi _{ij}^{\gamma } \simeq
\pi _{ij}$ if and only if $i = j$.   Let $\Pi _{ij}$ be an irreducible
representation of $G$ such that $\pi _{ij}$ is contained in the restriction of
$\Pi _{ij}$ to $G^0$.  Then if $i = j$ we have $$\Ind ^G_{G^0} (\pi _{ij}) = \Pi
_{ij} \oplus (\Pi _{ij} \otimes \eta )$$ where $\eta $ is the non-trivial
character of $G/G^0$ and $\Pi _{ij} \otimes \eta \not \simeq \Pi _{ij}$.  
  If $i \not = j$ we have $$\Ind ^G_{G^0}(\pi _{ij}) =
\Pi _{ij} \simeq \Pi _{ji} = \Ind ^G_{G^0}(\pi _{ji}).$$  In this case we also
have
$\Pi _{ij} \otimes \eta \simeq \Pi _{ij}$.  Thus we have
$$\Ind ^G_{B^0} (\chi ) = 
\Ind ^G_{G^0}(\pi _{11} \oplus \pi _{1s} \oplus \pi _{s1} \oplus \pi _{ss}) $$
$$= \Pi _{11} \oplus (\Pi _{11} \otimes \eta ) \oplus 
\Pi _{ss} \oplus (\Pi _{ss} \otimes \eta ) \oplus 2 \Pi _{1s}.$$

Note that the above decompositions are reflected in the $R$-groups as follows. 
First, the $R$-group for $\Ind ^{G^0}_{B^0}(\chi )$ is given by $R^0 = R_1 \times
R_1 = W(G^0,A^0) \simeq Z_2 \times Z_2$.  It has 4 characters $\rho ^0_{ij} =
\rho _i \otimes \rho _j, i,j \in \{1,s \}$ corresponding to the irreducible
representations $\pi _{ij}$.  The $R$-group for $\Ind ^G_{B^0}(\chi )$ is $R =
W(G, A_0) = R^0 \cup \gamma R^0$.  By abuse of notation we will denote the
non-trivial character of $R/R^0$ by the same letter $\eta $ as used above for
the non-trivial character of $G/G^0$ and below for the non-trivial character of
$M/M^0$.  The irreducible representations of $R$ are the characters $\rho _{ii}$
and $\rho _{ii} \otimes \eta $ where $$\Ind ^R_{R^0} (\rho ^0 _{ii}) = \rho _{ii}
\oplus (\rho _{ii} \otimes \eta )$$ and the two-dimensional irreducible
representation  $$\Ind ^R_{R^0} (\rho ^0 _{1s}) = \rho _{1s} \simeq \Ind ^R_{R^0}
(\rho ^0 _{s1}).$$   Thus we have the irreducible constituents of $\Ind ^G_{B^0}
(\chi )$ parameterized by the irreducible representations of $R$ and occurring
with multiplicities given by the degrees of the corresponding representations.  

Now we consider $\Ind ^M_{M^0} (\chi )$.  Since $M = M^0 \cup \gamma M^0$ and
$\chi ^{\gamma } = \chi $, we have $$\Ind ^M_{M^0} (\chi ) = \sigma \oplus
(\sigma
\otimes \eta )$$ where $\sigma , \sigma \otimes \eta $ are distinct
one-dimensional unitary representations of $M$ which restrict to $\chi $ on $M^0$.  Now,
using transitivity of induction and properties of tensor products, we can write
$$\Ind ^G_{ B^0} (\chi ) \simeq \Ind ^G_B \Ind ^M_{M^0} (\chi )
 $$ $$\simeq \Ind ^G_B (\sigma ) \oplus 
\Ind ^G_B (\sigma \otimes \eta ).$$   Now since 
$\Ind ^G_B (\sigma \otimes \eta ) \simeq \Ind ^G_B (\sigma ) \otimes \eta $,
we see that
$$\Ind ^G_B (\sigma ) \oplus (\Ind ^G_B (\sigma ) \otimes \eta ) \simeq $$
$$\Pi _{11} \oplus (\Pi _{11} \otimes \eta ) \oplus 
\Pi _{ss} \oplus (\Pi _{ss} \otimes \eta ) \oplus 2 \Pi _{1s}.$$
Thus we can assume that $\Pi _{11}, \Pi _{ss}$ were chosen so that
$$\Ind ^G_B (\sigma ) \simeq \Pi _{11} \oplus  
\Pi _{ss}  \oplus  \Pi _{1s}$$ and  
$$\Ind ^G_B (\sigma \otimes \eta ) \simeq (\Pi _{11} \otimes \eta ) \oplus  
(\Pi _{ss} \otimes \eta ) \oplus  \Pi _{1s}.$$

This example exhibits a number of unpleasant features.  First, we have
irreducible discrete series representations $\sigma _1 = \sigma $ and $\sigma _2 = \sigma
\otimes \eta $ of $B$ such that 
$\Ind ^G_B (\sigma _1)$ and $\Ind ^G_B (\sigma _2)$
have a non-trivial intertwining, but are not
equivalent. Second, $\Ind ^G_B (\sigma ) $ has 3 inequivalent irreducible
subrepresentations, each occuring with multiplicity one, so that the dimension of
its space of intertwining operators is 3.  There are no subgroups of any possible
Weyl groups here with order 3. 

\section*{References}

\bgroup\Resetstrings%
\def\Loccittest{}\def\Abbtest{ }\def\Capssmallcapstest{}\def\Edabbtest{ }\def\Edcapsmallcapstest{}\def\Underlinetest{ }%
\def\NoArev{0}\def\NoErev{0}\def\Acnt{1}\def\Ecnt{0}\def\acnt{0}\def\ecnt{0}%
\def\Ftest{ }\def\Fstr{1}%
\def\Atest{ }\def\Astr{J\Initper  Arthur}%
\def\Ttest{ }\def\Tstr{Unipotent automorphic representations: conjectures}%
\def\Jtest{ }\def\Jstr{Societ\'e Math\'e\-matique de France, Ast\'erisque}%
\def\Vtest{ }\def\Vstr{171-172}%
\def\Dtest{ }\def\Dstr{1989}%
\def\Ptest{ }\def\Pcnt{ }\def\Pstr{13--71}%
\Refformat\egroup%

\bgroup\Resetstrings%
\def\Loccittest{}\def\Abbtest{ }\def\Capssmallcapstest{}\def\Edabbtest{ }\def\Edcapsmallcapstest{}\def\Underlinetest{ }%
\def\NoArev{0}\def\NoErev{0}\def\Acnt{1}\def\Ecnt{0}\def\acnt{0}\def\ecnt{0}%
\def\Ftest{ }\def\Fstr{2}%
\def\Atest{ }\def\Astr{J\Initper  Arthur}%
\def\Ttest{ }\def\Tstr{On elliptic tempered characters}%
\def\Jtest{ }\def\Jstr{Acta Math.}%
\def\Vtest{ }\def\Vstr{171}%
\def\Dtest{ }\def\Dstr{1993}%
\def\Ptest{ }\def\Pcnt{ }\def\Pstr{73--138}%
\def\Astr{\Underlinemark}%
\Refformat\egroup%

\bgroup\Resetstrings%
\def\Loccittest{}\def\Abbtest{ }\def\Capssmallcapstest{}\def\Edabbtest{ }\def\Edcapsmallcapstest{}\def\Underlinetest{ }%
\def\NoArev{0}\def\NoErev{0}\def\Acnt{2}\def\Ecnt{0}\def\acnt{0}\def\ecnt{0}%
\def\Ftest{ }\def\Fstr{3}%
\def\Atest{ }\def\Astr{A\Initper  Borel%
  \Aand N\Initper  Wallach}%
\def\Ttest{ }\def\Tstr{Continuous Cohomology, Discrete Subgroups, and Representations of Reductive Groups}%
\def\Itest{ }\def\Istr{Princeton University Press}%
\def\Stest{ }\def\Sstr{Annals of Math. Studies}%
\def\Ntest{ }\def\Nstr{94}%
\def\Dtest{ }\def\Dstr{1980}%
\def\Ctest{ }\def\Cstr{Princeton, NJ}%
\Refformat\egroup%

\bgroup\Resetstrings%
\def\Loccittest{}\def\Abbtest{ }\def\Capssmallcapstest{}\def\Edabbtest{ }\def\Edcapsmallcapstest{}\def\Underlinetest{ }%
\def\NoArev{0}\def\NoErev{0}\def\Acnt{1}\def\Ecnt{0}\def\acnt{0}\def\ecnt{0}%
\def\Ftest{ }\def\Fstr{4}%
\def\Atest{ }\def\Astr{L\Initper  Clozel}%
\def\Ttest{ }\def\Tstr{Characters of non-connected reductive $p$-adic groups}%
\def\Jtest{ }\def\Jstr{Canad. J. Math.}%
\def\Vtest{ }\def\Vstr{39}%
\def\Ntest{ }\def\Nstr{1}%
\def\Dtest{ }\def\Dstr{1987}%
\def\Ptest{ }\def\Pcnt{ }\def\Pstr{149--167}%
\Refformat\egroup%

\bgroup\Resetstrings%
\def\Loccittest{}\def\Abbtest{ }\def\Capssmallcapstest{}\def\Edabbtest{ }\def\Edcapsmallcapstest{}\def\Underlinetest{ }%
\def\NoArev{0}\def\NoErev{0}\def\Acnt{2}\def\Ecnt{0}\def\acnt{0}\def\ecnt{0}%
\def\Ftest{ }\def\Fstr{5}%
\def\Atest{ }\def\Astr{S\Initper \Initgap S\Initper  Gelbart%
  \Aand A\Initper \Initgap W\Initper  Knapp}%
\def\Ttest{ }\def\Tstr{$L$--indistinguishability and $R$ groups for the special linear group}%
\def\Jtest{ }\def\Jstr{Adv. in Math.}%
\def\Vtest{ }\def\Vstr{43}%
\def\Dtest{ }\def\Dstr{1982}%
\def\Ptest{ }\def\Pcnt{ }\def\Pstr{101--121}%
\Refformat\egroup%

\bgroup\Resetstrings%
\def\Loccittest{}\def\Abbtest{ }\def\Capssmallcapstest{}\def\Edabbtest{ }\def\Edcapsmallcapstest{}\def\Underlinetest{ }%
\def\NoArev{0}\def\NoErev{0}\def\Acnt{1}\def\Ecnt{0}\def\acnt{0}\def\ecnt{0}%
\def\Ftest{ }\def\Fstr{6}%
\def\Atest{ }\def\Astr{D\Initper  Goldberg}%
\def\Ttest{ }\def\Tstr{Reducibility for non-connected $p$-adic groups with $G^\circ$ of prime index}%
\def\Jtest{ }\def\Jstr{Canad. J. Math.}%
\def\Vtest{ }\def\Vstr{47}%
\def\Dtest{ }\def\Dstr{1995}%
\def\Ptest{ }\def\Pcnt{ }\def\Pstr{344--363}%
\Refformat\egroup%

\bgroup\Resetstrings%
\def\Loccittest{}\def\Abbtest{ }\def\Capssmallcapstest{}\def\Edabbtest{ }\def\Edcapsmallcapstest{}\def\Underlinetest{ }%
\def\NoArev{0}\def\NoErev{0}\def\Acnt{1}\def\Ecnt{0}\def\acnt{0}\def\ecnt{0}%
\def\Ftest{ }\def\Fstr{7}%
\def\Atest{ }\def\Astr{Harish-Chandra}%
\def\Ttest{ }\def\Tstr{Harmonic analysis on reductive $p$-adic groups}%
\def\Dtest{ }\def\Dstr{1973}%
\def\Jtest{ }\def\Jstr{Proc. Sympos. Pure Math.}%
\def\Itest{ }\def\Istr{AMS}%
\def\Ctest{ }\def\Cstr{Providence, RI}%
\def\Vtest{ }\def\Vstr{26}%
\def\Ptest{ }\def\Pcnt{ }\def\Pstr{167--192}%
\Refformat\egroup%

\bgroup\Resetstrings%
\def\Loccittest{}\def\Abbtest{ }\def\Capssmallcapstest{}\def\Edabbtest{ }\def\Edcapsmallcapstest{}\def\Underlinetest{ }%
\def\NoArev{0}\def\NoErev{0}\def\Acnt{1}\def\Ecnt{0}\def\acnt{0}\def\ecnt{0}%
\def\Ftest{ }\def\Fstr{8}%
\def\Atest{ }\def\Astr{Harish-Chandra}%
\def\Ttest{ }\def\Tstr{Harmonic analysis on real reductive groups III. The Maass-Selberg relations and the Plancherel formula}%
\def\Jtest{ }\def\Jstr{Ann. of Math. (2)}%
\def\Vtest{ }\def\Vstr{104}%
\def\Dtest{ }\def\Dstr{1976}%
\def\Ptest{ }\def\Pcnt{ }\def\Pstr{117--201}%
\def\Astr{\Underlinemark}%
\Refformat\egroup%

\bgroup\Resetstrings%
\def\Loccittest{}\def\Abbtest{ }\def\Capssmallcapstest{}\def\Edabbtest{ }\def\Edcapsmallcapstest{}\def\Underlinetest{ }%
\def\NoArev{0}\def\NoErev{0}\def\Acnt{2}\def\Ecnt{0}\def\acnt{0}\def\ecnt{0}%
\def\Ftest{ }\def\Fstr{9}%
\def\Atest{ }\def\Astr{A\Initper \Initgap W\Initper  Knapp%
  \Aand G\Initper  Zuckerman}%
\def\Ttest{ }\def\Tstr{Classification of irreducible tempered representations of semisimple Lie groups}%
\def\Jtest{ }\def\Jstr{Proc. Nat. Acad. Sci. U.S.A.}%
\def\Vtest{ }\def\Vstr{73}%
\def\Ntest{ }\def\Nstr{7}%
\def\Ptest{ }\def\Pcnt{ }\def\Pstr{2178--2180}%
\def\Dtest{ }\def\Dstr{1976}%
\Refformat\egroup%

\bgroup\Resetstrings%
\def\Loccittest{}\def\Abbtest{ }\def\Capssmallcapstest{}\def\Edabbtest{ }\def\Edcapsmallcapstest{}\def\Underlinetest{ }%
\def\NoArev{0}\def\NoErev{0}\def\Acnt{1}\def\Ecnt{0}\def\acnt{0}\def\ecnt{0}%
\def\Ftest{ }\def\Fstr{10}%
\def\Atest{ }\def\Astr{R\Initper \Initgap P\Initper  Langlands}%
\def\Ttest{ }\def\Tstr{On the classification of irreducible representations of real algebraic groups}%
\def\Btest{ }\def\Bstr{Representation Theory and Harmonic Analysis on Semisimple Lie Groups}%
\def\Stest{ }\def\Sstr{Mathematical Surveys and Monograph}%
\def\Ntest{ }\def\Nstr{31}%
\def\Itest{ }\def\Istr{American Mathematical Society}%
\def\Ctest{ }\def\Cstr{Providence, RI}%
\def\Dtest{ }\def\Dstr{1989}%
\def\Ptest{ }\def\Pcnt{ }\def\Pstr{101--170}%
\Refformat\egroup%

\bgroup\Resetstrings%
\def\Loccittest{}\def\Abbtest{ }\def\Capssmallcapstest{}\def\Edabbtest{ }\def\Edcapsmallcapstest{}\def\Underlinetest{ }%
\def\NoArev{0}\def\NoErev{0}\def\Acnt{1}\def\Ecnt{0}\def\acnt{0}\def\ecnt{0}%
\def\Ftest{ }\def\Fstr{11}%
\def\Atest{ }\def\Astr{J\Initper \Initgap D\Initper  Rogawski}%
\def\Ttest{ }\def\Tstr{Trace Paley-Wiener theorem in the twisted case}%
\def\Jtest{ }\def\Jstr{Trans. Amer. Math. Soc.}%
\def\Vtest{ }\def\Vstr{309}%
\def\Dtest{ }\def\Dstr{1988}%
\def\Ptest{ }\def\Pcnt{ }\def\Pstr{215--229}%
\Refformat\egroup%

\bgroup\Resetstrings%
\def\Loccittest{}\def\Abbtest{ }\def\Capssmallcapstest{}\def\Edabbtest{ }\def\Edcapsmallcapstest{}\def\Underlinetest{ }%
\def\NoArev{0}\def\NoErev{0}\def\Acnt{1}\def\Ecnt{0}\def\acnt{0}\def\ecnt{0}%
\def\Ftest{ }\def\Fstr{12}%
\def\Atest{ }\def\Astr{D\Initper  Shelstad}%
\def\Ttest{ }\def\Tstr{$L$--indistinguishability for real groups}%
\def\Jtest{ }\def\Jstr{Math. Ann.}%
\def\Vtest{ }\def\Vstr{259}%
\def\Dtest{ }\def\Dstr{1982}%
\def\Ptest{ }\def\Pcnt{ }\def\Pstr{385--430}%
\Refformat\egroup%

\bgroup\Resetstrings%
\def\Loccittest{}\def\Abbtest{ }\def\Capssmallcapstest{}\def\Edabbtest{ }\def\Edcapsmallcapstest{}\def\Underlinetest{ }%
\def\NoArev{0}\def\NoErev{0}\def\Acnt{1}\def\Ecnt{0}\def\acnt{0}\def\ecnt{0}%
\def\Ftest{ }\def\Fstr{13}%
\def\Atest{ }\def\Astr{A\Initper \Initgap J\Initper  Silberger}%
\def\Ttest{ }\def\Tstr{Introduction to Harmonic Analysis on Reductive $p$-adic Groups}%
\def\Itest{ }\def\Istr{Princeton University Press}%
\def\Ctest{ }\def\Cstr{Princeton, NJ}%
\def\Stest{ }\def\Sstr{Mathematical Notes}%
\def\Ntest{ }\def\Nstr{23}%
\def\Dtest{ }\def\Dstr{1979}%
\Refformat\egroup%

\bgroup\Resetstrings%
\def\Loccittest{}\def\Abbtest{ }\def\Capssmallcapstest{}\def\Edabbtest{ }\def\Edcapsmallcapstest{}\def\Underlinetest{ }%
\def\NoArev{0}\def\NoErev{0}\def\Acnt{1}\def\Ecnt{0}\def\acnt{0}\def\ecnt{0}%
\def\Ftest{ }\def\Fstr{14}%
\def\Atest{ }\def\Astr{M\Initper  Tadic}%
\def\Ttest{ }\def\Tstr{Notes on representations of non-archimedean $SL(n)$}%
\def\Jtest{ }\def\Jstr{Pacific J. Math.}%
\def\Vtest{ }\def\Vstr{152}%
\def\Dtest{ }\def\Dstr{1992}%
\def\Ptest{ }\def\Pcnt{ }\def\Pstr{375--396}%
\Refformat\egroup%

\end{document}